\documentclass[12pt]{article} 

\usepackage[margin=2.5cm]{geometry}
\usepackage[english]{babel}
\usepackage{hyperref,xcolor}
\usepackage{amsmath}
\usepackage{amsthm}
\usepackage{amssymb}

\usepackage{mathtools}
\usepackage{ucs}

\usepackage{enumerate}
\usepackage[utf8]{inputenc}
\usepackage[T1]{fontenc}
\usepackage{ae,aecompl} 
\usepackage{tensor}

\usepackage{tikz-cd}

\def\defin{\textbf}

\DeclareMathOperator{\supp}{\mathrm{supp}}
\DeclareMathOperator{\Cost}{\mathrm{Cost}}

\newcommand{\1}{1}

\newcommand{\MAlg}{\mathrm{MAlg}}
\newcommand{\Aut}{\mathrm{Aut}}

\newcommand{\R}{\mathbb R}
\newcommand{\N}{\mathbb N}
\newcommand{\Q}{\mathbb Q}
\newcommand{\Z}{\mathbb Z}
\newcommand{\C}{\mathbb C}
\newcommand{\LL}{\mathrm L}

\newcommand{\id}{\mathrm{id}}  
\newcommand{\Stwo}{\mathfrak S_{2^\infty}}

\newcommand{\inv}{^{-1}}

\renewcommand{\leq}{\leqslant}
\renewcommand{\geq}{\geqslant}
\newcommand{\abs}[1]{\left\lvert #1\right\rvert}
\newcommand{\norm}[1]{\left\lVert #1\right\rVert}

\newcommand{\impl}{\Rightarrow}

\newcommand{\la}{\left\langle}
\newcommand{\ra}{\right\rangle}
\newcommand{\into}{\hookrightarrow}

\DeclareMathOperator{\Orb}{\mathrm{Orb}}

\newtheorem{thmi}{Theorem}

\newtheorem{thm}{Theorem}[section]
\newtheorem{theorem}[thm]{Theorem}
\newtheorem{cor}[thm]{Corollary}
\newtheorem{crl}[thm]{Corollary}
\newtheorem{corollary}[thm]{Corollary}
\newtheorem{lem}[thm]{Lemma}
\newtheorem{lemma}[thm]{Lemma}

\newtheorem{prop}[thm]{Proposition}
\newtheorem{proposition}[thm]{Proposition}

\theoremstyle{definition}

\newtheorem{qu}[thm]{Question}

\newtheorem*{ack}{Acknowledgements}

\newtheorem{df}[thm]{Definition}
\newtheorem{definition}[thm]{Definition}

\newtheorem{remark}[thm]{Remark}
\newtheorem{question}[thm]{Question}

\newtheorem{example}[thm]{Example}

\newcommand{\pmp}{p.m.p.\ }

\newcommand{\ph}{\varphi}
\newcommand{\eps}{\varepsilon}

\DeclareMathOperator{\Fix}{\mathrm{Fix}}
\DeclareMathOperator{\Stab}{\mathrm{Stab}}
\DeclareMathOperator{\Sub}{\mathrm{Sub}}
\DeclareMathOperator{\Iso}{\mathrm{Iso}}
\newcommand{\rel}{\mathcal R}

\title{Classification of non-free p.m.p.\ boolean  actions of ergodic full groups and applications}

\author{Alessandro Carderi, Alice Giraud and François Le Maître}

\begin{document}

\maketitle
\begin{abstract}

	We extend Dye's reconstruction theorem, which classifies isomorphisms between full groups,
	to a classification of homomorphisms between full groups.
	For full groups of ergodic p.m.p.\ equivalence relations, our result roughly says
	that such homomorphisms come only from  actions of the equivalence relation, or of one of its
	symmetric powers.
	This has several rigidity consequences for homomorphisms between full groups.
	Our main application is a characterization of property (T) for  ergodic p.m.p.\ equivalence
	relations purely in full group terms, without using their topology: an ergodic p.m.p.\ equivalence relation has (T) iff all non-free ergodic boolean actions of its full group are strongly ergodic.
\end{abstract}
\setcounter{tocdepth}{2}
\section{Introduction}

\subsection{General context}

\subsubsection{Dynamical systems and their orbits}

Given a bijection on a set,
the associated partition into orbits characterizes the bijection up to conjugacy. Indeed the conjugacy class of the bijection is completely determined by the number of orbits of size \(k\geq 1\) and the number of infinite orbits.
The picture gets much more intricate when the set \(X\) is endowed with additional structure,
for instance a topology or a measure. In this paper, we will work in a purely measurable context, but our results are inspired by results of Matte-Bon in the topological context \cite{mattebonRigidityPropertiesFull2018} as we will explain later on.

We begin with a standard probability space \((X,\mu)\), i.e.\ a
probability space isomorphic to the \([0,1]\) interval endowed with the Lebesgue measure.
The bijections that we consider are those which preserve this structure, called
\pmp (probability measure-preserving) bijections, and we will restrict to the \emph{aperiodic} ones, i.e.\ those all whose orbits are infinite.
A good example to have in mind is the following:
consider \(X=[0,1)\) endowed with the Lebesgue measure, let \(\alpha\in\R\setminus\Q\), then the irrational rotation \(T_\alpha\) is the \pmp bijection of \((X,\mu)\) defined by: for all \(x\in X\),
\[
	T_\alpha(x)=x+1\mod 1.
\]
Contrarily to bijections on sets without additional structure,
\pmp bijections are in general hard to classify up to conjugacy.
In the measurable context, the appropriate definition of conjugacy is the following: two \pmp bijections \(T_1\) and \(T_2\) of
a standard probability space \((X,\mu)\) are conjugate if there is another \pmp bijection \(S\) such that for \(\mu\)-\emph{almost all} \(x\in X\),
\[ST_1x=T_2Sx.\]
Note that this can be rewritten as the following equality of bijections up to measure zero: \(S\circ T_1\circ S\inv=T_2\), hence the name conjugacy.

In the context of irrational rotations, a good example of a conjugacy is any map \(S:[0,1)\to[0,1)\) such that for almost all \(x\in [0,1)\), \(S(x)=1-x\): such a map conjugates \(T_\alpha\) to \(T_{-\alpha}\).
This remark actually says it all for irrational rotations: the Halmos-von Neumann theorem implies
that \(T_\alpha\) is conjugate to \(T_\beta\) iff \(\alpha=\pm \beta \mod 1\) \cite{halmosOperatorMethodsClassical1942}.
Another fundamental class of transformations which are completely characterized up to conjugacy are \emph{Bernoulli shifts} \cite{kolmogorovNewMetricInvariant1958,ornsteinBernoulliShiftsSame1970},
which act by shifting coordinates on an infinite product of a fixed (possibly finite) probability space, 
the product being indexed by \(\Z\). 
In full generality, the problem of deciding whether two p.m.p.\ transformations are conjugate or not 
is intractable in a sense that can be made precise, see for instance
\cite{hjorthInvariantsMeasurePreserving2001,foremanconjugacyproblemergodic2011}.

A basic invariant of conjugacy is ergodicity, which says that every measurable set
which is a reunion of orbits must have measure \(0\) or \(1\).
By design, such a property is remembered by the partition of the space into orbits.
The irrational rotations are ergodic, and one can build aperiodic non-ergodic examples easily. For instance, noting
the \([0,1)\) interval is equal to the disjoint union \([0,1/2)\sqcup [1/2,1)\) one can define a \pmp bijection which acts as an irrational rotation on each of these intervals.
A finer invariant is then the \emph{ergodic decomposition} of a \pmp bijection, which roughly tells us the sizes of the ergodic parts, plus the size of a diffuse part where ergodic components have measure zero (think of \(T_\alpha\times \id_{[0,1]})\). Here is the fundamental result that started the study of dynamical systems through
their orbits in the measurable context.

\begin{thm}[{\cite[Thm.\ 5]{dyeGroupsMeasurePreserving1959}}]
	Two aperiodic p.m.p.\ bijections of a standard probability space can be conjugated so as to have the same orbits iff they have the same ergodic decomposition.
\end{thm}

When a \pmp bijection \(S\) establishes a bijection bewteen \(T_1\)-orbits and \(T_2\)-orbits,
let us say that \(S\) is an orbit equivalence.
Dye's theorem tells us that there is an orbit equivalence between two \pmp bijections
\(T_1\) and \(T_2\) iff they have the same ergodic decomposition.
In particular, \emph{any two ergodic \pmp bijections admit an orbit equivalence between them}.
\subsubsection{Full groups}

Dye's result was actually stated in a more abstract framework which is
central to us. In order to define it properly, let us first denote by
\(\Aut(X,\mu)\) the group of all \pmp bijections of a standard probability space \((X,\mu)\),
where we identity two such bijections if they coincide on a conull set.
Now Dye defines a subgroup \(\mathbb G\) of \(\Aut(X,\mu)\) to be \textbf{full}
if it is stable under \emph{countable cutting and pasting}: whenever \((A_n)\)
is a partition of \(X\), and \((T_n)\) is a sequence of elements of \(\mathbb G\) such that
\((T_n(A_n))\) is also a partition of \(X\), the \pmp bijection \(T:X\to X\) defined by
\[
	T(x)=T_n(x)\text{ for all }x\in A_n\text{ and all }n\in\N
\]
belongs to \(\mathbb G\).
The first example of full groups we want to give is that of the full group associated to
a p.m.p. bijection \(T\in\Aut(X,\mu)\), defined as
\[
	[T]\coloneqq\{U\in \Aut(X,\mu)\colon \forall x\in X, \, (x,U(x))\in \Orb_T(x)\},
\]
where \(\Orb_T(x)=\{T^n(x)\colon n\in\Z\}\) is the \(T\)-orbit of \(X\).
It is a good exercise to check that \([T]\) is the smallest full group containing \(T\),
so we might call such full groups \emph{singly generated}.

Moreover, it it another good exercise to check that \(S\in\Aut(X,\mu)\) is an orbit equivalence between
\(T_1,T_2\in\Aut(X,\mu)\) iff \(S[T_1]S\inv=[T_2]\).
So orbit equivalence is the study of full groups up to conjugacy, which
was Dye's original viewpoint. 
Note however that there are many more full groups, e.g.\ \(\Aut(X,\mu)\) itself is a full group. 
More importantly, replacing \pmp bijections by \emph{countable groups of \pmp bijections}, 
the situation can get much richer: there are many non conjugate ergodic countably generated full groups. 
We won't say more on this in this introductory section, 
but the reader can consult \cite{gaboriauOrbitEquivalenceMeasured2011} for a broad overview of orbit equivalence,
which has become a field in its own right.

Since conjugacy of full groups implies that they are isomorphic as abstract groups,
it is natural to wonder whether the converse holds.
As it turns out, Dye was also able to answer this question positively in the ergodic case
(see Sec.\ \ref{sec: ergodic full group} for the precise definition of
ergodicity for full groups), obtaining the following result,
also known as Dye's reconstruction theorem.

\begin{theorem}[{\cite[Thm.\ 2]{dyeGroupsMeasurePreserving1963}}]
	Let \(\mathbb G_1\) and \(\mathbb G_2\) be two ergodic full groups. Then for every group isomorphism \(\Psi: \mathbb G_1\to \mathbb G_2\), there is a unique \(S\in\Aut(X,\mu)\) such that
	for all \(T\in\mathbb G_1\),
	\[\Psi(T)=STS\inv.\]
\end{theorem}
The most general version of Dye's reconstruction theorem does not need ergodicity, but then \(S\) is only non-singular since \(\Psi\) might stretch ergodic components.

The main purpose of the present paper is to extend Dye's result by looking at group \emph{homomorphisms}
between ergodic full groups, rather than isomorphisms.
The fact that such a classification might be possible was hinted to us by Matte Bon's
remarkable work in the topological setup: he
was able to classify
actions of any  minimal topological full group over a Cantor space \(X\)
on another Cantor space \(Y\), which in particular allows him to describe homomorphism between topological full groups \cite{mattebonRigidityPropertiesFull2018}.
He proved that either the action has some free behaviour, or that it has to
come from the standard action of the full group on a symmetric power of \(X\).
We exhibit a similar classification for full
groups, and obtain applications which show how full groups
remember partitions into orbits as abstract groups.


\subsection{Results}
We will now state our main result after introducing some more terminology.

A \pmp \textbf{boolean action} of a full group \(\mathbb G\)
is simply a group homomorphism \(\mathbb G\to \Aut(Y,\nu)\), where \((Y,\nu)\) is another standard probability space.
Although full groups
do not admit non-trivial measure-preserving actions on standard probability space
because all their boolean actions are \emph{whirly} (see \cite{glasnerSpatialNonspatialActions2005}
and \cite[Sec.~5]{giordanoextremelyamenablegroups2007}),
they do admit boolean actions such as the one provided by their inclusion
into \(\Aut(X,\mu)\).
Other examples are provided by (symmetric) diagonal actions,
which we now define along with some useful auxiliary notation.

Given a standard Borel space \(X\) and \(n\in\N\),
we have a natural action of the symmetric group \(\mathfrak S_n\) on \(X^n\) by permuting
coordinates and the quotient \(X^{n}/\mathfrak S_n\) is still a standard Borel space
which we denote by \(X^{\odot n}\).
Moreover, if \((X,\mu)\) is a standard probability space,
we can endow the space \(X^{\odot n}\) with the pushforward measure via the
quotient map \(\pi_n:(X^n,\mu^{\otimes n})\to X^{\odot n}\), and we denote by
\(\mu^{\odot n}:=\pi_{n*}\mu^{\otimes n}\) this measure.
Given \(x_1,\ldots,x_n\in X\) we will use the notation
\([x_1,\ldots,x_n]:=\pi_n(x_1,\ldots,x_n)\in X^{\otimes n}\).
We have a group homomorphism \(\iota^{\odot n}:\Aut(X,\mu)\to \Aut(X^{\odot n},\mu^{\odot n})\)
defined by \(\iota^{\odot n}(T)[x_1,\ldots,x_n]=[Tx_1,\ldots,Tx_n]\), in particular
every full group admits countably many boolean actions obtained by restricting
\(\iota^{\odot n}\) to the full group in question.
It is not hard to see that these boolean actions are not
conjugate as soon as the full group is not the trivial group.

We can now state our measurable version of Matte Bon's theorem, which states roughly that every
boolean action
decomposes as a free part plus a countable union of boolean actions factoring
onto the above defined \(\iota^{\odot n}\) in a measure-preserving (and ``support-preserving'') manner.

\begin{thmi}\label{thm:main}
	Let \(\mathbb G\) be an ergodic full group over the standard probability space \((X,\mu)\) and let
	\(\rho\colon\mathbb G\to \Aut(Y,\nu)\) be a boolean action on another standard
	probability space \((Y,\nu)\). Then there is a unique measurable \(\rho(\mathbb G)\)-invariant
	partition\footnote{We allow partitions to contain empty sets throughout the whole paper.} \(\{Y_n\}_{n=0,\ldots,\infty}\) of \(Y\) such that
	the boolean actions  \(
		\rho_n\) induced on \((Y_n,\frac{\nu_{\restriction Y_n}}{\nu(Y_n)})\)
	are subject to the following conditions
	\begin{enumerate}
		\item \(\rho_0\) is the trivial action (it maps every element of \(\mathbb G\) to the identity map on \(Y_0\));
		\item for all \(n\in\mathbb N^*\), there is a (unique) measure preserving map
		      \(\pi_n:Y_n\rightarrow X^{\odot n}\) satisfying the
		      following two properties\begin{enumerate}
			      \item for all \(g\in\mathbb G\) and almost all \(y\in Y_n\)  we have
			            \(\pi_n(\rho_n(g)y)=\iota^{\odot n}(g)\pi_n(y)\);
			      \item for all \(g\in\mathbb G\) and almost all \(y\in Y_n\), if \(\rho_n(g)(y)\neq y\) then
			            \(\iota^{\odot n}(g)\pi_n(y)\neq\pi_n(y)\);
		      \end{enumerate}

		      \label{item2}
		\item \(\rho_\infty\) is free: for every \(g\in\mathbb G\setminus \{\id\}\) we have that
		      \(\{y\in Y_\infty\colon \rho_\infty(g)y=y\}\) is a null-set.
	\end{enumerate}
\end{thmi}

Let us briefly describe our approach to the above result.
An important distinction with the case of topological full groups is that
full groups of measure-preserving actions are uncountable. This difficulty is
counter-balanced by the presence of the so-called \emph{uniform metric}, which turns them
into topological groups endowed with a complete biinvariant metric, not separable in general.
Although our main result does not mention this group topology, it relies crucially on it,
as our first observation
towards it is an automatic continuity result which actually makes
our life much easier than in the topological context
(see Cor. \ref{cor: autocont action on malg}).

We also make a crucial
use of a classification of the invariant random subgroups of the group
of dyadic permutations due to Thomas and Tucker-Drob, while an important part
of Matte Bon's proof is to  first obtain an analogous classification
in the topological context, and to somehow extend it
to the full group using the topology. 
All in all, although our result is heavily
inspired by Matte Bon's, the proof is different and much easier in our setup.
We also use the following result
which might be of independent interest. Recall that a group topology is called \textbf{SIN} (Small Invariant Neighborhoods) when there is a basis of neighborhoods of the identity consisting of conjugacy invariant neighborhoods.

\begin{thmi}\label{thmi: two sin topologies}
	Let \(\mathbb G\) be an ergodic full group, let \(\tau\) be a SIN group topology on \(\mathbb G\).
	Then \(\tau\) is either the discrete or the uniform topology.
\end{thmi}

Note that the SIN hypothesis above is important, otherwise one can always endow \(\mathbb G\) with the
weak topology induced by \(\Aut(X,\mu)\). \\

Our main application of the  classification theorem for boolean action
is about full groups which are
\emph{separable} for the uniform topology, or equivalently which are countably generated.
Such full groups are better-known
as full groups of \emph{p.m.p.~equivalence relations} (with countable classes) and are thus closely related
to the well-developed field of orbit equivalence of p.m.p.\ actions.
This connection is particularly apparent in Dye's reconstruction theorem,
which can be restated as the fact that any abstract group isomorphism between full groups of
p.m.p.\ aperiodic equivalence relations must come from an isomorphism between
the equivalence relations themselves.
One thus expects group properties
of full groups to reflect properties of the equivalence relation. To our knowledge,
only three such examples were known so far for p.m.p. equivalence relations:
\begin{itemize}
	\item A non-trivial \pmp equivalence relation is ergodic if and only if its full group is simple \cite{eigenSimplicityFullGroup1981}.
	\item A non-trivial p.m.p.~equivalence relation is aperiodic if and only if
	      its full group has no index two normal subgroup \cite[Thm. 1.14]{lemaitrefullgroupsnonergodic2016}.
	\item An ergodic p.m.p.~equivalence relation is amenable if and only if every
	      action by homeomorphism of its full group on a compact metrizable space
	      admits a fixed point \cite[Cor.\ 1.4]{kittrellTopologicalPropertiesFull2010}.
\end{itemize}
We should also mention that the two last properties have natural continuous counterparts which are crucial in
establishing the above algebraic properties using automatic continuity. In particular, the last item above uses Giordano and Pestov's beautiful results on extreme amenability for full groups \cite{giordanoextremelyamenablegroups2007}.
Another topological
property which reflects properties of the p.m.p.\ equivalence relation is the topological
rank \cite{lemaitrenumbertopologicalgenerators2014},
but it has no natural purely algebraic counterpart.

Let us now explain how our main result allows us to obtain one more theorem of the above kind.
First note that every p.m.p.~action of an equivalence relation on a probability
space \((Y,\nu)\) induces a boolean action of its associated full group.
The main result of the present paper has the following consequence:
\emph{every non-free ergodic boolean action
	of the full group of a p.m.p.\ equivalence relation comes from a measure-preserving
	action of the equivalence relation itself or one of its symmetric powers} (see Corollary \ref{cor: main thm pmp to pmp}).

We can then obtain a dynamical characterization of property (T) for p.m.p.~equivalence relations,
purely in terms of their full group, thus adding a fourth item to the above list.

\begin{thmi}\label{thmi:charaTfullgroups}
	Let \(\mathcal R\) be a p.m.p. ergodic equivalence relation.
	Then \(\mathcal R\) has (T) if and only if all the non-free ergodic boolean actions
	of its full group on standard probability spaces are strongly ergodic.
\end{thmi}

Let us highlight the key ingredients of the above result,
apart from Theorem \ref{thm:main}.
The direct implication (shown in Section \ref{sec:directimpl}) relies on the fact that
every ergodic p.m.p.~action of a p.m.p.~equivalence relation with property (T) is strongly ergodic, and
that if \(\mathcal R\) has (T) then so do all its symmetric powers.
On the other hand, the converse relies on a natural generalization to p.m.p.~equivalence relations
of the Connes-Weiss result \cite{connesPropertyAsymptoticallyInvariant1980}
which states that if  all p.m.p.~ergodic actions of a countable group are strongly ergodic then the group has property (T) (see Theorem \ref{thm:connes-weiss}).
This result was also obtained independently by Grabowski, Jardón-Sánchez and Mellick in a work in preparation.

While the above results are to be expected,
we take the time to prove them in details in the last two sections. We also give a purely quantitative version of property (T) of
an ergodic p.m.p.~equivalence relations as a kind of strong property (T) of the full group, but with
respect to the restricted class of \emph{full} unitary representations (see~Definition~\ref{df: full urep} and Section \ref{sec:kazhdan pair}).
This quantitative version is required in order to show that if two ergodic \pmp equivalence relations  have property (T),
then so does their direct product (Proposition~\ref{prop: T carries over to products}). Since we are far from understanding all unitary representations of full groups, the following
question remains open:

\begin{qu}
	Let \(\mathcal R\) be a p.m.p.~ergodic equivalence relation with property (T). Does its full group \([\mathcal R]\)
	have strong property (T) ?
\end{qu}

Other applications of our main result provide restrictions on homomorphisms between full groups of \pmp equivalence relations, and we refer the reader to Section \ref{sec: applications} for details.
It seems natural to expect that such restrictions also hold in the more general setup of \emph{orbit full groups of \pmp actions of locally compact groups} (see \cite{carderiOrbitfullgroups2018}).
It would also be nice to remove the ergodicity hypothesis in our main result, but it is not known whether aperiodic full groups satisfy the automatic continuity property, and Theorem \ref{thmi: two sin topologies} is not true anymore (see \cite[Sec.~7]{lemaitrefullgroupsnonergodic2016} for a finer SIN topology which however does satisfy the automatic continuity property).

Another interesting question is to understand non \pmp boolean actions of ergodic full groups, and also non \pmp boolean actions of non \pmp ergodic full groups. This would probably require an understanding of quasi-invariant random subgroups of \(\mathfrak S_{2^\infty}\).

Finally, we point out that our main results implies Dye's reconstruction theorem for ergodic full groups in a rather straightforward manner (see Section \ref{sec:dyereconstruction}). However thanks to the work of Fremlin \cite[384D]{fremlinMeasureTheoryVol2002}, Dye's reconstruction theorem is now known to hold for a much wider class of groups of measure-preserving transformations, e.g.~the so-called \(\LL^1\) full groups considered in \cite{lemaitremeasurableanaloguesmall2018}. It would be interesting to know if Theorem \ref{thm:main} can be extended to such groups, possibly at the cost of requiring continuity of the boolean action.

\begin{ack}
	The results of the present paper first appeared as part of the second-named 
	author's PhD thesis \cite{giraudClassificationActionsErgodiques2018}. 
	We are grateful to the referees of the thesis Claire Anantharaman-Delaroche 
	and Robin Tucker-Drob for their comments which helped improve the exposition.
	In addition, we would like to warmly thank the referee for their
	very helpful remarks and suggestions.
\end{ack}

\tableofcontents

\section{Boolean actions}

In this work we will denote by \((X,\mu)\) or \((Y,\nu)\) atomless standard probability spaces, while  \((Z,\lambda)\) will usually be a standard probability space possibly with atoms. We denote by \(\MAlg(X,\mu)\)
the \textit{measure algebra} of \((X,\mu)\), that is the boolean algebra consisting of measurable
subsets of \(X\) up to measure zero. This algebra is complete, which means that
any family \((A_i)_{i\in I}\) of elements has a supremum and an infimum. We will
denote these by \(\bigvee_{i\in I}A_i\) and \(\bigwedge_{i\in I}A_i\) respectively. Whenever \(I\) is countable, they
coincide with set-theoretic union and intersection and we will also use the notations \(\bigcup_{i\in I}A_i\) and
\(\bigcap_{i\in I}A_i\).

For notational convenience, we will use the topology induced by the metric \(d_\mu\) on \(\MAlg(X,\mu)\)
given by \(d_\mu(A,B)=\mu(A\bigtriangleup B)\) (this metric is actually complete and separable, but we won't make direct use of this fact).
We will implicitly use the fact that all the boolean operations are uniformly continuous with respect to \(d_\mu\).
The topology induced by \(d_\mu\) on \(\MAlg(X,\mu)\) will be the only topology that we consider there, in particular we may omit to mention this topology explicitly in proofs.
Also note that measure subalgebras of \(\MAlg(X,\mu)\) are by definition
closed subalgebras, so that they are complete as well.

All the subsets of \(X\) that we will consider are seen as elements of
\(\MAlg(X,\mu)\).
In particular they are measurable and we will neglect what happens on sets of measure zero. We now give the pointwise definition of the group \(\Aut(X,\mu)\).

\begin{definition}
	A \emph{measure-preserving bijection} of \((X,\mu)\)
	is a Borel bijection \(T:X\to X\) such that \(\mu(T\inv(A))=\mu(A)\) for all Borel \(A\subseteq X\). Then \(\Aut(X,\mu)\) is the group of all measure-preserving bijections up to measure zero: we identify \(T_1\) and \(T_2\) as soon as
	\(\mu(\{x\in X\colon T_1(x)\neq T_2(x)\})=0\).
\end{definition}

Every measure-preserving bijection \(T\) defines naturally an automorphism of the measure algebra \(A\mapsto T(A)\), which gives rise to a faithful action by automorphisms of \(\Aut(X,\mu)\) on \(\MAlg(X,\mu)\).
Moreover, every measure-preserving automorphism of \(\MAlg(X,\mu)\) arises in this manner (see e.g.\ \cite[Sec.~1]{kechrisGlobalaspectsergodic2010}), so one can also view \(\Aut(X,\mu)\) as the group of all measure-preserving automorphisms of the measure algebra \(\MAlg(X,\mu)\).
While this second viewpoint is conceptually more satisfying, we will sometimes work with the first which is more concrete, especially when dealing when a single element of \(\Aut(X,\mu)\).

Given
\(T\in\Aut(X,\mu)\), we put \(\Fix(T):=\{x\in X:\ Tx=x\}\) and we denote by
\(\supp(T)\) its complement, called the \textbf{support} of \(T\). The latter can be defined purely in measure algebra terms as \(\supp T=\bigvee \{A\in\MAlg(X,\mu)\colon T(A)\cap A=\emptyset\}\).

The group \(\Aut(X,\mu)\) can be equipped with two metrizable group topologies which will play an important role in our arguments:
\begin{itemize}
	\item the \textbf{weak topology} for which a sequence \((T_n)_n\) converges to \(T\) if for every \(A\subset X\), we have that \((\mu(T_nA\Delta TA))_n\) converges to \(0\) as \(n\) tends to infinity;
	\item the \textbf{uniform topology} which is the topology induced by the \textit{uniform metric} \(d_u\) defined by
	      \[d_u(S,T)\coloneqq \mu(\lbrace x\in X\colon S(x)\neq T(x)\rbrace)=\mu(\supp(ST^{-1})).\]
\end{itemize}
The weak topology is a Polish topology, which means that it is separable and admits a compatible complete metric.
The uniform topology refines the weak topology and it is not separable. The metric \(d_u\) is complete and biinvariant: \(d_u(T_1ST_2,T_1TT_2)=d_u(S,T)\) for all \(S,T,T_1,T_2\in\Aut(X,\mu)\).

\begin{definition}
	A (p.m.p.) \textbf{boolean action} of a group \(G\) on a standard probability space \((X,\mu)\) is a homomorphism \(\rho\colon G\rightarrow \Aut(X,\mu)\). If \(G\) is a topological group, we say that a boolean action is \textbf{continuous} if it is continuous when we endow
	\(\Aut(X,\mu)\) with the weak topology.
\end{definition}


\begin{remark}
	Actions which are actually continuous when we endow \(\Aut(X,\mu)\) with the stronger uniform topology will play an important role in our work, but since the adjective \emph{uniformly continuous} could be confusing, we will not use it, and will explicitly state the topologies that we are considering instead.
\end{remark}

Whenever a group \(G\) acts on a standard probability space \((X,\mu)\) by measure-preserving bijections, we automatically obtain a boolean action on \((X,\mu)\).
Moreover, when the acting group is Polish and the action is Borel, the boolean action is continuous.
Given a Polish group \(G\), not every continuous boolean \(G\)-action on \((X,\mu)\) lifts to a Borel measure-preserving \(G\)-action in general, e.g.\ the natural (continuous) boolean \(\Aut(X,\mu)\)-action on \((X,\mu)\) does not have a so-called pointwise realization, see \cite{glasnerautomorphismgroupGaussian2005}.
However when \(G\) is locally compact second-countable, every continuous boolean \(G\)-action on \((X,\mu)\) lifts to a genuine Borel action on \((X,\mu)\) by measure-preserving bijections. We will often use without mention the much easier fact that when \(\Gamma\) is a countable discrete group, every boolean \(\Gamma\)-action lifts to a \(\Gamma\)-action on \((X,\mu)\) by measure-preserving bijections.

Given \(H\leq G\) and a boolean action \(\rho:G\to \Aut(X,\mu)\), we will denote the \textbf{restriction} of \(\rho\) to \(H\) by \(\rho_{\restriction H}:H\to\Aut(X,\mu)\), which is a boolean \(H\)-action. Moreover, if \(A\in\MAlg(X,\mu)\) is non zero and \(\rho\) (resp.\ \(\rho_{\restriction H}\))-invariant, we call \(\rho^A:G\to \Aut(A,\mu_A)\) (resp.~\(\rho^A_{\restriction H}:H\to \Aut(A,\mu_A))\) the \textbf{induced} boolean action on \((A,\mu_A)\), where \(\mu_A\) is the probability measure on \(A\) obtained by renormalizing \(\mu\) to a probability measure.

\begin{definition}
	A boolean action \(\rho\colon G\to \Aut(X,\mu)\) is \textbf{free} if the support of any element has full measure.
	It is \textbf{nowhere free} if its restriction to every invariant positive measure set is not free.
\end{definition}

Any boolean action \(\rho\colon G\to \Aut(X,\mu)\) splits as a free and a nowhere free boolean action. Indeed if we set \[X_{\text{free}}\coloneqq \bigwedge_{g\in G\setminus \{\id\}}\supp(\rho(g))\] then it is clear that \(X_{\text{free}}\) is \(\rho\)-invariant, the boolean \(G\)-action induced by \(\rho\) on \(X_{\text{free}}\) is free and the boolean \(G\)-action induced by \(\rho\) on \(X\setminus X_{\text{free}}\) is nowhere free.

\begin{definition}
	Given two boolean actions \(\rho_X: G\to \Aut(X,\mu)\) and \(\rho_Y: G\to \Aut(Y,\nu)\),
	we say that a measure-preserving map \(\pi: (Y,\nu)\to (X,\mu)\) is a \textbf{factor map}
	if for every \(g\in G\), and every \(A\in\MAlg(X,\mu)\) we have that \[\rho_Y(g)\pi^{-1}(A)=\pi^{-1}(\rho_X(g)A).\]
\end{definition}


From the pointwise viewpoint, a measure-preserving map \(\pi\) is a factor map if and only if for every \(g\in G\) and almost every \(y\in Y\),
\[
	\pi(\rho_Y(g)y)=\rho_X(g)\pi(y).\
\]
It follows that for any factor map \(\pi\) we have that \(\pi^{-1}(\supp(\rho_X(g))\subseteq\supp(\rho_Y(g))\). Factor maps where the reverse inclusion holds play a key role in our work, so let us give them a name.

\begin{definition}\label{df: support preserving}
	Given two boolean actions \(\rho_X: G\to \Aut(X,\mu)\) and \(\rho_Y: G\to \Aut(Y,\nu)\),
	we say that a factor map \(\pi: (Y,\nu)\to (X,\mu)\) is \textbf{support-preserving}
	when
	\(\pi^{-1}(\supp(\rho_X(g))=\supp(\rho_Y(g))\) for all \(g\in G\).
\end{definition}

\begin{proposition}\label{prop: easy support preserving}
	Given two boolean actions \(\rho_X: G\to \Aut(X,\mu)\) and \(\rho_Y: G\to \Aut(Y,\nu)\),
	Let  \(\pi: (Y,\nu)\to (X,\mu)\) be a factor map from
	\(\rho_Y: G\to \Aut(Y,\nu)\) onto \(\rho_X: G\to \Aut(X,\mu)\).
	Then \(\pi\) is support-preserving iff for all \(g\in G\), we have
	\(\mu(\supp(\rho_X(g))=\mu(\supp(\rho_Y(g))\) for all \(g\in G\).
\end{proposition}
\begin{proof}
	We already observed that \(\pi^{-1}(\supp(\rho_X(g))\subseteq\supp(\rho_Y(g))\),
	and since we are working at the measure algebra level,
	these two sets are equal iff they have the same measure.
	Since \(\pi\) is measure-preserving, the result follows.
\end{proof}

\begin{remark}\label{rmk: support pres is class bij}
	Let \(\Gamma\) be a countable group.
	Consider two boolean actions \(\rho_X: \Gamma\to \Aut(X,\mu)\) and \(\rho_Y: \Gamma\to \Aut(Y,\nu)\) viewed as \(\Gamma\)-actions by measure-preserving bijections, and let \(\pi:(Y,\nu)\to(X,\mu)\) be a factor map.
	Straightforward considerations on actions on sets show that the following are equivalent:
	\begin{itemize}
		\item \(\pi\) is support-preserving.
		\item for almost every \(y\in Y\), \(\pi\) induces a bijection from the orbit \(\rho_Y(\Gamma)y\) to the orbit \(\rho_X(\Gamma)\pi(y)\).
		\item for almost every \(y\in Y\), the \(\rho_Y\)-stabilizer of \(y\) is equal to the \(\rho_X\)-stabilizer of \(\pi(y)\).
	\end{itemize}
\end{remark}

We end this section by noting that being a support-preserving
factor map can be directly read at the level of the actions.

\begin{lemma}\label{lem: support preserving factor}
	Let \(\rho_X: G\to \Aut(X,\mu)\) and \(\rho_Y: G\to \Aut(Y,\nu)\) be two boolean \(G\)-actions,
	consider a factor map \(\pi: (Y,\nu)\to (X,\mu)\). Then \(\pi\) is support-preserving
	iff for all \(g\in G\), we have \(\mu(\supp \rho_X(g))=\mu(\supp\rho_Y(g))\).
\end{lemma}
\begin{proof}
	The direct implication is an immediate consequence of the fact that \(\pi\)
	is measure-preserving and the definition of being support preserving.
	To see the converse, recall that \(\pi^{-1}(\supp(\rho_X(g))\subseteq\supp(\rho_Y(g))\)
	so if they have the same measure they must be equal.
\end{proof}

\subsection{The symmetric-diagonal action}
\label{sec: sym diag}
Given a standard probability space \((X,\mu)\) and a natural number \(n\geq 1\),
we have a natural action of the symmetric group \(\mathfrak S_n\) on \(X^n\) by permuting
coordinates. We denote the quotient by \(X^{\odot n}\coloneqq X^{n}/\mathfrak S_n\) and the quotient map by \(\pi_n:X^n\rightarrow X^{\odot n}\). The standard Borel space \(X^{\odot n}\) can be endowed with the pushforward probability measure \(\mu^{\odot n}:=\pi_{n*}\mu^{\otimes n}\). The standard probability space \((X^{\odot n},\mu^{\odot n})\) is called  the \(n\)-\textit{symmetric power} of \((X,\mu)\). We will also use the notation \((X^{\odot 0},\mu^{\odot 0})\) for the trivial measure space, that is \(X^{\odot 0}\) is a point. Finally, given a sequence \((\alpha_i)_{0\leq i<\infty}\), we will set \[(X^{\odot (\alpha_i)_i},\mu^{\odot (\alpha_i)_i})\coloneqq \left(\bigsqcup_i X^{\odot i},\sum_i \alpha_i\mu^{\odot i}\right).\]

Given \(A_1,\ldots,A_n\subseteq X\) measurable  we will use the notation
\[[A_1,\ldots,A_n]:=\pi_n(A_1\times\ldots\times A_n)\subseteq X^{\odot n}.\]
Clearly \([A_1,\ldots,A_n]=[A_{\sigma (1)},\ldots,A_{\sigma(n)}]\) for every \(\sigma\in\mathfrak S_n\).

We have a boolean action \(\iota^{\odot n}:\Aut(X,\mu)\to \Aut(X^{\odot n},\mu^{\odot n})\)
uniquely defined by \[\iota^{\odot n}(T)[A_1,\ldots,A_n]=[TA_1,\ldots,TA_n],\]
From the pointwise viewpoint it is given by \(\iota^{\odot n}(T) \pi_n(x_1,\ldots,x_n)=\pi_n(Tx_1,\ldots,Tx_n)\).
For \(n=0\), we let \(\iota^{\odot 0}\) be the unique homomorphism from \(\Aut(X,\mu)\) to the trivial group \(\Aut(X^{\odot 0},\mu^{\odot 0})\), which is the trivial homomorphism. Note that for all \(n\geq 1\) we have and all \(T\in\Aut(X,\mu)\) we have
\begin{align*}
	\Fix(\iota^{\odot n}(T))=  & [\Fix(T),\ldots,\Fix(T)] \\
	\supp(\iota^{\odot n}(T))= & [\supp(T),X,\ldots,X],
\end{align*}
hence \(\mu^{\odot n}(\Fix(\iota^{\odot n}(T)))=\mu(\Fix(T))^n\). In particular, 	the map \(\pi_n:X^n\to X^{\odot n}\) is a support-preserving factor map from \(\iota^n\) to \(\iota^{\odot n}\), where \(\iota^n:\Aut(X,\mu)\to\Aut(X^n,\mu^{\otimes n})\) is the diagonal boolean action.

Gluing the boolean actions \(\iota^{\odot n}\) together,
for every sequence of non-negative reals \((\alpha_i)_{i<\infty}\) which sums to \(1\),
we obtain a boolean action \(\iota^{\odot (\alpha_i)_i}\) as the diagonal embedding
\begin{gather*}\iota^{\odot (\alpha_i)_i}:\Aut(X,\mu)\rightarrow \prod_{i\geq 0} \Aut(X^{\odot i},\alpha_i\mu^{\odot i})\leq \Aut(X^{\odot (\alpha_i)_i},\mu^{\odot (\alpha_i)_i}),\\
	\text{where }\iota^{\odot (\alpha_i)_i}(T)(A)\coloneqq \iota^{\odot i}(T)A\text{ if }A\subseteq X^{\odot i}.
\end{gather*}

\begin{df}\label{df:symdiagsum}
	Let \(G\leq \Aut(X,\mu)\).
	Given a sequence of non-negative reals \((\alpha_i)_{i<\infty}\)
	which sums to \(1\), the \textbf{symmetric diagonal sum} boolean action \(\iota^{\odot (\alpha_i)_i}\) is the boolean action \(\iota^{\odot (\alpha_i)_i}:G\rightarrow \Aut(X^{\odot (\alpha_i)_i},\mu^{\odot (\alpha_i)_i})\) obtained by restricting the above defined map to \(G\).
\end{df}

Observe that symmetric diagonal sums are uniform-to-uniform continuous. We will see that for ergodic full groups, they arise naturally from the non-free part of any boolean action.

\begin{remark}
	Our work provides a reasonable classification for nowhere free boolean actions of \(\Aut(X,\mu)\), and more generally of its full subgroups.
	Let us point out that there are however a priori many free boolean actions of \(\Aut(X,\mu)\) (and hence of full groups). Indeed, given any sequence of non negative integers \((m_n)_{n\geq 1}\), where infinitely many \(m_n\)'s are non zero, the natural boolean diagonal action \(\prod_n (\iota^{\odot n})^{m_n}\) of \(\Aut(X,\mu)\)  on the infinite direct product
	\[\prod_{n\geq 1}(X^{\odot n},\mu^{\odot n})^{m_n}\]
	is free. We wonder how many such actions there are up to conjugacy.
\end{remark}
\begin{remark}
	Let us also note that since \(\mu\) is atomless, the underlying set
	of a
	symmetric diagonal sum boolean action can be viewed as the set of
	all finite subsets of \(X\),
	where \(X^{\odot i}\) corresponds to all such subsets of cardinality \(i\)
	(see also the proof of Lemma \ref{lem: measure algebra generation for sym diag sum}).
\end{remark}

\subsection{Diagonally support-preserving factorizable actions}

We will be interested in boolean actions which factorize over a symmetric diagonal sum action in a support-preserving manner, which we will abbreviate as follows.

\begin{df}
	Let \(G\leq \Aut(X,\mu)\) be a subgroup. A boolean action \(\rho\colon G\rightarrow \Aut(Y,\nu)\) is \textbf{diagonally support-preserving factorizable} if there is \(Y_\infty\subseteq Y\) such that
	\begin{itemize}
		\item \(Y_\infty\) is \(G\)-invariant and the action induced by \(\rho\) on \(Y_\infty\) is a free boolean action,
		\item the action induced by \(\rho\) on \(Y\setminus Y_\infty\) factorizes over a symmetric diagonal sum in a support-preserving way.
	\end{itemize}
\end{df}

The point of the above definition is to give a compact way of stating Theorem \ref{thm:main}. Note that it does not cover the uniqueness in it, but we will now see  why uniqueness is automatic for the groups that we are considering.

Here is the general condition on the inclusion \(G\leq \Aut(X,\mu)\)
which will imply the uniqueness of such factorizations.
This condition was introduced by Dudko and Grigorchuk in
\cite[Def.~11]{dudkoDiagonalActionsBranch2018} (in a previous version
of this paper, we used a slight strengthening of this notion and later realized
that it was not needed).

\begin{df}
	A subgroup \(G\leq \Aut(X,\mu)\) is called  \textbf{absolutely non free} if for every measurable partition of \(X\) into \(n\) pieces \(A_1,...,A_n\) and every \(\epsilon>0\), there are \(T_1,...,T_n\in G\) with disjoint support such that for all \(i\in\{1,...,n\}\) we have
	\[\mu(A_i\bigtriangleup \supp T_i)<\epsilon.\]
\end{df}

This condition is met by \(\Aut(X,\mu)\), more generally ergodic full groups as well as some natural countable subgroups such as
the group of dyadic permutations, see Lemma \ref{lem:S2infty hanf}.
For examples coming from branch groups, see \cite{dudkoDiagonalActionsBranch2018}.

The next section will be devoted to the proof of the following proposition, which shows the uniqueness of the symmetric diagonal sum factor appearing in Theorem \ref{thm:main}.

\begin{prop}\label{prop:diagfactdec}
	Let \(G\leq \Aut(X,\mu)\) be an absolutely non free subgroup and consider a diagonally support-preserving factorizable boolean action \(\rho\colon G\rightarrow \Aut(Y,\nu)\). Then there is a unique partition up to null-sets \((Y_i)_{i=0,\ldots,\infty}\) of \(Y\) which is \(\rho\)-invariant and such that
	\begin{itemize}
		\item the action of \(G\) on \(Y_\infty\) is free;
		\item for every \(i<\infty\), the action of \(G\) on \(Y_i\) factorizes onto the action \(\iota^{\odot i}\) on \(X^{\odot i}\) in a support-preserving manner, and the corresponding support-preserving factor map is unique.
	\end{itemize}
\end{prop}

\subsection{Proof of the key proposition for uniqueness in Theorem \ref{thm:main}}

Before proving Proposition \ref{prop:diagfactdec}, we need a notion of independent interest.

Let \(G\leq \Aut(X,\mu)\) be a subgroup. The associated isotropy subalgebra \(\Iso_G(X,\mu)\)
is the measure subalgebra generated by the supports (or equivalently the set of fixed points)
of all the elements of \(G\), that is \(\Iso_G(X,\mu)\coloneqq\la\supp g: g\in G\ra\subseteq \MAlg(X,\mu)\).

\begin{df}[{\cite{vershikTotallyNonfreeActions2012}}]
	A subgroup \(G\leq \Aut(X,\mu)\) is said to be \textbf{totally non free} (TNF) if \(\Iso_G(X,\mu)=\MAlg(X,\mu)\).
\end{df}

Clearly absolutely non free actions are totally non free.
The following lemma is a simple observation which will yield Proposition \ref{prop:diagfactdec}.
It can be viewed as a generalization of Vershik's result that TNF actions
admit no nontrivial automorphisms.

\begin{lem}\label{lem:unique support preserving for TNF}
	Let \(\rho_1:G\to\Aut(Y_1,\nu_1)\) and \(\rho_2: G\to \Aut(Y_2,\nu_2)\) be two boolean actions of a group \(G\).
	Suppose \(\rho_2(G)\) is totally non free, and that
	\(\rho_1\) factors in a support-preserving manner onto \(\rho_2\) via \(\pi: (Y_1,\nu_1)\to (Y_2,\nu_2)\).
	Then \(\pi\) is the only factor map from \(\rho_1\) onto \(\rho_2\) and
	\(\rho_2\) is the unique totally non-free action that \(\rho_1\) factors onto
	in a support-preserving manner.
\end{lem}
\begin{proof}
	Suppose \(\pi\) is a support-preserving factor map from \(\rho_1\) onto \(\rho_2\),
	then by Proposition \ref{prop: easy support preserving}
	for all \(g\in G\) we have \(\nu_1(\supp \rho_1(g))=\nu_2(\supp \rho_2(g))\),
	and hence all factor maps from \(\rho_1\) onto \(\rho_2\) must be
	support-preserving.

	By total non freeness the measure algebra generated by \(\{\supp \rho_2(g): g\in G\}\) is equal to
	\(\MAlg(Y_2,\nu_2)\).
	So if \(\ph\) is any other factor map, \(\ph\inv\)
	is completely determined by the values it takes on
	\(\{ \supp \rho_2(g): g\in G\}\), and since it is support-preserving
	we must have
	\[
		\ph\inv(\supp\rho_2(g))=\supp\rho_1(g)=\pi\inv(\supp\rho_2(g)),
	\]
	so that \(\ph=\pi\) a.e., which proves the uniqueness of \(\pi\).
	Finally, to see that \(\rho_2\) is the only (up to isomorphism)
	totally non-free action that \(\rho_1\) factors onto in a support-preserving
	manner, simply observe that any such action
	is by definition isomorphic to the restriction of
	\(\rho_1\) to the subalgebra generated by all the \(\supp \rho_1(g)\)
	where \(g\in G\).
\end{proof}
\begin{remark}
	Note that once there is a support-preserving factor map
	\(\pi\) from \(\rho_1\) to \(\rho_2\), then
	all factor maps from \(\rho_1\) to \(\rho_2\)
	are support-preserving thanks to Lemma \ref{lem: support preserving factor},
	so \(\pi\) is actually the unique
	factor map from \(\rho_1\) to \(\rho_2\).
\end{remark}

We also need the following useful fact on symmetric diagonal sums,
which follows from the proof of  \cite[Thm.~14]{dudkoDiagonalActionsBranch2018}.
\begin{lemma}\label{lem: measure algebra generation for sym diag sum}
	Let \((\alpha_i)_{i<\infty}\) be a sequence of non-negative reals which sums to \(1\), consider the
	associated symmetric diagonal sum \((X^{\odot (\alpha_i)_i}, \mu^{\odot(\alpha_i)_i})\). Consider the map
	\[
		\begin{array}{crcl}
			\tilde\iota^{\odot (\alpha_i)_i}: & \MAlg(X,\mu) & \to     & \MAlg(X^{\odot (\alpha_i)_i}, \mu^{\odot(\alpha_i)_i}) \\
			                                  & A            & \mapsto & \bigsqcup_i \tilde \iota_i(A),
		\end{array}
	\]
	where for all \(i\geq 1\), we put \(\tilde\iota_i(A)=[A,X,\dots,X]\in\MAlg(X^{\odot i},\mu^{\odot i})\),
	and \(\tilde\iota_0(A)=\emptyset\).
	Then \(\tilde \iota^{\odot (\alpha_i)_i}\) is continuous and its image generates
	\(\MAlg(X^{\odot (\alpha_i)_i}, \mu^{\odot(\alpha_i)_i})\).
\end{lemma}
\begin{proof}
	The continuity of \(\tilde \iota^{\odot (\alpha_i)_i}\) is a direct consequence of
	the fact that each \(\tilde\iota_i\) is \(i\)-Lipschitz, which can be shown
	by noting that
	\[\tilde \iota_i(A)\bigtriangleup \tilde \iota_i(B)=
		[B\setminus A,X\setminus A,\dots,X\setminus A] \cup
		[A\setminus B, X\setminus B,\dots, X\setminus B],
	\]
	has \(\mu^{\odot i}\)-measure at most \(i\mu(A\bigtriangleup B)\).

	Let us now prove that the image of \(\tilde\iota^{\odot (\alpha_i)_i}\) generates the measure algebra.
	Observe that since \(\mu\) is atomless, for all \(i\),
	\(\mu^{\otimes i}\)-almost every \(i\)-tuple consists of distinct points.
	So we may view \(X^{\odot i}\) as the set of \(i\)-element subsets of \(X\),
	so that \(X^{\odot(\alpha_i)_i}\) is the set of finite subsets of \(X\).
	Then \(\tilde \iota^{\odot (\alpha_i)_i}(A)\) is simply the set of finite subsets which
	intersect \(A\).

	Without loss of generality, we assume \(X=[0,1]\), and denote by
	\(\mathcal I_\Q\) the countable set of all intervals with rational endpoints.
	Now suppose \(\bar x=\{x_1,...,x_k\}\) and \(\bar y=\{y_1,...,y_l\}\) are
	two distinct elements of \(X^{\odot(\alpha_i)_i}\), without loss
	of generality we may assume \(\bar x\setminus \bar y\) contains \(x_1\),
	and then if \(I_1\) is a rational interval which contains \(x_1\) but
	none of \(x_2,\dots, x_k,y_1,\dots,y_k\), we have that
	\(\bar x\in\tilde \iota^{\odot (\alpha_i)_i}(I)\) but \(\bar y\notin\tilde\iota^{\odot (\alpha_i)_i}(I)\).
	It follows that \(\{\tilde \iota^{\odot (\alpha_i)_i}(I)\colon I\in\mathcal I_{\Q}\}\)
	is a countable family of Borel sets which separates points,
	in particular it generates the measure algebra
	of \((X^{\odot(\alpha_i)_i}, \mu^{\odot(\alpha_i)_i})\) as wanted.
\end{proof}

\begin{prop}\label{prop:TNFprod}
	Let \(G\leq\Aut(X,\mu)\) be an absolutely non free subgroup.
	Let \((\alpha_i)_{0\leq i<\infty}\) be a sequence of non-negative reals such
	that \(\sum_{i}\alpha_i=1\), and let \(\iota^{\odot (\alpha_i)_i}\) be the associated
	symmetric diagonal sum boolean action.
	Then \(\iota^{\odot (\alpha_i)_i}(G)\leq\Aut(X^{\odot (\alpha_i)_i},\mu^{\odot (\alpha_i)_i})\)
	is not discrete for the uniform topology and totally non free.
\end{prop}
\begin{proof}
	Observe that by absolute non freeness
	we have a sequence of non-trivial
	elements \(T_n\in G\) such that \(\mu(\supp T_n)\to 0\). Since \(\sum_i\alpha_i=1\), we can exchange the limit and obtain
	\begin{align*}
		\lim_n\mu^{\odot (\alpha_i)_i}(\supp \iota^{\odot (\alpha_i)_i}(T_n))= &
		\lim_n\sum_i\alpha_i(1-\mu(\Fix(T_n)^i) =
		                                                                       & \sum_i\alpha_i\lim_n(1-\mu(\Fix(T_n)^i)=0
	\end{align*}
	so that \(\rho(G)\) is not discrete as claimed.

	Finally, observe that for all \(T\in G\) we have, using the notation
	of Lemma \ref{lem: measure algebra generation for sym diag sum},
	\[
		\supp \iota^{\odot (\alpha_i)_i}(T)=
		\tilde \iota^{\odot (\alpha_i)_i}(\supp T).
	\]
	By Lemma \ref{lem: measure algebra generation for sym diag sum},
	every element of \(\MAlg(X^{\odot(\alpha_i)_i}, \mu^{\odot(\alpha_i)_i})\)
	is a limit of boolean combinations of elements of
	\(\tilde \iota^{\odot (\alpha_i)_i}(\MAlg(X,\mu))\), and since
	\(\{\supp T\colon T\in G\}\) is dense in \(\MAlg(X,\mu)\) and
	\(\tilde \iota^{\odot (\alpha_i)_i}\) is continuous,
	every element of \(\MAlg(X^{\odot(\alpha_i)_i}, \mu^{\odot(\alpha_i)_i})\)
	is a limit of boolean combinations of supports of elements of
	\(\iota^{\odot(\alpha_i)_i}(G)\).
	This shows that \(\iota^{\odot(\alpha_i)_i}(G)\) is totally non free
	as wanted.
\end{proof}

The proof of Proposition \ref{prop:diagfactdec} can now be given.

\begin{proof}[Proof of Proposition \ref{prop:diagfactdec}]
	Set \(Y_\infty\coloneqq \bigwedge_{g\in G\setminus\{\id\}}\supp(\rho(T))\).
	Then the boolean action induced by \(\rho\) on \(Y_\infty\) is free.
	Let us now set \(Y_{<\infty}\coloneqq Y\setminus Y_\infty\).
	By assumption, there is a sequence \((\alpha_i)_{0\leq i<\infty}\) of non-negative reals such that \(\sum_{i\geq 0} \alpha_i=1-\mu(Y_\infty)\)
	and a support-preserving factor map \(\pi:(Y_{<\infty},\nu)\rightarrow (X^{\odot (\alpha_i)_i},\mu^{\odot (\alpha_i)_i})\).

	Let us first show that the sequence \((\alpha_i)_{0\leq i<\infty}\) is unique:
	indeed
	for all \(T\in G\), since the factor map \(\pi\)
	is support-preserving,
	the fixed set of \(\rho(T)\) has \(\nu\)-measure
	\(\sum_{i}\alpha_i\mu(\Fix T)^i\),
	which determines each coefficient \(\alpha_i\) uniquely
	because \(\{\mu(\Fix T)\colon T\in G\}\) is dense in \([0,1]\)
	by absolute non freeness.

	Now since the symmetric diagonal sum \(\iota^{\odot(\alpha_i)_i}\) is
	totally non free and \(\pi\) is support-preserving,
	\(\pi\) is the only factor map from \(\rho\) onto \(\iota^{\odot(\alpha_i)_i}\).
	In particular, we can set non-ambiguously \(Y_i\coloneqq \pi^{-1}(X^{\odot i})\),
	whose uniqueness follows from the previous discussion.
\end{proof}

\subsection{Diagonal factorization for the dyadic symmetric group}\label{sec: dyadic perm}

Fix the standard Cantor space \((X,\mu)\coloneqq(\{0,1\}^\N,\mathcal B(1/2)^{\otimes\N})\),
which we equip with the product of \(1/2\) Bernoulli measures \(\mathcal B(1/2)\coloneqq \frac 12(\delta_0+\delta_1)\).
We view the symmetric groups \(\mathfrak S_{2^n}:=\mathfrak S(\{0,1\}^n)\) as a
subgroup of \(\Aut(X,\mu)\) as follows:
for each \((x_k)_k\in X\) and each \(\sigma\in\mathfrak S_{2^n}\),
we set
\[(\sigma\cdot (x_k)_k)_i\coloneqq\left\lbrace\begin{array}{cc}(\sigma((x_k)_{k\leq n})))_i&\text{ if }i\leq n,\\ x_i&\text{ if }i>n;\end{array}\right.\]
that is \(\sigma\) acts as the concatenation of the finite permutation in \(\mathfrak S_{2^n}\) and the identity.
Note that with this identification in mind,
for each \(n\in\N\) we have \(\mathfrak S_{2^n}\leq\mathfrak S_{2^{n+1}}\leq \Aut(X,\mu)\).
We then define the \textbf{dyadic symmetric group} \(\mathfrak S_{2^\infty}\) as the union of all these finite groups:
\[\mathfrak S_{2^\infty}\coloneqq\bigcup_{n\in\N}\mathfrak S_{2^n}\leq \Aut(X,\mu).\]

\begin{lemma}\label{lem:S2infty hanf}
	The dyadic symmetric group \(\mathfrak S_{2^\infty}\leq \Aut(X,\mu)\) is absolutely non free.
\end{lemma}
\begin{proof}
	This follows from the fact that any cylinder subset of \(X\) is the support
	of some element of \(\mathfrak S_{2^\infty}\).
\end{proof}

\begin{prop}\label{prop:actionss2}
	Every measure-preserving action \(\rho:\mathfrak S_{2^\infty}\to\Aut(Y,\nu)\) is diagonally support-preserving factorizable.
\end{prop}

The above proposition will be deduced from a theorem of Thomas and Tucker-Drob
\cite{thomasInvariantRandomSubgroups2014}. This theorem is stated in the
context of invariant random subgroups (IRS) and we will need to introduce a
little bit of terminology before applying it.

For a countable group \(\Gamma\), we denote by \(\Sub(\Gamma)\subseteq \{0,1\}^\Gamma\)
the Borel set of subgroups of \(\Gamma\).
The group \(\Gamma\) acts on \(\Sub(\Gamma)\) by conjugation and an IRS of \(\Gamma\)
is by definition a conjugacy invariant probability measure on \(\Sub(\Gamma)\).
Note that the stabilizer of any \(\Lambda\in\Sub(\Gamma)\) for the \(\Gamma\)-action
by conjugacy is equal to its normalizer.

\begin{df}
	An IRS \(\zeta\in\mathrm{Prob}(\Sub(\Gamma))\) is \textbf{self-normalizing} if for \(\zeta\)-almost every
	\(\Lambda\in\Sub(\Gamma)\) we have that the normalizer of \(\Lambda\) (defined as \(N_\Gamma(\Lambda)\coloneqq \{\gamma\in\Gamma\colon \gamma \Lambda\gamma\inv =\Lambda\}\)) is equal to \(\Lambda\).
\end{df}



Let us fix a p.m.p.~action of the countable group \(\Gamma\) on \((X,\mu)\).
We define the measurable stabilizer map \(\Stab_\Gamma: X\to\Sub(\Gamma)\)
which maps \(x\in X\) to its stabilizer \(\Stab_\Gamma(x):=\{\gamma\in \Gamma: \gamma x=x\}\).
The pushforward measure \((\Stab_\Gamma)_*(\mu)\) is an IRS, which we call the \textbf{IRS of the action}.

\begin{lem}[{\cite[Prop.\ 4]{vershikTotallyNonfreeActions2012}}]\label{lem:TNF implies self-normalizing}
	The IRS of any totally non free action is self-normalizing.
\end{lem}

Also note that any self-normalizing
IRS (viewed as a p.m.p.~action) is totally non free (but as observed
by Vershik, there exists totally non free IRS which are not self-normalizing).
Here is another useful observation.

\begin{lem}\label{lem: TNF+SN has nice supp}
	Consider a p.m.p.~action  \(\alpha:\Gamma\rightarrow \Aut(X,\mu)\) of a countable group \(\Gamma\),
	and denote by \(\beta:\Gamma\rightarrow \Aut(\Sub(\Gamma),(\Stab_\Gamma)_*(\mu))\) the action
	associated to the IRS of \(\alpha\).
	Assume that the IRS \((\Stab_\Gamma)_*(\mu)\) is self-normalizing.
	Then the stabilizer map is support-preserving: for every \(\gamma\in \Gamma\)
	\[\Stab_\Gamma\inv(\supp \beta(\gamma))=\supp\alpha(\gamma).\]
	In particular the stabilizer map is the unique factor map from \(\alpha\)
	onto \(\beta\) by Lemma \ref{lem:unique support preserving for TNF}.
\end{lem}
\begin{proof}
	Since \(\beta\) is a factor of \(\alpha\) via the stabilizer map,
	we have \(\Stab\inv(\supp \beta(\gamma))\subseteq\supp \alpha(\gamma)\) for every \(\gamma\in \Gamma\).
	For the converse, the hypothesis that IRS is self-normalizing yields
	\[\supp \beta(\gamma)=\{\Lambda\leq \Gamma: \gamma\not\in N_\Gamma(\Lambda)\}=\{\Lambda\leq \Gamma:\ \gamma\not\in \Lambda\}.\]
	Using the definition of the stabilizer map, we finally have
	\[\mu(\Stab\inv(\supp \beta(\gamma)))=\mu(\{x\in X: \alpha(\gamma) x\neq x\})=\mu(\supp\alpha(\gamma))\]
	so the inclusion \(\Stab\inv(\supp \beta(\gamma))\subseteq\supp \alpha(\gamma)\) must be an equality and we are done.
\end{proof}

Let us go back to the proof of Proposition \ref{prop:actionss2}.
Assume that \((X,\mu):=(\{0,1\}^\N,\mathcal B(1/2)^{\otimes\N})\) and
consider \(\mathfrak S_{2^\infty}\leq \Aut(X,\mu)\).
For every \(i\geq 1\), we denote by \(\zeta_i\) the IRS associated to the symmetric diagonal action of \(\mathfrak S_{2^\infty}\)
on \((X^{\odot i}, \mu^{\odot i})\).

Note that for every non negative sequence \((\alpha_i)_{i\geq 0}\) such that \(\sum_{i\geq 0}\alpha_i=1\), the IRS \(\mathfrak S_{2^\infty}\curvearrowright(\Sub(\Stwo),\sum_{i\geq0} \alpha_i\zeta_i)\) is conjugate to the diagonal sum action \(\iota^{(\alpha_i)_i}\) of \(\mathfrak S_{2^\infty}\) as a consequence of Lemma \ref{lem:S2infty hanf} and  Proposition~\ref{prop:TNFprod} .

Finally we denote by \(\zeta_\infty\) the Dirac measure on the identity, which is the IRS associated to any free action and by \(\zeta_0\) the Dirac measure on \(\mathfrak S_{2^\infty}\), which is the IRS associated to any trivial action. We now state Thomas and Tucker-Drob's main result from \cite{thomasInvariantRandomSubgroups2014}.

\begin{thm}\label{thm:thomastucker}
	Every IRS of \(\mathfrak S_{2^\infty}\) can be written uniquely as an infinite convex combination
	of elements of the family \((\zeta_i)_{i=0,\ldots,\infty}\).
\end{thm}

Note that the theorem can be restated by saying that the IRS associated to any nowhere free action of \(\mathfrak S_{2^\infty}\) is a diagonal sum action. We are now ready to give the proof of Proposition \ref{prop:actionss2}.
\begin{proof}[Proof of Proposition \ref{prop:actionss2}]
	Consider a p.m.p.~action \(\rho\) of \(\mathfrak S_{2^\infty}\) on \((Y,\nu)\) and assume that the action is nowhere free. Denote by \(\zeta:=(\Stab_\Gamma)_*(\nu)\) the IRS associated to the action, which we view as a \pmp action using the same letter \(\zeta\). Denote by \((\alpha_i)_{0\leq i\leq \infty}\) the coefficients provided by Theorem \ref{thm:thomastucker}, since the action \(\rho\) is nowhere free \(\alpha_\infty=0\). As remarked above, it is a consequence  of Lemma \ref{lem:S2infty hanf} and  Proposition~\ref{prop:TNFprod}  that \(\zeta\) is a diagonal sum action, and \(\rho\) factorizes onto \(\zeta\) by construction. The only thing we have to show is that the factor map is support preserving.

	Proposition \ref{prop:TNFprod} and Lemma \ref{lem:S2infty hanf} imply that \(\zeta\) is totally non-free. Lemma \ref{lem:TNF implies self-normalizing} tells us that totally non free IRS are self normalizing, in particular so is \(\zeta\). Therefore we can use Lemma \ref{lem: TNF+SN has nice supp} to obtain that the factor map from \(\rho\) to \(\zeta\) is support-preserving, which finishes the proof.
\end{proof}

We end this section with another consequence of Theorem \ref{thm:thomastucker},
which was noted by the referee towards a simpler proof for Lemma \ref{lem:Aisinvariant} that they suggested.

\begin{cor}\label{cor: free part coincides with that of local subgroups}
	For every non empty clopen \(O\subseteq X\),
	denote by \((\mathfrak S_{2^\infty})_O\) the subgroup of all
	\(\sigma\in\mathfrak S_{2^\infty}\) such that \(\supp \sigma\subseteq O\).
	Then for any p.m.p.~\(\mathfrak S_{2^\infty}\)-action \(\rho:\mathfrak S_{2^\infty}\to \Aut(Y,\nu)\), the free
	part of \(\rho\) coincides with the free part of
	the restriction of \(\rho\) to \((\mathfrak S_{2^\infty})_O\).
\end{cor}
\begin{proof}
	By definition the free part of \(\rho\) is \(\Stab\inv (\{e\})\).
	Observe that for all \(i\geq 0\), the IRS \(\zeta_i\) has the property
	that almost every subgroup intersects \((\mathfrak S_{2^\infty})_O\)
	nontrivially,
	while \(\zeta_\infty\) is supported on \(\{e\}\).
	Applying Theorem \ref{thm:thomastucker}, we conclude that
	for almost all \(y\) not in the free part of \(\rho\), the stabilizer of \(y\)
	intersects \((\mathfrak S_{2^\infty})_O\)
	nontrivially. In other words, the free parts of \(\rho\) and of its restriction to \((\mathfrak S_{2^\infty})_O\) coincide as wanted.
\end{proof}

\section{Ergodic full groups}\label{sec: ergodic full group}

As in the last section, we will denote by \((X,\mu)\) a standard probability space without atoms.
Let us start by recalling Dye's definition of full groups.

\begin{df}\label{df: full group}
	A \textbf{full group} is a subgroup \(\mathbb G\) of \(\Aut(X,\mu)\) which is
	stable under \textit{cutting and pasting}, that is: whenever \((A_n)_n\)
	is a measurable partition of \(X\) and \((T_n)_n\) is a sequence of elements of \(\mathbb G\) such
	that \((T_n(A_n))_n\) is also a partition of \(X\), then the element
	\(T\in\Aut(X,\mu)\) defined by \(T(x)\coloneqq T_n(x)\) for \(x\in A_n\) belongs
	to \(\mathbb G\).
\end{df}

Given \(G\leq \Aut(X,\mu)\), the smallest full group containing \(G\) is denoted by \([G]\) and called the full group generated by \(G\);
it can be constructed by cutting and pasting the elements of \(G\).
A subgroup \(G\leq \Aut(X,\mu)\) is \textbf{ergodic} if the only \(G\)-invariant elements
of \(\MAlg(X,\mu)\) are \(X\) and \(\emptyset\). 
Recall that ergodicity can also be detected at the level of the \(\LL^2\) space:

\begin{proposition}[{see e.g.~\cite[Prop.~2.7]{kerrErgodicTheoryIndependence2016}}]\label{prop:ergo via L2} A subgroup \(G\leq\Aut(X,\mu)\) is ergodic iff the only functions \(f\in\LL^2(X,\mu)\) such that \(f\circ T=f\) for all \(T\in G\) are the constant functions.
\end{proposition}

The following theorem can be proven exactly as in \cite{fathiGroupeTransformationsQui1978} (see \cite[Thm.~3.11]{lemaitreGroupesPleinsPreservant2014} for the detailed proof in our setup).

\begin{theorem}\label{thm:ergodicvssimple}
	A full group is ergodic if and only if it is simple.
\end{theorem}

We will use the following important property of ergodic full groups.

\begin{prop}[{\cite[Lem. 3.2]{dyeGroupsMeasurePreserving1959}}] \label{prop: transitive on equal measure}
	Let \(\mathbb G\leq\Aut(X,\mu)\) be an ergodic full group. Then for every \(A,B\in\MAlg(X,\mu)\)
	such that \(\mu(A)=\mu(B)\), there is an involution \(T\in\mathbb G\) such that \(T(A)=B\) and \(\supp(T)=A\cup B\).
\end{prop}

Since \((X,\mu)\) is atomless, every measurable set can be written as the disjoint union of two measurable sets of equal measure. The previous proposition thus implies the following.

\begin{cor}\label{cor:every set is some support}
	Let \(\mathbb G\leq \Aut(X,\mu)\) be an ergodic full group and let \(A\in \MAlg(X,\mu)\). Then there
	is an involution \(U\in\mathbb G\) whose support is equal to \(A\). 
\end{cor}

Corollary \ref{cor:every set is some support} clearly implies that the inclusion \(\mathbb G\leq \Aut(X,\mu)\) is a highly absolutely non free and hence totally non free.

We will often use the following other well-known consequence of Proposition \ref{prop: transitive on equal measure}.

\begin{prop}\label{prop:conjugacy of invol}
	Let \(\mathbb G\leq \Aut(X,\mu)\) be an ergodic full group. Then two involutions in \(\mathbb G\)
	are conjugate if and only if their supports have the same measure.
\end{prop}
\begin{proof}
	Since all the elements of \(\mathbb G\) preserve the measure, if two involutions are conjugate
	then their supports must have the same measure.

	Conversely, let \(U,V\in\mathbb G\) be two involutions whose supports have the same measure. Since all standard Borel spaces are isomorphic, we can assume that \(X=[0,1]\) and in particular that we have a Borel linear order \(<\) on \(X\). Set \(A\coloneqq\{x\in X: U(x)>x\}\) and \(B\coloneqq\{x\in X: V(x)>x\}\),
	then it is straightforward to check that \(\supp U=A\sqcup U(A)\) and \(\supp V=B\sqcup V(B)\).
	Since \(U\) and \(V\) preserve the measure and \(\mu(\supp U)=\mu(\supp V)\), we deduce that \(\mu(A)=\mu(B)\).
	Applying Proposition \ref{prop: transitive on equal measure}, we find \(T\in \mathbb G\) such that \(T(A)=B\).
	Define a partial bijective measure-preserving map \(W\) by
	\[
		W(x)\coloneqq\left\{\begin{array}{cl}
			T(x)   & \text{if }x\in A      \\
			VTU(x) & \text{if } x\in U(A).
		\end{array}\right.
	\]
	The element \(W\) is only partially defined, but since its domain and range have equal measure, using Proposition \ref{prop: transitive on equal measure}
	one can easily extend \(W\) to an element \(\tilde W\) of \(\mathbb G\). Then \(\tilde WU\tilde W\inv=V\).
\end{proof}

The following result is key to all automatic continuity results on full groups.
It is due to Ryzhikov \cite{ryzhikovRepresentationTransformationsPreserving1985}
(see also the very neat proof in Miller's thesis \cite{millerFullGroupsClassification2004}).

\begin{thm}[Ryzhikov]\label{thm: 3invols}
	Let \(T\in\Aut(X,\mu)\). Then \(T\) can be written as the product of three involutions which belong
	to the full group generated by \(T\).
\end{thm}

Here is a well-known consequence (we include the proof since
we were not able to find that exact statement in the literature).

\begin{lem}\label{lem:generation in full groups}
	Let \(\mathbb G\leq \Aut(X,\mu)\) be an ergodic full group, let \(A_1,A_2\subseteq X\)
	such that \(\mu(A_1\cap A_2)>0\) and \(X=A_1\cup A_2\).
	Then \(\mathbb G=\la (\mathbb G)_{A_1},(\mathbb G)_{A_2}\ra\), where for all \(A\subseteq X\)
	we let \((\mathbb G)_A=\{T\in\mathbb G\colon \supp T\subseteq A\}\).
\end{lem}
\begin{proof}
	By Rhyzikov's theorem, it suffices to show that any involution \(U\in\mathbb G\)
	is a product of involutions in \(\mathbb G\) supported on either \(A_1\) or \(A_2\).
	So let \(U\in\mathbb G\) be an involution, and let \(\epsilon=\mu(A_1\cap A_2)>0\).
	Let \(B\) be a fundamental domain for the restriction of \(U\) to its support, so that
	\(\supp U=B\sqcup U(B)\).
	We partition \(B\) into finitely many subsets \(B_1,\dots B_k\) of measure \(\leq \epsilon\) such that each \(B_i\) is contained in either \(A_1\) or \(A_2\).

	Then we can write \(U=U_1\cdots U_k\) where
	\[U_i(x)=
		\left\{
		\begin{array}{cl}
			U(x) & \text{if } x\in B_i\sqcup U(B_i) \\
			x    & \text{otherwise.}
		\end{array}
		\right.
	\]
	So without loss of generality we may as well assume that \(U\) admits a fundamental domain \(B\)
	of measure \(\leq \epsilon\)
	contained in either \(A_1\) or \(A_2\).
	By symmetry, we may as well assume \(B\subseteq A_1\).
	Now, let \(C=\{x\in B \colon U(x)\notin A_1\}\), let \(D=U(C)\subseteq A_2\).
	By ergodicity and Proposition \ref{prop: transitive on equal measure}
	there is an involution \(V\in (\mathbb G)_{A_2}\) such that \(V(U(C))\subseteq A_1\cap A_2\),
	so that \(VU\in (\mathbb G)_{A_1}\), which finishes the proof.
\end{proof}

Finally, we note that the notion of support-preserving factor map as defined in Definition \ref{df: support preserving} behaves well with respect to full groups. In what follows, given a subgroup \(G\leq \Aut(X,\mu)\), we denote by \(\iota_G\) its inclusion in \(\Aut(X,\mu)\), which is a boolean action.
\begin{prop} \label{prop: support preserving extend to full group}
	Let \(G\leq \Aut(X,\mu)\) and \(\rho_Y: G\to \Aut(Y,\nu)\) be boolean action, suppose \(\pi:(Y,\nu)\to(X,\mu)\) is a support-preserving factor map from \(\rho_Y\) to \(\iota_G\).
	Then there is a unique extension \([\rho_Y]:[G]\to\Aut(Y,\nu)\)
	of \(\rho_Y\)  such that \(\pi\) is still support-preserving from \([\rho_Y]\) to \(\iota_{[G]}\).
\end{prop}
\begin{proof}
	Since \(\pi\) is support-preserving, if \(U_1,U_2\in G\) satisfy that
	\(U_{1\restriction A}=U_{2\restriction A}\) for some measurable \(A\subseteq X\), we have \(\rho_Y(U_1)_{\restriction\pi\inv(A)}=\rho_Y(U_2)_{\restriction \pi\inv(A)}\).
	Now take \(T\in [G]\), let \((A_n)\) be measurable partition of \(X\) and \((T_n)\) a sequence of elements of \(G\) such that \(T_{\restriction A_n}=T_{n\restriction A_n}\) for every \(n\in\N\), we are forced to define \([\rho_Y](T)\) by
	\[ [\rho_Y](T)_{\restriction \pi\inv(A_n)}=\rho(T_n)_{\restriction \pi\inv(A_n)}\]
	so that \(\pi\) is still support-preserving.
\end{proof}

\subsection{Topologies, automatic continuity and proof of Theorem \ref{thmi: two sin topologies}}

Any full group \(\mathbb G\leq \Aut(X,\mu)\) can be equipped with the uniform topology coming from \(\Aut(X,\mu)\), that is the group topology induced by the \textit{uniform metric} \(d_u\) defined by
\[d_u(S,T)= \mu(\lbrace x\in X\colon S(x)\neq T(x)\rbrace)=\mu(\supp(ST^{-1})).\]

Full groups are always closed for the uniform topology, so the uniform metric restricts to a
complete metric on them \cite[Lem. 5.4]{dyeGroupsMeasurePreserving1959}. The full groups which are moreover separable for the uniform topology
are exactly the full groups of \emph{countable p.m.p. equivalence relations}
\cite[Prop. 3.8]{carderiMorePolishfull2016}.
The following result was proved by Kittrell and Tsankov  for
such full groups \cite{kittrellTopologicalPropertiesFull2010}, but their proof extends verbatim to the general case,
as was already observed by Ben Yaacov, Berenstein and Melleray in the case \(\mathbb G=\Aut(X,\mu)\)
\cite{benyaacovPolishTopometricGroups2013} (see \cite[Thm.~3.18]{lemaitreGroupesPleinsPreservant2014} for an explicit general statement).

\begin{thm}\label{thm: automatic cont}
	Let \(\mathbb G\leq\Aut(X,\mu)\) be an ergodic full group.
	Then every group homomorphism \(\mathbb G\to H\), where \(H\) is a Polish group, has
	to be continuous with respect to the uniform topology on \(\mathbb G\).
\end{thm}
We will only need the following corollary.

\begin{cor}\label{cor: autocont action on malg}
	Let \(\mathbb G\) be an ergodic full group, suppose that we have a boolean action
	\(\rho:\mathbb G\rightarrow \Aut(Y,\nu)\). Then \(\rho\) is uniform to weak continuous.
\end{cor}
\begin{proof}
	The corollary follows from the above theorem and the fact that \(\Aut(X,\mu)\) equipped with the weak topology is a Polish group.
\end{proof}

The uniform metric on \(\Aut(X,\mu)\) is bi-invariant, which implies that the
uniform topology is SIN (the identity element admits a basis of conjugacy invariant neighborhoods).
We now prove that the discrete topology is the only other SIN topology on any ergodic full group.

\begin{thm}\label{thm: SIN topos}
	Let \(\mathbb G\leq\Aut(X,\mu)\) be an ergodic full group, let \(\tau\) be a SIN group topology on \(\mathbb G\).
	Then \(\tau\) is either the discrete or the uniform topology.
\end{thm}

We will use the fact that the uniform topology can be interpreted as a topology of uniform convergence:

\begin{prop}[{\cite[Lem. 5.4]{dyeGroupsMeasurePreserving1959}}]\label{prop: uniform cv}
	The metric \(\delta_u\) on \(\Aut(X,\mu)\) defined by
	\[
		\delta_u(T,U)\coloneqq \sup_{A\in\MAlg(X,\mu)}\mu(T(A)\bigtriangleup U(A))=2\sup_{A\in\MAlg(X,\mu)}\mu(T(A)\setminus U(A))
	\]
	is equivalent to the uniform metric \(d_u\).
\end{prop}

In the proof of the theorem we will need the following lemma.

\begin{lemma}\label{lem:commutatortrick}
	Let \(\mathbb G\) be an ergodic full group and
	let \(\mathcal V\subseteq \mathbb G\) be a conjugacy-invariant symmetric subset.
	Let  \(T\in\mathcal V\) and \(A\subseteq X\) measurable.
	Then \(\mathcal V^2\) contains every involution whose support has measure at most \[2 \mu(TA\setminus A).\]
\end{lemma}
\begin{proof}
	Let  \(\eps\in[0,\mu(TA\setminus A)]\). We have to show that \(\mathcal V^2\) contains all involutions whose support has measure \(2\eps\).

	Set \(B\coloneqq A\setminus TA\). Let \(U\) be an involution whose support has measure \(\eps\) and is contained in \(B\). Clearly \(TUT^{-1}\) is an involution whose support is contained in \(TB\).
	Since \(B\cap TB=\emptyset\), \(U(TUT\inv)\) is an involution and its support has measure \(2\mu(\supp(U))=2\eps\). Since \(\mathcal V^2\) is conjugacy invariant \(UTUT^{-1}\in \mathcal V T\mathcal VT\inv=\mathcal V^2\).
	Therefore \(\mathcal V^2\) contains an involution whose support has measure \(2\eps\). Finally, since \(\mathcal V\) is conjugacy invariant, so is \(\mathcal V^2\),
	and Proposition \ref{prop:conjugacy of invol} implies that \(\mathcal V^2\) must contain all involution whose support has measure \(2\eps\) as wanted.
\end{proof}

\begin{proof}[Proof of Theorem \ref{thm: SIN topos}]
	We start by proving \(\tau\) has to refine the uniform topology.
	By Proposition \ref{prop: uniform cv}, a basis of neighborhoods of the identity for the uniform topology is given by sets of the form \[\mathcal V_\eps\coloneqq \{T\in\mathbb G\colon \mu(T(A)\setminus A)<\eps \text{ for all }A\subseteq X\text{ measurable}\}.\]
	Therefore let us fix \(\eps\in (0,\frac 12)\). We have to exhibit a \(\tau\)-neighborhood of the identity \(\mathcal W\) such that \(\mathcal W\subseteq \mathcal V_\eps\).

	Take an involution \(U\) whose
	whose support has
	measure \(2\eps\). Since the topology is Hausdorff and SIN, we
	find a conjugacy-invariant \(\tau\)-neighborhood
	of the identity \(\mathcal W'\) which does not contain \(U\). Choose a conjugacy-invariant and symmetric \(\tau\)-neighborhood of the identity \(\mathcal W\) such that \(\mathcal W^2\subseteq\mathcal W'\).
	Then Lemma \ref{lem:commutatortrick} and the fact that \(\mathcal W^2\) does not contain \(U\) imply that \(\mu(TA\setminus A)<\eps\) for every \(T\in \mathcal W\), that is \(\mathcal W\subseteq \mathcal V_\eps\). So \(\tau\) does refine the uniform topology.

	Let us now show that if \(\tau\) is not discrete, it has to be equal to the uniform topology. So given a \(\tau\)-neighborhood of the identity \(\mathcal W\), we have to show that there is \(\eps>0\) such that \(\mathcal W\supseteq \{T\in\mathbb G:\ d_u(T,\id)\leq \eps\}.\)

	Since \(\tau\) is not discrete, there there exists a conjugacy-invariant, symmetric
	\(\tau\)-neighborhood of the identity \(\mathcal V\subseteq \mathcal W\) such that \(\mathcal V\neq\{\id\}\) and \(\mathcal V^6\subset\mathcal  W\). Let \(T\in \mathcal V\setminus \{\id\}\) and let \(A\subseteq X\) such that \(\mu(TA\setminus A)>0\). Then Lemma \ref{lem:commutatortrick}, implies that \(\mathcal V^2\) contains every involution whose support has measure at most \(\epsilon\coloneqq 2\mu(TA\setminus A)\). By Ryzhikov theorem, Theorem \ref{thm: 3invols}, every element of \(\mathbb G\) whose support has measure at most \(\epsilon\) is the product of \(3\) involution whose support also has measure at most \(\epsilon\). So \(\mathcal W\supseteq \mathcal V^6\) contains every element whose support has measure at most than \(\epsilon\) as claimed.
\end{proof}

We will only use the following corollary when the target group \(G\) is \(\Aut(Y,\nu)\) equipped with the uniform topology.

\begin{crl}\label{crl:automemb}
	Let \(\mathbb G\) be an ergodic full group and let \(G\) be a group with a SIN group topology \(\tau\).
	Assume that we have a homomorphism \(\rho:\mathbb G\rightarrow G\).
	If \(\rho(\mathbb G)\) is not discrete, then \(\rho\) is uniform-to-\(\tau\) continuous.
\end{crl}
\begin{proof}
	We may as well assume there is some \(T\in\mathbb G\) such that \(\rho(T)\neq 1\) because
	otherwise \(\rho\) is clearly continuous. It follows that \(\rho\) is injective since \(\mathbb G\)
	is simple by Theorem \ref{thm:ergodicvssimple}.
	Consider the pullback topology \(\rho^*(\tau)\) on \(\mathbb G\),
	that is the topology generated by the pre-images of open sets in \(G\).
	Since \(\rho\) is injective and \(\tau\) is Hausdorff, \(\rho^*(\tau)\) is Hausdorff.
	Moreover \(\tau\) is SIN so \(\rho^*(\tau)\) is a SIN group topology.
	So by the above theorem, if the image of \(\mathbb G\) is not discrete,
	\(\rho^*(\tau)\) has to be the uniform topology and hence \(\rho\) is continuous
	(actually an embedding) of topological groups.
\end{proof}

We conclude this section by recalling a fundamental tool in the study of full groups, namely induced bijections.
Given \(T\in\Aut(X,\mu)\) and \(A\subseteq X\) measurable, for almost all \(x\in A\) there is \(n\geq 1\) such that \(T^n(x)\in A\).
If we define \(\tau_A(x)\) as the smallest such \(n\), and let \(\tau_A(x)=0\) whenever \(x\not\in A\), the bijection \(T_A\) \textbf{induced} by \(T\) on \(A\) is by definition given by: for all \(x\in X\),
\[
	T_A(x)=T^{\tau_A(x)}(x).
\]
Observe that \(T_A\) always belongs to the full group \([T]\) generated by \(T\).
The following well-known continuity property follows from \cite[Lem.~3]{keaneContractibilityAutomorphismGroup1970}.
\begin{proposition}\label{lem: continuity induced transfo}
	Let \(T\in\Aut(X,\mu)\).
	Endow the measure algebra \(\MAlg(X,\mu)\) with the metric \(d_\mu(A,B)\coloneqq \mu(A\bigtriangleup B)\) and the full group of \(T\) with the uniform topology.
	Then the map
	\[
		\begin{array}{rcl}
			\MAlg(X,\mu) & \to     & [T] \\
			A            & \mapsto & T_A
		\end{array}
	\]
	is continuous.
\end{proposition}

\subsection{P.m.p.\ equivalence relations and their full groups}\label{sec: pmp eq rel and full groups}

Let us start by recalling some terminology; the unfamiliar reader is refered to \cite{kechrisTopicsOrbitEquivalence2004}.
Let \(X\) be a standard Borel space.
A \textbf{Borel equivalence relation} \(\rel\) on \(X\) is an equivalence relation which is Borel when one sees it as a subset of the product space \(\rel\subseteq X\times X\) equipped with the product \(\sigma\)-algebra. We will be mainly interested in equivalence relation with \textit{countable classes}, also called countable equivalence relations.
The prototypical example of Borel countable equivalence relations is obtained as follows: given a countable group \(\Gamma\) acting on \(X\) by Borel bijections, we get a Borel equivalence relation \(\mathcal R_\Gamma\) on \(X\) whose classes are the \(\Gamma\)-orbits, defined by
\[(x,y)\in\mathcal R_\Gamma\text{ if } y\in\Gamma\cdot x.\]
The Feldman-Moore theorem ensures us that conversely every countable Borel equivalence relation comes from a Borel action of some countable group \(\Gamma\) on \(X\).

Define the Borel full group
of a countable Borel  equivalence relation \(\mathcal R\)  as the group of all Borel bijections \(T\) of \(X\) such that \((x,T(x))\in\mathcal R\) for every \(x\in X\). Note that if \(\mathcal R\) is a countable Borel equivalence relation and \(T:X\to X\) is a Borel function such that \(T\) restricts to a bijection on each \(\mathcal R\)-class, then \(T\) belongs to the Borel full group of \(\mathcal R\).
Also observe that if \(\mathcal R=\mathcal R_\Gamma\), then  \(T:X\to X\) is in the Borel full group of \(\mathcal R\) if and only if there is a partition of \(X\) into Borel subsets \((A_\gamma)_{\gamma\in\Gamma}\) such that \((\gamma A_\gamma)_{\gamma\in\Gamma}\) is also a partition of \(X\) and for all \(\gamma\in\Gamma\) and all \(x\in A_{\gamma}\), \(T(x)=\gamma\cdot x\).

Now suppose that \((X,\mu)\) be a standard probability space, i.e.\ \(X\) is a standard Borel space equipped with a Borel probability measure. It follows from the previous discussion that if we now have a Borel action by measure-preserving bijections of a countable group \(\Gamma\) on \((X,\mu)\), then every element of the Borel full group of \(\mathcal R_\Gamma\) preserves the measure \(\mu\). This means that \(\mathcal R_\Gamma\) is a \pmp equivalence relation as defined below.

\begin{definition}
	A countable Borel equivalence relation \(\mathcal R\) on a standard probability space \((X,\mu)\) is \textbf{measure-preserving} when the elements of its Borel full group preserve the measure \(\mu\).
	We will simply call \textbf{\pmp}equivalence relations such equivalence relations.
\end{definition}

\begin{definition}
	The \textbf{full group} of a Borel equivalence relation \(\rel\) is the group
	\[[\rel]\coloneqq\lbrace T\in \Aut(X,\mu)\colon (T(x),x)\in\rel\text{ for all }x\in X \rbrace.\]
\end{definition}

Observe that when \(\mathcal R\) is p.m.p., the group \([\rel]\) is the quotient of the Borel full group of \(\mathcal R\) by the normal subgroup of elements whose support has measure zero.
It is not hard to check that \([\rel]\) is always a full group on the sense of Definition \ref{df: full group}. As recalled before, the restriction of the uniform distance \(d_u\) on \([\rel]\) is complete separable and thus defines a Polish group topology on \([\rel]\).

\begin{remark}
	As discovered in \cite{carderiMorePolishfull2016}, there are actually many full groups of uncountable Borel equivalence relations whose full group is a Polish group.
\end{remark}

We say that a \pmp equivalence relation \(\rel\) is \textbf{ergodic} if every measurable  \(A\subseteq X\) which is a union of \(\rel\)-classes has measure either \(0\) or \(1\). When \(\mathcal R=\mathcal R_\Gamma\), this is equivalent to every \(\Gamma\)-invariant set having measure \(0\) or \(1\), which is equivalent to every almost everywhere invariant set having measure \(0\) or \(1\). It follows that a p.m.p.\ equivalence relation is ergodic if and only if its full group is.

\begin{definition}
	Say that a \pmp equivalence relation \(\mathcal R\) is \textbf{finite} when it has only finite classes, and  \textbf{hyperfinite} when it can be written as an increasing union of \pmp finite subequivalence relations.
\end{definition}

Here is the prototypical of an ergodic hyperfinite equivalence relation.
Let us work in the standard Cantor space \((X,\mu)\coloneqq(\{0,1\}^\N,\mathcal B(1/2)^{\otimes\N})\),
equipped with the standard Bernoulli product measure.
\begin{definition}\label{def:R0} Denote by \(\mathcal R_0\) the equivalence relation on \(\{0,1\}^\N\)
	defined by \(((x_n),(y_n))\in\mathcal R_0\) if and only if there is \(N\in\N\) such that \(x_n=y_n\) for all \(n\geq N\).
\end{definition}

Recall from Section \ref{sec: dyadic perm} that we have a natural inclusion of the dyadic symmetric group \(\mathfrak S_{2^\infty}\) into \(\Aut(X,\mu)\).
Observe that the equivalence relation \(\mathcal R_0\) is the equivalence relation associated to this \pmp action, in particular \(\rel_0\) is  a \pmp equivalence relation.
Note that \(\mathcal R_0\) can be written as the increasing union of the finite \pmp equivalence relations \(\mathcal S_k\)  where \(\mathcal S_k\) is
given by \(((x_n),(y_n))\in\mathcal S_k\) if and only if \(x_n=y_n\) for all \(n\geq k\).
We thus see that \(\mathcal R_0\) is hyperfinite, and it is not hard to check that it is ergodic.

We finally define the important notion of orbit equivalence between \pmp equivalence relations.

\begin{definition}
	Let \(\mathcal R\) and \(\mathcal S\) be two \pmp equivalence relations on \((X,\mu)\) and \((Y,\nu)\) respectively.
	A measure-preserving bijection \(\varphi:(X,\mu)\to(Y,\nu)\) is called an \textbf{orbit equivalence} between \(\mathcal R\) and \(\mathcal S\) if up to restricting \(\mathcal R\) and \(\mathcal S\) to full measure subsets, we have \(\varphi\times\varphi(\mathcal R)=\mathcal S\).
	When there is such an orbit equivalence between \(\mathcal R\) and \(\mathcal S\) we say that \(\mathcal R\) and \(\mathcal S\) are orbit equivalent.
\end{definition}

By a theorem of Dye, the equivalence relation \(\mathcal R_0\) is the unique hyperfinite ergodic equivalence relation up to orbit equivalence (\cite[Thm.~3]{dyeGroupsMeasurePreserving1959}, see also \cite[Sec.~7]{kechrisTopicsOrbitEquivalence2004} for a statement and proof in a setup closer to ours. Note that Dye only consider full groups and that the notion of hyperfiniteness is then stated as \emph{approximate finiteness} of the full group; see \cite[Prop.~1.57]{lemaitreGroupesPleinsPreservant2014} for the equivalence between the two notions).

It is well-known that \(\varphi:(X,\mu)\to(Y,\nu)\) is an orbit equivalence between \pmp equivalence relations \(\mathcal R\) on \((X,\mu)\) and \(\mathcal S\) on \((Y,\nu)\) if and only if \(\varphi [\mathcal R]\varphi\inv =[\mathcal S]\) (see for instance \cite[Prop.~1.30]{lemaitreGroupesPleinsPreservant2014}).
Moreover, Dye's reconstruction theorem \cite[Thm.~2]{dyeGroupsMeasurePreserving1963} yields that any abstract isomorphism between ergodic full groups must come from an isomorphism of the underlying measured spaces so that full groups of ergodic \pmp equivalence relations completely remember the equivalence relation up to orbit equivalence. Our main result allows us to recover Dye's reconstruction theorem with a completely different proof (see Section \ref{sec:dyereconstruction}).

\subsection{Products and symmetric products of \pmp equivalence relations}
\begin{df}
	Let \(\rel\) be a countable p.m.p.\ equivalence relation on \((X,\mu)\) and let \(n\in\mathbb N^*\) a positive integer. Then we define the \textbf{product equivalence relation} \(\rel^n\) on \((X^n,\mu^n)\) by declaring \[((x_i)_{i=1}^{n},(y_i)_{i=1}^n)\in\rel^n\text{ if for every }i\in\{1,\dots,n\},\ (x_i,y_i)\in\rel. \]

	Recall the definition of the symmetric \(n\)-power \((X^{\odot n},\mu^{\odot n})\) from Section \ref{sec: sym diag}. We then define the \textbf{symmetric product equivalence relation}  \(\rel^{\odot n}\) on \((X^{\odot n},\mu^{\odot n})\) by declaring \[([x_i]_{i=1}^n,[y_i]_{i=1}^n)\in\rel^{\odot n}\text{ if there is }\sigma\in \mathfrak S_n\text{ such that for all }i\in\{1,\dots,n\},\ (x_i,y_{\sigma(i)})\in\rel.\]
\end{df}

One can easily check that both \(\rel^n\) and \(\rel^{\odot n}\) are countable p.m.p.\ equivalence relations. Here is one way of seeing this.

Assume that the equivalence relation \(\rel\) is generated by an action \(\alpha\) of the countable group \(\Gamma\) on the probability space \((X,\mu)\). Then the action \(\alpha^n\) of \(\Gamma^n\) on \(X^n\) defined by \[\alpha^n(\gamma_1,\dots,\gamma_n)(x_1,\dots,x_n)=(\alpha(\gamma_1)x_1,\dots,\alpha(\gamma_n)x_n)\]
generates the equivalence relation \(\rel^n\).
Now consider the group \(
	\Gamma^n\rtimes \mathfrak S_n\) where \(\mathfrak S_n\) acts on \(\Gamma^n\) by permuting the copies. Then the \(\Gamma^n\)-action \(\alpha^{n}\) on \(X^n\) naturally extends to a \(\Gamma^n\rtimes \mathfrak S_n\)-action by declaring \[\alpha^{ n}((\gamma_1,\dots,\gamma_n)\sigma)(x_1,\dots,x_n)\coloneqq \left(\alpha(\gamma_1)x_{\sigma\inv(1)},\dots\alpha(\gamma_n)x_{\sigma\inv(n)}\right).\]
Let \(D\subseteq X^n\) be a measurable fundamental domain of the action of \(\mathfrak S_n\) on the full measure subset \(X^{(n)}=\{(x_1,\dots,x_n)\colon x_i\neq x_j\text{ for all }i\neq j\}\) (e.g.~fix a Borel linear order \(<\) on  \(X\) and let \(D\coloneqq\{(x_1,\dots,x_n)\in X^n\colon x_1<x_2<\ldots<x_n\}\)).
The projection \(\pi_n\colon X^{(n)}\rightarrow X^{\odot n}\) restrict to an isomorphism \(\pi_n\colon D\to X^{\odot n}\) which scales the measure by \(n!\).
Then, it is easy to show that \(\pi_n\) defines an \textit{orbit equivalence}\footnote{Two \pmp equivalence relations \(\mathcal R\) and \(\mathcal S\) are orbit equivalent if there is a \pmp map which induces an isomorphism between \(\mathcal R\) and \(\mathcal S\) up to a null set.}
between the restriction of the orbit equivalence relation of \(\alpha^n(\Gamma^n\rtimes \mathfrak S_n)\) on \(D\) and \(\rel^{\odot n}\), that is given \(x,y\in D\)
\[\text{there is }g\in \Gamma^n\rtimes \mathfrak S_n\text{ such that } y=\alpha^n(g)x\text{ if and only if }(\pi_n(x),\pi_n(y))\in\rel^{\odot n}.\]
This actually shows the following lemma.
\begin{lemma}\label{lem:SOE for R odot n}
	The \pmp equivalence relation \(\mathcal R^{\odot n}\) is stably orbit equivalent\footnote{Two \pmp equivalence relations \(\mathcal R\) and \(\mathcal S\) on \((X,\mu)\) are stably orbit equivalent if there are positive measure sets \(A\) and \(B\) such that their respective restrictions to \(A\) and \(B\) (equipped with appropriate renormalized measures) are orbit equivalent.} to the equivalence relation
	\[
	\mathcal R^n\rtimes \mathfrak S_n\coloneqq
		\{((x_i)_{i=1}^n,(y_i)_{i=1}^n)\colon \exists \sigma\in\mathfrak S_n
		\forall i\in\{1,\dots,n\}, (x_i,y_{\sigma(i)})\in\mathcal R\}. \qed \]
\end{lemma}

The construction of \(\mathcal R^{\odot n}\) is naturally compatible with the symmetric diagonal action as follows.

\begin{lemma}\label{lem: fg of R gens fg of Rn}
	Let \(\mathcal R\) be an ergodic countable p.m.p.~equivalence relation and let \(n\in\N^*\) be a positive integer.
	Consider the symmetric diagonal boolean action \(\iota^{\odot n}\colon [\rel]\to\Aut(X^{\odot n},\mu^{\odot n})\).
	Then \(\iota^{\odot n}([\rel])\leq [\rel^{\odot n}]\) and the full group generated by \(\iota^{\odot n}([\mathcal R])\) is equal to
	\([\mathcal R^{\odot n}]\).
\end{lemma}
\begin{proof}
	By the Feldman-Moore theorem, we may and do fix a measure-preserving action \(\alpha\) of a countable group \(\Gamma\) such that \(\mathcal R=\mathcal R_\Gamma\). We thus get an action \(\alpha^{ n}\) of \(\Gamma^n\rtimes \mathfrak S_n\) on \(X^n\) as explained above.

	Let again \(D\subseteq X^n\) be a measurable fundamental domain of the action of \(\mathfrak S_n\) on \(X^n\). By the discussion which precedes the lemma, for every \(T\in [\rel^{\odot n}]\) there is a measurable partition \(\{D^g\}_{g\in \Gamma^n\rtimes \mathfrak S_n}\) of \(D\) such that
	\[T\pi_n(x_1,\ldots,x_n)=\pi_n(\alpha^n(g)(x_1,\ldots,x_n))\] for every \((x_1,\ldots,x_n)\in D^g\).
	Up to discarding the null set of non-injective \(n\)-uples in \(X^n\), we have a countable cover \[\{A_1^i\times\ldots\times A_n^i\}_{i\in\mathbb N}\] of \(X^n\) such that for every \(i\in\mathbb N\), the sets \(A_1^i,\ldots,A_n^i\) are pairwise disjoint.

	Let us now fix \(g=(\gamma_1,\ldots,\gamma_n)\sigma\in \Gamma^{\odot n}\). Observe that for all \(A_1,\ldots,A_n\subseteq X\) measurable we have
	\[
	\alpha^n(g)(A_1\times\ldots\times A_n)=
	\alpha(\gamma_1)A_{\sigma\inv(1)}\times\ldots\times \alpha(\gamma_n)A_{\sigma\inv(n)}.
	\]
	Now for every \(i,j\in\mathbb N\) and \(k\in\{1,\ldots,n\}\), we set \[C^{g,i,j}_{k}\coloneqq A^i_k\cap \alpha(\gamma_k)^{-1}A^j_{\sigma^{-1}(k)}.\]

	Since \(\{A_1^i\times\ldots\times A_n^i\}_{i\in\mathbb N}\) is a cover of \(X^n\), we first have that \[\left\lbrace\alpha^n(\gamma_1,...,\gamma_n)\inv (A^j_{\sigma\inv(1)}\times \ldots\times A^j_{\sigma\inv(n)})\right\rbrace_{j\in\N}\] is a cover of \(X^n\).
	We deduce that
	\(\{C^{g,i,j}_1\times\ldots\times C^{g,i,j}_n\}_{i,j\in\N}\) forms a cover of \(X^n\) and hence of \(D^g\).
	Let us now fix \(i,j\in\mathbb N\). By construction for all \(k\in\{1,\ldots, n\}\) we have
	\(\alpha(\gamma_k)C^{g,i,j}_k\subseteq A^j_{\sigma\inv(k)}\), so the sets \[\alpha(\gamma_1)C^{g,i,j}_{1},\dots,\alpha(\gamma_n)C^{g,i,j}_{n}\]
	are pairwise disjoint. Therefore using Proposition \ref{prop: transitive on equal measure} we can find \(U\in [\rel]\) such that for all \(k\in\{1,\ldots,n\}\),
	\[
	U(x)=\alpha(\gamma_k)x\text{ for every }
	x\in C_{k}^{g,i,j}.
	\]
	Finally remark that by construction, for every \((x_1,\ldots,x_n)\in D^g\cap C_1^{g,i,j}\times\ldots\times C_n^{g,i,j}\) we have
	\begin{align*}
		\iota^{\odot n}(U)\pi_n(x_1,\ldots,x_n)
		 & =\pi_n(\alpha(\gamma_1)x_1,\ldots,\alpha(\gamma_n)x_n) \\
		 & =\pi_n(\alpha^n(g)(x_1,\ldots,x_n))                    \\
		 & =T\pi_n(x_1,\ldots,x_n).\end{align*}
	Therefore \(T\) can be obtained by cutting and pasting elements of \(\iota^{\odot n}([\rel])\) and hence the proof is completed.
\end{proof}

\subsection{Boolean actions of full groups coming from actions of the equivalence relations}	\label{sec: boolean action from eq rel}

Let \(\mathcal R\) be a p.m.p.\ equivalence relation on the probability space \((X,\mu)\). Let \(\mathcal S\) be a p.m.p.\ equivalence relation on \((Y,\nu)\) and let \(\pi:Y\rightarrow X\) a measure-preserving map.
\begin{itemize}
	\item We say that \(\pi\) is a \textbf{factor map} if for every \((y_1,y_2)\in\mathcal S\), one has that \((\pi(y_1),\pi(y_2))\in\mathcal R\).
\end{itemize}

Clearly \(\pi\) is a factor map if and only if the map \(\pi\times \pi:Y\times Y\rightarrow X\times X\) restrict to a map from \(\mathcal S\) to \(\mathcal R\). Therefore we will also use the convention to denote a factor map as a map \(\pi\colon \mathcal S\rightarrow \mathcal R\). Observe also that \(Y\) can be seen inside \(\mathcal S\) as the set \(\Delta_Y\coloneqq\lbrace (y,y)\colon y\in Y\rbrace\) and \(\pi\) is compatible with this restriction.

\begin{itemize}
	\item  We say that a factor map \(\pi\colon \mathcal S\rightarrow \mathcal R\) is \textbf{class-bijective}, if there is \(Y_0\subseteq Y\) of full measure such that for all \(y_1,y_2\in Y_0\), if \((y_1,y_2)\in \mathcal S\) and \(y_1\neq y_2\), then \(\pi(y_1)\neq \pi(y_2)\).
\end{itemize}

If \(\pi\colon\mathcal S\rightarrow\mathcal R\) is a class-bijective factor map, then \(\mathcal S\) is often called a \textbf{class-bijective extension} of \(\mathcal R\) and we say that \(\pi\) is a \textbf{measure-preserving action} of \(\mathcal R\).

\begin{example}
	Let \(\mathcal R\) be a \pmp equivalence relation on \((X,\mu)\), let \(n\geq 1\).
	Then \(\mathcal R^n\) is a class-bijective extension of \(\mathcal R^{\odot n}\) via \(\pi: X^n\to X^{\odot n}\).
	Indeed for almost every \((x_1,\dots,x_n)\), we have that the \(x_i\)'s belong to distinct \(\mathcal R\)-classes. This yields that the restriction of \(\pi\) to the \(\mathcal R^n\) equivalence class of \((x_1,\dots,x_n)\) (which is equal to \([x_1]_{\mathcal R}\times\cdots\times [x_n]_{\mathcal R}\)) is injective as wanted.
\end{example}

Measure-preserving actions of equivalence relations give rise to boolean actions of the associated full groups.

\begin{lemma}\label{lemma: def rho pi}
	Let \(\mathcal R\) be a \pmp equivalence relation on \((X,\mu)\) and \(\mathcal S\) be a \pmp equivalence relation on \((Y,\nu)\).
	Let \(\pi\) be a class-bijective factor map from \(\mathcal S\) to \(\mathcal R\) and for \(T\in [\mathcal R]\) and \(y\in Y\) define
	\[\rho_\pi(T)y\text{ as the unique }y'\in Y\text{ such that }\pi(x)=T\pi(y)\text{ and }(y,y')\in \mathcal S.\]
	Then \(\rho_\pi\colon [\mathcal R]\rightarrow [\mathcal S]\) is a group homomorphism. Moreover \(\rho_\pi\) is support-preserving: for all \(T\in[\mathcal R]\) we have \(\supp(\rho_\pi(T))=\pi^{-1}(\supp(T))\), and the full group generated by \(\rho_\pi([\mathcal R])\) is equal to \([\mathcal S]\).
\end{lemma}
\begin{proof}
	Fix \(T\in [\mathcal R]\), then remark that
	\[\pi^{-1}(\lbrace (x,Tx)\in\mathcal R\colon x\in X\rbrace)=\lbrace (y,\rho_\pi(T)y)\in \mathcal S\colon y\in Y\rbrace\]
	so the graph of \(\rho_\pi(T)\) is Borel, which implies that \(\rho_\pi(T)\) is a Borel map (see \cite[Thm.~14.12]{kechrisClassicaldescriptiveset1995}).
	By construction \(\rho_\pi(T)\) preserves \(\mathcal S\)-classes, and since \(\pi\) is class-bijective, \(\rho_\pi(T)\) induces a bijection on each \(\mathcal S\)-class. It follows that \(\rho_\pi(T)\) is in the full group of \(\mathcal S\).
	To see that \(\rho_\pi\) is a group homomorphism, observe that
	\[\pi(\rho_\pi(T_1)\rho_\pi(T_2)y)=T_1\pi(\rho_\pi(T_2)y))=T_1T_2\pi(y).\]
	We have by definition \(\rho_\pi(T)(y)=y\) if and only if \(T(\pi(y))=\pi(y)\), so \(\rho_\pi\) is support-preserving.

	Finally, to see that the full group generated by \(\rho_\pi([\mathcal R])\) is equal to \([\mathcal S]\), let \(\Gamma\leq[\mathcal R]\) be a countable group such that \(\mathcal R=\mathcal R_\Gamma\), then since \(\pi\) is class-bijective we have that \(\mathcal S\) is generated by the action \(\rho_\pi\) of \(\Gamma\). In particular the full group of \(\mathcal S\) is contained in the full group generated by \(\rho_\pi(\Gamma)\).
	Since the full group generated by \(\rho_\pi(\Gamma)\) is contained in the full group generated \(\rho_\pi([\mathcal R])\), which is itself contained in \([\mathcal S]\) by construction, the conclusion follows.
\end{proof}

The above lemma has a converse that will play an important role.

\begin{lemma}\label{lem:pmp action from support pres}
	Let \(\mathcal R\) be a p.m.p.\ equivalence relation on \((X,\mu)\) and consider a boolean action \(\rho:[\mathcal R]\rightarrow \Aut(Y,\nu)\) such that \(\rho\) factors in a support-preserving way on the inclusion \(\iota:[\mathcal R]\rightarrow \Aut(X,\mu)\). Then there is a p.m.p.\ equivalence relation \(\mathcal S\) on \(Y\) and a class-bijective factor map \(\pi:\mathcal S\rightarrow \mathcal R\) such that \(\rho=\rho_\pi\).
\end{lemma}
\begin{proof}
	Let \(\Gamma\) be a countable group acting on \((X,\mu)\) so that \(\mathcal R=\mathcal R_\Gamma\).
	Denote by \(\pi:Y\rightarrow X\) the support-preserving factor map.
	Let \(\mathcal S\) be the orbit equivalence relation of the action \(\rho\) of \(\Gamma\).
	The map \(\pi\) intertwines by definition the two actions of \(\Gamma\),
	that is \(\pi\colon\mathcal S\rightarrow\mathcal R\) is a factor map.
	Moreover the hypothesis that \(\pi\) is support-preserving implies that 
	\(\pi\) is class-bijective (see Remark \ref{rmk: support pres is class bij}).

	Finally observe that \(\rho_\pi(\gamma)=\rho(\gamma)\) for every \(\gamma\in \Gamma\).
	Since \(\rho\) and \(\rho_\pi\) are both support-preserving, they agree on the full group generated by \(\alpha(\Gamma)\), which finishes the proof.
\end{proof}

\begin{example}
	Let \(\mathcal R\) be a \pmp equivalence relation on \((X,\mu)\), let \(n\geq 1\).
	The diagonal action \(\iota^n:[\mathcal R]\to [\mathcal R^n]\) factors  onto
	\(\iota^{\odot n}: [\mathcal R]\to [\mathcal R^{\odot n}]\) via the quotient map \(\pi:X^n\to X^{\odot n}\),
	and this factor map is support-preserving since
	\[\pi\inv(\Fix(\iota^{\odot n}(T)))=\pi\inv([\Fix T,\dots,\Fix T])=\Fix T\times\cdots\times\Fix T=\Fix(\iota^n(T)).\]
	The full group generated by \(\iota^{\odot n}([\mathcal R])\) is equal to \([\mathcal R^{\odot n}]\) so by Proposition \ref{prop: support preserving extend to full group} \(\iota^n\) extends in a support preserving manner to \(\widetilde{\iota^n}:[\mathcal R^{\odot n}]\to[\mathcal R^n]\), which can be rewritten as \(\widetilde{\iota^n}=\rho_\pi\) when viewing \(\pi\) as a class-bijective factor map \(\mathcal R^n\to \mathcal R^{\odot n}\).
\end{example}

We refer the reader to Section \ref{sec:boolean actions} for the notion of boolean action of a \pmp equivalence relation, which as we prove there is essentially the same as a \pmp action (this is needed for our characterization of property (T) in terms of full groups).

\section{Diagonal support-preserving factorization for boolean actions}\label{sec:proof of main thm}

The aim of this section is to prove Theorem \ref{thm:main} which we restate here.

\begin{thm}\label{thm: main - in section}
	Let \(\mathbb G\) be an ergodic full group over the standard probability space \((X,\mu)\) and let
	\(\rho\colon\mathbb G\to \Aut(Y,\nu)\) be a boolean action on another standard
	probability space \((Y,\nu)\). Then there is a unique measurable \(\rho(\mathbb G)\)-invariant
	partition \(\{Y_n\}_{n=0,\ldots,\infty}\) of \(Y\) such that
	the boolean actions  \(
		\rho_n=\rho_{\restriction \left(Y_n,\frac{\nu_{\restriction Y_n}}{\nu(Y_n)}\right)}\)
	are subject to the following conditions
	\begin{enumerate}
		\item the map \(\rho_0\) maps every element to the identity on \(Y_0\);
		\item for all \(n\in\mathbb N^*\), there is a (unique) measure-preserving map
		      \(\pi_n\colon Y_n\rightarrow X^{\odot n}\) which is a support-preserving factor map from the boolean action \(\rho_n\) to \(\iota^{\odot n}\colon \mathbb G\hookrightarrow \Aut(X^{\odot n},\mu^{\odot n})\);
		\item for every \(T\in\mathbb G\setminus \{\id\}\) we have that
		      \(\{y\in Y_\infty:\ \rho_\infty(T)y=y\}\) is a null-set.
	\end{enumerate}
\end{thm}

Let \(\mathbb G\) be an ergodic full group. We first observe that Corollary \ref{cor:every set is some support} yields that the inclusion \(\mathbb G\leq \Aut(X,\mu)\) is highly absolutely non-free. Therefore we can apply Proposition \ref{prop:diagfactdec}: Theorem \ref{thm: main - in section} is implied by the following theorem.

\begin{theorem}\label{thm:diagonally factorizable}
	Let \(\mathbb G\leq \Aut(X,\mu)\) be an ergodic full group and let \(\rho\colon\mathbb G\rightarrow \Aut(Y,\nu)\) be a boolean action. Then \(\rho\) is diagonally support-preserving factorizable.
\end{theorem}

The rest of this section is devoted to the proof of Theorem \ref{thm:diagonally factorizable}.
We will first prove it in the hyperfinite case although it is not needed for the general case.
Indeed, the proof in the hyperfinite case is much easier and we find it of independent interest.
Nevertheless, the hasty reader should jump to Section \ref{sec: general case setup}.

\subsection{Proof in the hyperfinite case}

First recall that by Dye's theorem, there is up to isomorphism only one full group of ergodic \pmp hyperfinite equivalence relation, namely the full group \(\mathbb G=[\mathcal R_0]\) (see Definition \ref{def:R0} and the paragraphs below).
We will thus prove Theorem \ref{thm:diagonally factorizable} when \(\mathbb G=[\mathcal R_0]\).
This will rely on the following fact
several times.

\begin{proposition}[{see \cite[Prop.\ 2.8]{kechrisGlobalaspectsergodic2010})}]\label{prop:closuredyadic}
	The closure of \(\mathfrak S_{2^\infty}\leq \Aut(X,\mu)\) in the uniform topology is the full group \([\mathcal R_0]\).
\end{proposition}


Let us fix a boolean action \(\rho\colon [\mathcal R_0]\rightarrow \Aut(Y,\nu)\). By automatic continuity (see Theorem \ref{thm: automatic cont}) we obtain that \(\rho\) is uniform-to-weak continuous. We then have the following straightforward consequence.

\begin{lem}\label{lem:hyperfinite same invariant}
	Let \(A\in\MAlg(Y,\nu)\) such that \(\rho(\sigma)A=A\) for every \(\sigma\in\mathfrak S_{2^\infty}\). Then \(\rho(T)A=A\) for every \(T\in [\mathcal R_0]\).
\end{lem}
\begin{proof}
	By uniform-to-weak continuity, the set of all \(T\in[\mathcal R_0]\) such that \(\rho(T)A=A\) is closed in the uniform topology, and by assumption it contains the dense subgroup \(\mathfrak S_{2^\infty}\) (see Proposition \ref{prop:closuredyadic}), so it has to be equal to \([\mathcal R_0]\).
\end{proof}

\begin{lem}\label{lem:hyperfinite restriction free}
	It the boolean action \(\rho: [\mathcal R_0]\to\Aut(Y,\nu)\) is not free,
	then its restriction  to the subgroup \(\mathfrak S_{2^\infty}\) is not free either.
\end{lem}
\begin{proof}
	We prove the contrapositive; assume that the restriction of \(\rho\) to \(\mathfrak S_{2^\infty}\) is free.
	Since \(\mathcal R_0\) is ergodic, all involutions in \([\mathcal R_0]\) with support of the same measure are conjugate, see Proposition \ref{prop:conjugacy of invol}. Remark that for every \(n\geq 1\), there is an involution in \(\mathfrak S_{2^\infty}\) whose support in \(\Aut(X,\mu)\) has measure \(2^{-n}\). Since the action of \(\mathfrak S_{2^\infty}\) is free and all involutions of support of measure \(2^{-n}\) are conjugate in \([\mathcal R_0]\), we get that \(\supp\rho(U)\) has full measure for every involution \(U\in [\mathcal R_0]\) whose support has measure \(2^{-n}\).

	Take \(T\in [\mathcal R_0]\) different from the identity and consider a subset of positive measure \(A\subset X\) such that \(TA\cap A=\emptyset\).
	For any large enough \(n\geq 1\), consider an involution \(U_n\) with support of measure \(2^{-n-1}\) and contained in \(A\). Remark that \(TU_nT^{-1}\) is an involution with support contained in \(TA\) and in particular \(U_n\) and \(TU_nT^{-1}\) are involutions with disjoint supports.
	So \(U_nTU_nT^{-1}\) is an involution with support of measure \(2^{-n}\) and hence \(\supp\rho(U_nTU_nT^{-1})\) has full measure.
	On the other hand \(U_n\) tends to the identity in \([\mathcal R_0]\). By uniform-to-weak continuity, this implies that \(\rho(U_n)\supp\rho(T)\) tends to \(\supp\rho(T)\) and hence
	\[Y=\supp\rho(U_nTU_nT^{-1})\subseteq \rho(U_n)\supp\rho(T)\cup\supp\rho(T)\rightarrow \supp\rho(T)\]
	which implies that \(\supp T=Y\). The boolean action \(\rho\) is  thus free.
\end{proof}

Recall now that by Proposition \ref{prop:actionss2} the action \(\rho:\mathfrak S_{2^\infty}\to\Aut(Y,\nu)\) is diagonally support-preserving factorizable. We thus have a unique sequence \((Y_i)_{i=0,\ldots,\infty}\) of measurable subsets of \(Y\) such that the action restricted to \(Y_\infty\) is free and the action restricted to \(Y_i\) for \(i<\infty\) factorizes onto the action \(\iota^{\odot i}\) on \(X^{\odot i}\) in a support preserving manner, see Proposition \ref{prop:diagfactdec}. For \(i<\infty\), we denote by \(\pi_i\colon Y_i\rightarrow X^{\odot i}\) the factor maps.

Remark now that Lemma \ref{lem:hyperfinite same invariant} implies that \(Y_i\) is \(\rho([\mathcal R_0])\)-invariant for every \(i=0,\ldots,\infty\). Moreover Lemma \ref{lem:hyperfinite restriction free} implies that the restriction of \(\rho\) on \(Y_\infty\) is free, seen as a boolean action of \([\mathcal R_0]\). Let us first show that \(\pi_i\) induces a factor map for the entire group \([\mathcal R_0]\).

\begin{lem}\label{lem:hyperfinite factor}
	For every \([B_1,\ldots,B_n]\subset X^{\odot i}\) and \(T\in [\mathcal R_0]\) we have that
	\[\rho(T)\pi_i^{-1}([B_1,\ldots,B_n])=\pi_i\inv([TB_1,\ldots,TB_i]).\]
\end{lem}
\begin{proof}
	Fix a sequence \((\sigma_n)_{n\geq 0}\) of elements of \(\mathfrak S_{2^{\infty}}\) converging to \(T\). Then for every \(n\), we have that
	\[\rho(\sigma_n)\pi_i^{-1}([B_1,\ldots,B_n])=\pi_i\inv([\sigma_nB_1,\ldots,\sigma_nB_i]).\]
	The lemma now follows from the continuity of \(\rho\).
\end{proof}

The following lemma is well-known.

\begin{lem}\label{lem:hyperfinite support measure algebra}
	Let \((T_n)_{n\geq 0}\) be a sequence of elements of \(\Aut(X,\mu)\) which converges weakly to \(T\in\Aut(X,\mu)\). Assume moreover that there is \(B\subseteq X\) such that \(\lim_n\mu(B\Delta \supp T_n)=0\). Then \(\supp T\subseteq B\).
\end{lem}
\begin{proof}
	Recall that \(\supp T\) is the supremum of the family of  \(A\in\MAlg(X,\mu)\) such that \(TA\cap A=\emptyset\), so it suffices to show that every such \(A\) is contained in \(B\).

	So take \(A\subseteq X\) such that \(TA\cap A=\emptyset\).
	Since \(T_n\) tends weakly to \(T\), we have that \(\lim_n \mu(T_n A\bigtriangleup TA)=0\). Since \(A\) is disjoint from \(TA\) we obtain \(\lim_n\mu(T_nA\cap A)=0\), and hence \(\lim_n \mu(A\setminus \supp T_n)=0\).
	Since \(\lim_n \mu(B\bigtriangleup \supp T_n)=0\), we conclude that \(\mu(A\setminus B)=0\), so
	\(A\subseteq B\) as wanted.
\end{proof}

We can now finish the proof by showing that the factor map is support preserving.

\begin{lem}
	For every \(T\in [\mathcal R_0]\), we have that \[Y_i\cap \supp\rho(T)=\pi_i^{-1}([\supp T,X,\ldots,X]).\]
\end{lem}
\begin{proof}
	Fix a sequence \((\sigma_n)_{n\geq 0}\) of elements of \(\mathfrak S_{2^{\infty}}\) converging to \(T\) in the uniform topology. By
	Lemma \ref{lem:hyperfinite factor} we have that \[\supp\rho(T)\supseteq \pi_i^{-1}([\supp T,X,\ldots, X]).\]
	Also observe that since \(\sigma_n\) converges to \(T\) in the uniform topology,
	\[\pi_i^{-1}([\supp T,X,\ldots,X])=\lim_n \pi_i^{-1}([\supp \sigma_n,X,\ldots,X])=\lim_n \supp \rho (\sigma_n).\]
	Finally since \(\rho\) is uniform-to-weak continuous, we can use Lemma \ref{lem:hyperfinite support measure algebra}
	to obtain the following diagram of inclusions and equalities:

	\tikzset{%
		symbol/.style={
				draw=none,
				every to/.append style={
						edge node={node [sloped, allow upside down, auto=false]{\(#1\)}} },},}
	\[
		\begin{tikzcd}
			\lim_n\supp\rho(\sigma_n)\cap Y_i\ar[r,symbol=\supseteq]\ar[equal]{d}& \supp\rho(T)\cap Y_i
			\\
			\lim_n\pi_i^{-1}([\supp \sigma_n,X,\ldots, X]) \ar[equal]{r}&
			\pi_i^{-1}([\supp T,X,\ldots, X]). \arrow[u,symbol=\subseteq]
		\end{tikzcd}
	\]
	So the four terms above are  equal, in particular \(\supp \rho(T)\cap Y_i=\pi_i^{-1}([\supp T,X,\ldots, X])\) as wanted.
\end{proof}

Therefore the proof of Theorem \ref{thm:diagonally factorizable} for the case \(\mathbb G=[\mathcal R_0]\) is concluded. We now proceed to the general case. Recall that the general case does not require what we have done in the present section.

\subsection{Setup in the general case}\label{sec: general case setup}

Let \(\mathbb G\) be an ergodic full group over the standard probability space \((X,\mu)\).
Since all the standard probability spaces are isomorphic,
we can assume that \(X=\{0,1\}^\N\) and \(\mu=\mathcal B(1/2)^{\otimes\N}\).
Moreover we can always assume that \(\mathbb G\leq \Aut(X,\mu)\)
contains the standard copy of \(\mathfrak S_{2^\infty}\) acting on \((X,\mu)\).
Indeed by different results of Dye,
every ergodic full group contains the full group of a hyperfinite \pmp ergodic equivalence relation \cite[Thm.~4]{dyeGroupsMeasurePreserving1959},
all these full groups are conjugated in \(\Aut(X,\mu)\) \cite[Thm.~3]{dyeGroupsMeasurePreserving1959} and we already observed that \(\mathfrak S_{2^\infty}\)
is contained the full group of the hyperfinite \pmp ergodic equivalence relation \(\mathcal R_0\).
Therefore from now on, we will always assume that \(\mathfrak S_{2^\infty}\leq \mathbb G\).

Fix a boolean action \(\rho:\mathbb G\rightarrow \Aut(Y,\nu)\).
Let \(A_\infty\coloneqq \bigcap_{\sigma\in\mathfrak S_{2^\infty}\setminus \{\id\}}\supp(\rho(\sigma))\) be the free part of the restriction of \(\rho\) to \(\mathfrak S_\infty\).
Proposition \ref{prop:actionss2} applied to \(\rho_{\restriction \mathfrak S_{2^\infty}}\) yields that the action induced by \(\rho_{\restriction \mathfrak S_{2^\infty}}\) is diagonally support-preserving factorizable.
Let \((\alpha_i)_{i\geq 0}\) such that the action	 \(\rho^{Y\setminus A_\infty}_{\restriction\mathfrak S_{2^\infty}}\)  induced on \(Y\setminus A_\infty\) by the restriction of \(\rho\) to \(\Stwo\) factorises in a support-preserving manner on the diagonal sum action \(\iota^{\odot (\alpha_i)_i}\) (see Definition \ref{df:symdiagsum}; the \(\alpha_i\)'s are well defined thanks to Proposition \ref{prop:diagfactdec}).

We now consider the associated character
\(\chi\colon \mathbb G\rightarrow \mathbb R,\)
defined by \(\chi(T)\coloneqq \nu(\Fix(\rho(T)))\). The definition of the diagonal sum action \(\iota^{\odot (\alpha_i)_i}\) and the fact that the factor map onto it is support-preserving imply that
\begin{equation}\label{equation:character}
	\chi(T)=\sum_{i\geq 0} \alpha_i\mu(\Fix(T))^i
\end{equation}
for every \(T\in\mathfrak S_{2^\infty}\). We will show in Proposition \ref{prop:fakechar} that the equation holds for every \(T\in \mathbb G\). 

\begin{remark}
	We do not know at this moment whether the character \(\chi\) is continuous, and hence we cannot apply the classification of characters on full groups \cite[Prop.\ 5.3]{carderiMorePolishfull2016}. Indeed \(\rho\) is a priori not uniform-to-uniform continuous and the map \(T\mapsto \mu(\Fix(T))\) is not continuous for the weak topology.
\end{remark}

We will first understand the free part of the action.

\subsection{The free part is the same as the free part of the restriction to \texorpdfstring{\(\mathfrak S_{2^\infty}\)}{S2infty}}\label{ssection:free}

Our aim now is to prove the following.

\begin{prop}\label{prop:equalpart}
	Given a boolean action \(\rho:\mathbb G\rightarrow \Aut(Y,\nu)\), we have that
	\[\bigcap_{\sigma\in\mathfrak S_{2^\infty}\setminus \{\id\}}\supp(\rho(\sigma))=\bigwedge_{T\in\mathbb G\setminus \{\id\}}\supp(\rho(T));\]
	that is, the free parts of \(\rho\) and of the restriction of \(\rho\) to \(\mathfrak S_{2^\infty}\) coincide.
\end{prop}

We first show an intermediate result.
The proof we present was suggested by the referee and simplifies nicely
our initial argument by removing the use of automatic continuity. Recall that in the previous section we defined  the free part of the restriction of \(\rho\) to \(\Stwo\) as \(A_\infty= \bigcap_{\sigma\in\mathfrak S_{2^\infty}\setminus \{\id\}}\supp(\rho(\sigma))\), which is clearly \(\rho(\Stwo)\)-invariant.

\begin{lem} \label{lem:Aisinvariant}
	The subset
	\(A_\infty\) is \(\rho(\mathbb G)\)-invariant.
\end{lem}
\begin{proof}
	Given any clopen subset \(O\) of \(X\), denote by \((\mathbb G)_O\)
	the subgroup of all \(T\in \mathbb G\) such that \(\supp T\subseteq O\).
	By Corollary \ref{cor: free part coincides with that of local subgroups},
	for any nonempty clopen subset \(O\),
	the free part of \(\rho\) restricted \((\mathfrak S_{2^\infty})_O\) is equal to \(A_\infty\).
	Since \((\mathfrak S_{2^\infty})_O\) commutes with \((\mathbb G)_{X\setminus O}\),
	we deduce that \(A_\infty\) is \(\rho((\mathbb G)_{X\setminus O})\)-invariant.
	Taking complements, we conclude that
	for all \(O\subseteq X\) proper clopen subset, \(A_\infty\) is
	\(\rho((\mathbb G)_{O})\)-invariant.
	The conclusion now follows by fixing two proper clopen subsets
	\(O_1\) and \(O_2\) such that \(O_1\cap O_2\neq\emptyset\) and \(X=O_1\cup O_2\): then
	\(A_\infty\) is both \(\rho((\mathbb G)_{O_1})\) and \(\rho((\mathbb G)_{O_2})\)-invariant,
	and \(\mathbb G=\la (\mathbb G)_{O_1},(\mathbb G)_{O_2}\ra\)
	by Lemma \ref{lem:generation in full groups}
	so that \(A_\infty\) is \(\rho(\mathbb G)\)-invariant as claimed.
\end{proof}

\begin{proof}[Proof of Proposition \ref{prop:equalpart}]
	We have to show
	\[\bigcap_{\sigma\in\mathfrak S_{2^\infty}\setminus \{\id\}}\supp(\rho(\sigma))=\bigwedge_{T\in\mathbb G\setminus \{\id\}}\supp(\rho(T));\]
	The left hand side is equal to \(A_\infty\)  and clearly contains the right hand side. By Lemma \ref{lem:Aisinvariant}, the set \(A_\infty\) is \(\mathbb G\)-invariant. Therefore we can assume that \(Y=A_\infty\) and need to prove that for every \(T\in\mathbb G\) the support of \(\rho(T)\) is the entire space \(X\), whenever \(T\neq \id\).

	Proposition \ref{prop:conjugacy of invol} tells us that
	every involution \(U\in\mathbb G\) whose support has measure \(1/2^n\)
	for some \(n\in\N\) is conjugate to an element of \(\mathfrak S_{2^\infty}\),
	and hence \(\rho(U)\) has support of full measure.
	We now argue as in Lemma \ref{lem:hyperfinite restriction free}: let \(T\in\mathbb G\) be a nontrivial element,
	we can find \(B\subseteq X\) measurable non null such that \(TB\cap B=\emptyset\).
	For large enough \(n\), we then find an involution \(U_n\) supported on \(B\)
	such that its support has measure \(1/2^n\).
	The commutator \([T,U_n]\) is an involution whose support has measure \(2/2^n\),
	and so \(\rho([T,U_n])\) has full measure support.

	But the support of \(\rho([T,U_n])\) is contained in \(\supp \rho(T)\cup \rho(U_n)(\supp \rho(T))\).
	By automatic continuity, Corollary\ \ref{cor: autocont action on malg}, the sequence \((\rho(U_n))_n\) tends weakly to the identity, therefore \(\rho(U_n)\supp T\) tends to \(\supp T\). This implies that \(Y=\supp \rho([T,U_n])\subseteq \supp \rho(T)\) as wanted.
\end{proof}

\subsection{Continuity and support dependency on the non-free part}

We will now show that \(\supp \rho(T)\) only depends on \(\supp T\). To this end, we first need to know that
\(\rho\) is uniform-to-uniform continuous, which is an easy consequence of what we have done so far.

\begin{prop}\label{prop:continuity}
	Given a boolean action \(\rho:\mathbb G\rightarrow \Aut(Y,\nu)\), we have that
	the action induced by \(\rho\) on its non-free part
	is uniform-to-uniform continuous.
\end{prop}
\begin{proof}
	By Proposition \ref{prop:actionss2}, the restriction \(\rho_{\restriction\mathfrak S_{2^\infty}}\)
	factors in a support preserving manner on a symmetric diagonal sum which by
	Proposition \ref{prop:TNFprod} is not discrete in the uniform topology.
	Since the factor map is support-preserving, this implies \(\rho\) has
	non-discrete image.
	We can now apply Corollary \ref{crl:automemb} to obtain that \(\rho\)
	is indeed uniform-to-uniform countionus.
\end{proof}

To see that \(\supp\rho(T)\) only depends on \(\supp T\), we will use again
the character associated to our boolean action. Since we now know that \(\rho\) is uniform-to-uniform continuous, \cite[Prop.\ 5.3]{carderiMorePolishfull2016} implies that Equation \eqref{equation:character} holds for every \(T\in \mathbb G\). However, there is a direct proof of this fact in the same line of that of \cite[Prop.\ 5.3]{carderiMorePolishfull2016}, and we now present it for the convenience of the reader.

\begin{prop}\label{prop:fakechar}
	Let \(\chi:\mathbb G\rightarrow \R\) be defined by \(\chi(T)=\nu(\Fix \rho(T))\).
	Then there is a unique sequence \((\alpha_i)_{0\leq i<\infty}\) such that
	for every \(T\in \mathbb G\) we have that
	\[\chi(T)=\sum_{i\geq 0}\alpha_i\mu(\Fix T)^i.\]
\end{prop}
\begin{proof}
	First observe that the uniqueness of \((\alpha_i)\) is a consequence of the uniqueness of the coefficients of the power series \(\sum_{i\geq 0}\alpha_i x^i\) whose radius of convergence is at least \(1\) and of the fact that \(\mu(\Fix(T))\) can take any value in \([0,1]\).

	Without loss of generality, we assume that the free part of \(\rho\) is trivial and by Proposition \ref{prop:equalpart},
	the free part of the restriction of \(\rho\) to \(\mathfrak S_{2^\infty}\) is also trivial.
	By Equation \eqref{equation:character} the desired result holds for elements of \(\mathfrak S_{2^\infty}\) (see also Prop.~\ref{prop:actionss2}).

	Fix \(T\in \mathbb G\). By Rokhlin's lemma and ergodicity, \(T\) is approximately conjugate to an element
	of the full group \([\mathcal R_0]\) and Proposition \ref{prop:closuredyadic} tells us that \([\mathcal R_0]=\overline{\mathfrak S_{2^\infty}}\).
	So there are sequences \(U_n\in \mathbb G\) and \(\sigma_n\in \mathfrak S_{2^\infty}\)
	such that \(U_n\sigma_nU_n^{-1}\) converges to \(T\).
	By uniform-to-uniform continuity (Proposition \ref{prop:continuity}), \(\rho(U_n\sigma_nU_n^{-1})\) converges to \(\rho(T)\) in the uniform topology and hence \[\nu(\Fix\rho(U_n\sigma_nU_n^{-1})\bigtriangleup \Fix\rho(T))\to 0.\]
	In particular, \(\nu(\Fix\rho(U_n\sigma_nU_n^{-1}))\to \nu(\Fix\rho(T))\).
	Now observe that
	\[
		\nu(\Fix \rho(U_n\sigma_nU_n^{-1}))=
		\nu(\Fix \rho(\sigma_n))=
		\sum_{i\geq 0}\alpha_i\mu(\Fix \sigma_n)^i=
		\sum_{i\geq 0}\alpha_i\mu(\Fix U_n\sigma_nU_n^{-1})^i.
	\]
	Clearly \(\mu(\Fix U_n\sigma_nU_n^{-1})^i\leq 1\), so using  \(\mu(\Fix\rho(U_n\sigma_nU_n^{-1}))\to \mu(\Fix\rho(T))\) we obtain that the right-hand term converges to \(\sum_i \alpha_i \mu(\Fix(\rho(T)))^i\). Since the left-hand term converges to \(\nu(\Fix\rho(T))=\chi(T)\), the conclusion follows.
\end{proof}

The exact formula in the previous proposition is not really important, all that
matters is that we now know that \(\nu(\supp \rho(T))\) depends only on \(\mu(\supp T)\). 
Moreover this dependency is continuous.

\begin{crl}\label{crl:rhoaper}
	If \(T\) is aperiodic on its support, then so is \(\rho(T)\).
\end{crl}
\begin{proof}
	Observe that a measure-preserving bijection \(T\) is aperiodic on its support if and only if
	for every \(n\), \(\mu(\Fix(T^n))=\mu(\Fix(T))\).
	Now if \(T\) is aperiodic on its support, then for every \(n\) we have
	\(\nu(\Fix(\rho(T)^n))=\chi(T^n)=\chi(T)=\nu(\Fix(\rho(T))\) so
	\(\rho(T)\) is aperiodic on its support.
\end{proof}

Here is another useful observation.
Note that \(\supp UT^n\subseteq\supp T\cup \supp U\).

\begin{lem}\label{lem:supputn}
	Let \(T\) and \(U\) be measure-preserving bijections,
	with \(T\) aperiodic when restricted to its support.
	Then \(\lim_n \left(\supp(UT^n)\Delta (\supp T\cup \supp U)\right)=0\).
\end{lem}
\begin{proof}
	Observe that if \(x\in\supp U\setminus \supp T\), then \(UT^n(x)\neq x\).
	For each \(n\in\N\), let \(A_n:=\{x\in \supp T: UT^n(x)=x\}\).
	Then \(x\in A_n\) implies \(T^n(x)=U\inv(x)\) and hence the aperiodicity of \(T\)
	implies that the measurable subsets \(\{A_n\}_n\) are pairwise disjoint.
	Hence \(\lim_n\mu(A_n)=0\) and the lemma is proved.
\end{proof}

We can now prove that \(\supp \rho(T)\) only depends on \(\supp T\).

\begin{prop}\label{prop:contofsupp}
	Let  \(T,U\in\mathbb G\) have the same support, then \(\rho(T)\) and \(\rho(U)\) also have the same support.
\end{prop}
\begin{proof}
	As a first step, let us assume that \(T\) is aperiodic on its support.
	By Lemma \ref{lem:supputn}, we have \[\lim_n\mu(\supp(UT^n))=\mu(\supp T\cup \supp U)=\mu(\supp T).\]
	Proposition \ref{prop:fakechar}  then
	implies \(\lim_n\nu(\supp(\rho(UT^n)))=\nu(\supp \rho(T))\).
	Corollary \ref{crl:rhoaper} tells us that \(\rho(T)\) is aperiodic on its support.
	We can therefore use a second time Lemma \ref{lem:supputn} to get
	\[\lim_n\nu(\supp\rho(UT^n))=\nu(\supp \rho(T)\cup \supp\rho(U)).\]
	Hence \(\nu(\supp\rho(U)\setminus \supp\rho(T))=0\).
	On the other hand, using Proposition \ref{prop:fakechar} again
	we have \(\nu(\supp\rho(U))=\nu(\supp\rho(T))\)
	and hence \(\supp \rho(U)=\supp \rho(T)\).

	Now if \(T\) is not aperiodic on its support, we can find an  element
	\(T'\in\mathbb G\) which is aperiodic on its support and has same support as \(T\) and \(U\).
	The above argument then yields \(\supp \rho(T)=\supp \rho(T')=\supp \rho(U)\) as wanted.
\end{proof}

\subsection{End of the proof of Theorem \ref{thm:main}}

As we already mentioned at the beginning of Section \ref{sec:proof of main thm}, the uniqueness part in Theorem \ref{thm:main} is a direct consequence of Proposition \ref{prop:diagfactdec} so in order to finish the proof of Theorem \ref{thm:main}, we  have to show that Theorem \ref{thm:diagonally factorizable} holds, namely that every boolean \(\mathbb G\)-action is diagonally support-preserving factorizable.

Let \(\mathbb G\leq \Aut(X,\mu)\) be an ergodic full group and fix a boolean action \(\rho\colon\mathbb G\rightarrow \Aut(Y,\nu)\).
We can assume that the free part of the action of \(\mathbb G\) on \((Y,\nu)\) is trivial and by Proposition \ref{prop:equalpart} this implies that the restriction of \(\rho\)
to \(\mathfrak S_{2^\infty}\) also has trivial free part.
Proposition \ref{prop:actionss2}
provides a symmetric diagonal sum boolean action \(\iota^{\odot (\alpha_i)_i}\) of \(\mathfrak S_{2^\infty}\) on
\((X^{\odot (\alpha_i)_i},\mu^{\odot (\alpha_i)_i})\), and a support-preserving factor map \(\pi\colon Y\rightarrow Z\)
of \(\rho_{\restriction \mathfrak S_{2^\infty}}\) onto \( \iota^{\odot (\alpha_i)_i}\).

We will show that \(\pi\) is still a measure preserving support preserving factor map of \(\rho\)
onto the same diagonal sum action \(\iota^{\odot (\alpha_i)_i}\), now seen as a boolean action of \(\mathbb G\).

Given \(\sigma\in\mathfrak S_{2^\infty}\), we already know that \(\pi^{-1}(\supp\iota^{\odot (\alpha_i)_i}(\sigma))=\supp\rho(\sigma)\) and by uniform-to-uniform continuity (see Proposition \ref{prop:continuity}), the same is true for elements of \([\mathcal R_0]=\overline{\mathfrak S_{2^\infty}}\), cf.\ Proposition \ref{prop:closuredyadic}. Fix \(T\in\mathbb G\). By Corollary \ref{cor:every set is some support},
there exists \(T'\in[\mathcal R_0]\) such that \(\supp(T')=\supp(T)\) and hence by Proposition \ref{prop:contofsupp} we have
\[
	\pi^{-1}(\supp \iota^{\odot (\alpha_i)_i}(T))=\pi^{-1}(\supp \iota^{\odot (\alpha_i)_i}(T'))=\supp \rho(T')=\supp\rho(T).
\]
We conclude that \(\pi\) is indeed support-preserving. Let us now show it is a factor map. Take \(T,U\in \mathbb G\) then we have
\begin{align*}
	\pi^{-1}(\iota^{\odot (\alpha_i)_i}(T)\supp \iota^{\odot (\alpha_i)_i}(U)) & =\pi^{-1}(\supp \iota^{\odot (\alpha_i)_i}(TUT\inv)) \\
	                                                                           & = \supp \rho(TUT\inv)                                \\
	                                                                           & =\rho(T)\supp \rho(U),
\end{align*}
so \(\pi\) is equivariant on elements of the form \(\supp \iota^{\odot (\alpha_i)_i}(U)\). The action \(\iota^{\odot (\alpha_i)_i}\) of \(\mathbb G\) is totally non free, Proposition \ref{prop:TNFprod}, so by definition these sets generate the measure algebra of \((X^{\odot (\alpha_i)_i},\mu^{\odot (\alpha_i)_i})\). Therefore \(\pi\) is equivariant and the proof of Theorem \ref{thm:diagonally factorizable} is concluded.

\subsection{Recovering Dye's reconstruction theorem}\label{sec:dyereconstruction}

We can now recover Dye's reconstruction theorem in the ergodic case \cite[Thm.~2]{dyeGroupsMeasurePreserving1963}.

\begin{thm}\label{thm:dyereconstructed}
	Let \(\mathbb G\) and \(\mathbb H\) be two ergodic full groups on \((X,\mu)\), let \(\varphi:\mathbb G\to\mathbb H\) be a group isomorphism.
	Then there is a unique \(S\in \Aut(X,\mu)\) such that for all \(T\in \mathbb G\),
	\[\varphi(T)=STS\inv.\]
\end{thm}
\begin{proof}
	Let us start by proving uniqueness: suppose that \(S\in\Aut(X,\mu)\) satisfies that for all \(T\in \mathbb G\),
	\(\varphi(T)=STS\inv\). Let \(A\in\MAlg(X,\mu)\), then by Corollary \ref{cor:every set is some support} there is \(T\in \mathbb G\) whose support is equal to \(A\), and then
	\(\supp \varphi(T)=\supp STS\inv=S(\supp T)=S(A)\) so we see that \(S(A)\) is completely determined by \(\varphi\), which implies the uniqueness of \(S\).

	We now prove the existence.
	Let \(X_0,\dots,X_\infty\) be the sets arising from the decomposition of \(\varphi\) as a diagonal support-preserving factorizable boolean action from Theorem \ref{thm:diagonally factorizable}.
	Since \(\mathbb H\) is ergodic, \(\varphi\) is an ergodic boolean action, so one of these sets is equal to \(X\).
	Since \(\mathbb H\) contains elements of arbitrary support and \(\varphi\) is an isomorphism, both \(X_0\) and \(X_\infty\) must be empty.

	We thus have \(n\geq 1\) such that \(\varphi\) factors in a support-preserving manner
	onto the symmetric diagonal action \(\iota_n\) of \(\mathbb G\) on \(X^{\odot n}\).
	Denote by \(\pi_n\) the associated (unique) support-preserving factor map.
	Observe that if \(n\geq 2\), the supports of elements of \(\iota_n(\mathbb G)\)
	form a proper subset of \(\MAlg(X^{\odot n},\mu^{\odot n})\).
	Moreover, for every \(h\in\mathbb H\), if we let \(g=\varphi\inv(h)\) we have
	\[\supp h=\supp \varphi(g)=\pi_n\inv(\supp \iota_n(g)).\]
	Since \(\pi_n\inv\) is an embedding of measure algebras, we conclude that
	the supports of the elements of \(\mathbb H\) form a proper
	subset of \(\MAlg(X,\mu)\), a contradition.

	We thus have \(n=1\).
	Again every \(A\in\MAlg(X,\mu)\) is the support of some \(h\in\mathbb H\), and hence of the form
	\(\pi_1\inv(\supp g)\) for some \(g\in\mathbb G\),
	so \(\pi_1\inv\) is a surjective measure-preserving
	measure algebra homomorphism
	\(\MAlg(X,\mu)\to\MAlg(X,\mu)\).
	We conclude that \(\pi_1\in\Aut(X,\mu)\),
	in other words \(\pi_1\) is an isomorphism between \(\varphi\)
	and the inclusion \(\iota_1\) of \(\mathbb G\) in \(\Aut(X,\mu)\).
	This means that
	for all \(T\in \mathbb G\) and almost all \(x\in X\)
	\(\pi_1 \varphi(T) x=\iota_1(T)\pi_1(x)\).
	If we let \(S=\pi_1\inv\), we finally get \(\varphi(T)=S T S\inv\) as wanted.
\end{proof}

\section{Boolean actions of full groups of p.m.p.\ equivalence relations}\label{sec: applications}

We will now explain how Theorem \ref{thm: main - in section} implies that for every measure-preserving
boolean action of a full group \([\mathcal R]\),
the non-free part of the action comes from measure-preserving actions
of the equivalence relation and its symmetric tensor powers, yielding a
statement more similar to the main result of Matte Bon \cite{mattebonRigidityPropertiesFull2018}. The rest of the section is then devoted to applications of this result.

\begin{thm}\label{thm: main pmp equivalence relations}
	Let \(\mathcal R\) be a \pmp ergodic equivalence relation on \((X,\mu)\). Let \(\rho\colon[\mathcal R]\to \Aut(Y,\nu)\) be a boolean action.
	Then there is a unique partition \(Y=\bigsqcup_{n\in\N} Y^\rho_n\sqcup Y^\rho_\infty\) into (possibly null) \(\rho([\mathcal R])\)-invariant sets such that:
	\begin{itemize}
		\item the restriction of \(\rho\) to \(Y^\rho_0\) is the trivial action;
		\item for all \(n\geq 1\) such that \(\nu(Y_n^\rho)>0\), there is a \pmp equivalence relation \(\mathcal S_n\) on \(Y_n\) and a class-bijective factor \(\pi\colon \mathcal S_n\rightarrow \mathcal R^{\odot n}\) such
		      that
		      \[
			      \rho_{\restriction Y^\rho_n}=\rho_\pi\circ\iota_n ,
		      \]
		      where \(\iota_n: [\mathcal R]\to[\mathcal R^{\odot n}]\) is the diagonal embedding and \(\rho_\pi\) is obtained \(\pi\) as defined in Lemma \ref{lemma: def rho pi};
		\item the restriction of \(\rho\) to \(Y^\rho_\infty\) is free.
	\end{itemize}
\end{thm}
\begin{proof}
	We know from Theorem \ref{thm: main - in section} that \(\rho\) is diagonally support-preserving factorizable, which yields the unique partition \(Y=\bigsqcup_{n\in\N} Y^\rho_n\sqcup Y^\rho_\infty\) that we seek. The two first items are satisfied by definition, let us see why the third holds. We have that \(\rho\) restricted to \(Y^\rho_n\) factors onto \(\iota^{\odot n}:[\mathcal R]\to [\mathcal R^{\odot n}]\) in a support-preserving manner.

	Let \(n\geq 1\), assume that \(\nu(Y_n^\rho)>0\) and denote by \(\nu_{Y_n^\rho}\) the probability measure on \(Y_n^\rho\) obtained by renormalizing the restriction of \(\nu\) to \(Y_n^\rho\).
	Let us define \(\rho_n:\iota^{\odot n}([\mathcal R])\to \Aut(Y^\rho_n,\nu_{Y_n^\rho})\) by \(\rho_n(\iota^{\odot n}(T))=\rho(T)\).
	By Lemma \ref{lem: fg of R gens fg of Rn} the full group generated by \(\iota^{\odot n}([\mathcal R])\) is equal to \([\mathcal R^{\odot n}]\) so by Proposition \ref{prop: support preserving extend to full group} we may extend \(\rho'\) to a support-preserving boolean action \(\rho':[\mathcal R^{\odot n}]\to\Aut(Y_n^\rho,\nu_{Y_n^\rho})\).
	It now follows from Lemma \ref{lem:pmp action from support pres} that there is a measure-preserving action \(\pi:\mathcal S_n\to\mathcal R^{\odot n}\) of  \(\mathcal R^{\odot n}\) on \((Y,\nu)\) for some measure-preserving equivalence relation \(\mathcal S_n\) on \((Y,\nu)\) such that \(\rho'=\rho_\pi\).
	Since \(\rho'(\iota^{\odot n}(T))=\rho_{\restriction Y_n}(T)\) for all \(T\in [\mathcal R]\), this finishes the proof.
\end{proof}

Let \(\mathcal R\) be a \pmp ergodic equivalence relation on \((X,\mu)\) and let \(\rho\colon[\mathcal R]\to \Aut(Y,\nu)\) be a boolean action. Denote by \(\rho_\infty\) the restriction of \(\rho\) to the measurable subset \(Y^\rho_\infty\) given by Theorem \ref{thm: main pmp equivalence relations}. Then \(\rho_\infty([\mathcal R])\leq \Aut(Y_{\infty}^\rho,\nu_{Y_\infty^\rho})\) is discrete. In particular, if \(\rho([\mathcal R])\) is separable with respect to the uniform topology, we must have that \(Y_\infty^\rho\) has measure zero. We thus get the following corollary.

\begin{corollary}\label{cor: main thm pmp to pmp}
	Let \(\mathcal R\) and \(\mathcal S\) be \pmp equivalence relations on \((X,\mu)\) and \((Y,\nu)\) respectively, suppose that \(\mathcal R\) is ergodic. Let \(\rho\colon[\mathcal R]\to [\mathcal S]\) be a boolean action. Then there is a unique partition \(Y=\bigsqcup_{n\in\N} Y_n^\rho\) into \(\rho([\mathcal R])\)-invariant sets such that:
	\begin{itemize}
		\item the restriction of \(\rho\) to \(Y_0^\rho\) is the trivial action;
		\item for all \(n\geq 1\), there is a \pmp equivalence relation \(\mathcal S_n\subset \mathcal S\) on \(Y_n^\rho\) and a class-bijective factor \(\pi\colon \mathcal S_n\rightarrow \mathcal R^{\odot n}\) such
		      that
		      \[
			      \rho_{\restriction Y_n}=\rho_\pi\circ\iota_n ,
		      \]
		      where \(\iota_n: [\mathcal R]\to[\mathcal R^{\odot n}]\) is the diagonal embedding and \(\rho_\pi\) is obtained \(\pi\) as defined in Lemma \ref{lemma: def rho pi}.\qed
	\end{itemize}
\end{corollary}

\begin{df}
	Let \(\mathcal R\)  and \(\mathcal S\) be \pmp ergodic equivalence relations on \((X,\mu)\). Let \(\rho\colon [\mathcal R]\rightarrow [\mathcal S]\) be a boolean action. We say that \(\rho\) is \defin{extension-like} when the set \(Y_0^\rho\sqcup Y^\rho_1\) from Corollary \ref{cor: main thm pmp to pmp} has full measure.
\end{df}

In other words, a boolean action \(\rho:[\mathcal R]\to[\mathcal S]\) is extension-like if once we remove its trivial part \(Y_0^\rho\), it comes from an extension \(\pi:Y\setminus Y_0^\rho\to X\) of \(\mathcal R\) by subequivalence relation \(\mathcal S'\) of the restriction of \(\mathcal S\) to \(Y\setminus Y_0^\rho\).

\subsection{Amenability}

A \pmp equivalence relation \(\mathcal R\) on \((X,\mu)\) is called \textbf{amenable} when up to restricting it to a full measure subset, we can find a sequence of Borel maps \(f^n: \mathcal R\to\R^{\geq 0}\) such that for all \(x\in X\), \(\sum_{x'\in[x]_{\mathcal R}} f_n(x,x')=1\), and if we denote by \(\norm{\cdot}_1\) the \(\ell^1\) norm on \(\ell^1([x]_{\mathcal R})\), then for all \(y,z\in [x]_{\mathcal R}\), we have
\[\norm{f_n(y,\cdot)-f_n(z,\cdot)}_1\to 0.\]
In other words, we have a measurable way to assign to every \(x\in X\) a sequence of probability measures on \([x]_{\mathcal R}\) so that the probability measures assigned to points \(y\) and \(z\) in the same equivalence class are asymptotically the same. Such a sequence \((f_n)\) is called a \textbf{Reiter sequence} for \(\mathcal R\).

Our main result (Corollary \ref{cor: main thm pmp to pmp}) interacts with amenability as follows.

\begin{prop}\label{prop: S amenable}
	Let \(\mathcal R\)  and \(\mathcal S\) be \pmp equivalence relations on \((X,\mu)\) and \((Y,\nu)\) respectively with \(\mathcal R\) ergodic. Let us assume that there is a non-trivial group homomorphism \(\rho\colon [\mathcal R]\rightarrow [\mathcal S]\). If \(\mathcal S\) is amenable, then \(\mathcal R\) is.
\end{prop}

In order to prove this proposition, we need a few well-known lemmas.

\begin{lemma}
	Every subequivalence relation of an amenable \pmp equivalence relation is amenable.
\end{lemma}
\begin{proof}[Sketch of proof]
	Since every subequivalence relation of a hyperfinite equivalence relation has to be hyperfinite, a quick way to see this is to use the Connes-Feldman-Weiss theorem \cite{connesAmenableEquivalenceRelation1981} and the fact that hyperfinite equivalence relations are amenable. See \cite[Sec.~6 and 10]{kechrisTopicsOrbitEquivalence2004} for details.
\end{proof}
The following lemma is the \pmp equivalence relation version of the fact that if a countable group \(\Gamma\) is acting freely \pmp on \((X,\mu)\), then the associated equivalence relation \(\mathcal R_\Gamma\) is amenable iff \(\Gamma\) is amenable.

\begin{lemma}\label{lem: amenable iff ext is}
	Let \(\mathcal R\)  and \(\mathcal S\) be \pmp equivalence relations on \((X,\mu)\) and \((Y,\nu)\) respectively. Let \(\pi\colon Y\rightarrow X\) be a class-bijective factor. Then \(\mathcal S\) is amenable if and only if \(\mathcal R\) is.
\end{lemma}
\begin{proof}
	If \(\mathcal R\) is amenable, let \((f_n)\) be a Reiter sequence witnessing it.
	We build a Reiter sequence \((g_n)\) for \(\mathcal S\) by letting \(g_n(y,y')=f_n(\pi(y),\pi(y'))\).
	Since for every \(y\in Y\), \(\pi\) induces a bijection between \([y]_{\mathcal S}\) and \([\pi(y)]_{\mathcal R}\), it is not hard to check that \((g_n)\) is indeed a Reiter sequence.

	Conversely, let \((g_n)\) be a Reiter sequence for \(\mathcal S\). For all \((x,x')\in\mathcal R\) and \(y\in Y\) such that \(\pi(y)=x\), let \(\pi_x\inv(x')\) denote the unique \(y'\in Y\) such that \((y,y')\in\mathcal S\) and \(\pi(y')=x'\). Let \((\nu_x)_{x\in X}\) be the disintegration of \(\nu\) w.r.t.\ \(\pi\).   We then  let \(f_n(x,x')=\int_{\pi\inv(x)} g_n(y,\pi_x\inv(x'))d\nu_x(y)\).
	Observe that \[\norm{f_n(x,\cdot)-f_n(x',\cdot)}_1\leq\int_{\pi\inv(x)}\norm{g_n(y,\cdot)-g_n(\pi_x\inv(x'),\cdot)}_1d\nu_x(y).\]
	Since \((g_n)\) is Reiter, the integrand converges to zero pointwise, and it is bounded by \(2\) so
	by the Lebesgue dominated convergence we conclude that \((f_n)\) is Reiter as wanted.
\end{proof}

\begin{lemma}\label{lem: R amenable iff Rn is}
	Let \(\mathcal R\) be a \pmp ergodic equivalence relations on \((X,\mu)\). Then the following are equivalent:
	\begin{enumerate}[(i)]
		\item \label{cond: R amenable}\(\mathcal R\) is amenable;
		\item \label{cond: Rodotn am for all n}\(\mathcal R^{\odot n}\) is amenable for all \(n\geq 1\);
		\item \label{cond: Rodotn am for some n}\(\mathcal R^{\odot n}\) is amenable for some \(n\geq 1\).
	\end{enumerate}
\end{lemma}
\begin{proof}
	Since \(\mathcal R^n\) is a class-bijective extension of \(\mathcal R^{\odot n}\), by Lemma \ref{lem: amenable iff ext is} it suffices to show the same equivalences as above where \(\mathcal R^{\odot n}\) is replaced by \(\mathcal R^n\). Also note that the implication \eqref{cond: Rodotn am for all n}\(\impl\)\eqref{cond: Rodotn am for some n} is clear.

	Let us show \eqref{cond: R amenable}\(\impl\)\eqref{cond: Rodotn am for some n}. Assume that \(\mathcal R\) is amenable and let \((f_k)\) be a Reiter sequence, then
	\[((x_1,...,x_n),(y_1,...,y_n))\mapsto f_k(x_1,y_1)\cdots f_k(x_n,y_n)\] is easily seen to define a Reiter sequence for \(\mathcal R^n\).

	We finally show \eqref{cond: Rodotn am for some n}\(\impl\)\eqref{cond: R amenable}.
	Suppose \(\mathcal R^n\) is amenable. Observe that the projection on the first coordinate \(\pi:X^n\to X\) is a class-bijective factor map from the subequivalence \(\mathcal S\) of \(\mathcal R^n\) given by \((((x_1,...,x_n),(y_1,...,y_n))\in\mathcal S\) iff \((x_1,y_1)\in\mathcal R\) and \(x_i=y_i\) for all \(i\in\{2,\dots,n\}\).
	Since every subequivalence relation of an amenable equivalence relation is amenable,
	\(\mathcal S\) is amenable, so one more application of Lemma \ref{lem: amenable iff ext is} yields that \(\mathcal R\) is amenable as wanted.
\end{proof}

\begin{proof}[Proof of Proposition \ref{prop: S amenable}]
	Since \(\rho\) is non-trivial, Corollary \ref{cor: main thm pmp to pmp} yields the existence of \(n\geq 1\) such that some subequivalence relation of \(\mathcal S\) is a class-bijective extension of \(\mathcal R^{\odot n}\).  Since \(\mathcal S\) is amenable, we get that \(\mathcal R^{\odot n}\) is amenable and hence by Lemma \ref{lem: R amenable iff Rn is}, \(\mathcal R\) is amenable as well.
\end{proof}

\begin{remark}\label{rmk: haagerup vs T}
	One can restate Proposition \ref{prop: S amenable} as the fact that as soon as \(\mathcal R\) is ergodic nonamenable and \(\mathcal S\) is amenable, there are no non-trivial group homomorphisms \([\mathcal R]\to[\mathcal S]\).
	One can probably also rule out homomorphisms \([\mathcal R]\to[\mathcal S]\) when \(\mathcal R\) has (T) and \(\mathcal S\) has the Haagerup property.
	However proving this would require to develop the theory of property (T) and the Haagerup property for \pmp equivalence relations in the non-ergodic context, which is beyond the scope of this paper.
\end{remark}
\subsection{Cost and treeability}

Given a \pmp equivalence relation \(\mathcal R\), a \textbf{graphing} of \(\mathcal R\) is a Borel subset \(\mathcal G\subseteq \mathcal R\) which is symmetric: for all \((x,y)\in\mathcal G\), we have \((y,x)\in\mathcal G\). Given such a graphing, the \textbf{\(\mathcal G\)-degree} of \(x\in X\) is \(\deg_{\mathcal G}(x)\coloneqq \abs{\{y\in X\colon (x,y)\in\mathcal G\}}\). The \textbf{cost} of a graphing \(\mathcal G\) is
\[\Cost(\mathcal G)\coloneqq\frac 12\int_X\deg_{\mathcal G}(x).\]
A graphing of \(\mathcal R\) \textbf{generates} \(\mathcal R\) when \(\mathcal R\) is the smallest equivalence relation containing \(\mathcal G\), and the \textbf{cost} of \(\mathcal R\) is the infimum of the costs of the graphings which generate it.

A \textbf{treeing} of \(\mathcal R\) is a graphing of \(\mathcal R\) whose connected components are trees (have no circuits).
A \pmp equivalence relation \(\mathcal R\) is \textbf{treeable} when it can be generated by a treeing. A fundamental theorem of Gaboriau shows that the cost of \(\mathcal R\) has then to be equal to the cost of such a treeing \cite[Thm.~1]{gaboriauCoutRelationsEquivalence2000}.

We now explain how Theorem \ref{thm: main pmp equivalence relations} yields restrictions for homomorphisms between full groups of \pmp ergodic equivalence relations related to the above definitions. Here is the result.

\begin{prop}\label{prop: S treeable}
	Let \(\mathcal R\)  and \(\mathcal S\) be non-amenable \pmp equivalence relations on \((X,\mu)\) and \((Y,\nu)\) respectively, with \(\mathcal R\) ergodic. Let us assume that there is a non-trivial group homomorphism \(\rho\colon [\mathcal R]\rightarrow [\mathcal S]\). If \(\mathcal S\) is treeable, then the cost of \(\mathcal R\) is strictly greater than \(1\) and \(\rho\) is extension-like.
\end{prop}

\begin{lemma}
	Let \(\mathcal R\) be a \pmp aperiodic equivalence relation on \((X,\mu)\). Then for every \(n\geq 2\), the cost of \(\mathcal R^{\odot n}\) is \(1\).
\end{lemma}
\begin{proof}
	Since \(n\geq 2\), the equivalence relation \(\mathcal R^n\) has cost \(1\) by
	\cite[Thm.\ 24.9]{kechrisTopicsOrbitEquivalence2004}. 
	Moreover \(\mathcal R^n\) has index \(n!\) in
	\[\mathcal R^n\rtimes \mathfrak S_n\coloneqq
		\{((x_i)_{i=1}^n,(y_i)_{i=1}^n)\colon \exists \sigma\in\mathfrak S_n
		\forall i\in\{1,\dots,n\}, (x_i,y_{\sigma(i)})\in\mathcal R\},\]
	so we deduce from Gaboriau's result on cost for direct products \cite[Prop.~V.1]{gaboriauCoutRelationsEquivalence2000} (see also \cite[Prop.~25.1]{kechrisTopicsOrbitEquivalence2004}) that \(\mathcal R^n\rtimes \mathfrak S_n\) has cost at most \(1\). But \(\mathcal R^n\rtimes \mathfrak S_n\) has cost at least \(1\) by aperiodicity (this dates back to Levitt's paper where he introduces cost \cite[Thm.~2]{levittCostGeneratingEquivalence1995}, see also \cite[Cor.~21.3]{kechrisTopicsOrbitEquivalence2004}) so it has cost \(1\).
	Since \(\mathcal R^{\odot n}\) and \(\mathcal R^n\rtimes\mathfrak S_n\) are stably orbit equivalent,
	Gaboriau's induction formula for stable orbit equivalence \cite[Invariance~II.11.]{gaboriauCoutRelationsEquivalence2000} (see also \cite[Thm.~21.1]{kechrisTopicsOrbitEquivalence2004}) yield that \(\mathcal R^{\odot n}\) has cost \(1\) as wanted.
\end{proof}

\begin{proof}[Proof of Proposition \ref{prop: S treeable}]
	Suppose by contradiction that the boolean action \(\rho:[\mathcal R]\to[\mathcal S]\) is not extension-like, then since it is not trivial we must have
	that \(Y_n^\rho\) has positive measure for some finite \(n\geq 2\) as per Corollary \ref{cor: main thm pmp to pmp}.
	In particular, the restriction of \(\mathcal S\) to \(Y_n^\rho\) contains a class-bijective extension \(\mathcal S_n\) of \(\mathcal R^{\odot n}\). Since \(\mathcal S\) is treeable, so is its restriction to \(Y_n^\rho\) by \cite[Prop.~II.10.]{gaboriauCoutRelationsEquivalence2000}.

	Moreover, \(\mathcal S_n\) must have cost at most \(1\) since  it is a class-bijective extension of a cost \(1\) \pmp equivalence relation (see the second remark after Proposition 10.14 in \cite{kechrisGlobalaspectsergodic2010}), so it has cost \(1\) by aperiodicity of \(\mathcal S_n\) (which itself follows from the aperiodicity of its factor \(\mathcal R^{\odot n}\)). Also note that \(\mathcal S_n\) is treeable by \cite[Thm.~IV.4]{gaboriauCoutRelationsEquivalence2000}, and not amenable since it is a class-bijective extension of a the non-amenable \pmp equivalence relation \(\mathcal R^{\odot n}\). We arrive at the desired contradiction by \cite[Cor.~IV.2]{gaboriauCoutRelationsEquivalence2000}.

	So our action is extension-like, and we thus get that the restriction of \(\mathcal S\) to \(Y^\rho_1\) is a class-bijective extension of \(\mathcal R\). Again if the cost of \(\mathcal R\) were equal to \(1\), we would get that the treeable non-amenable equivalence relation \(\mathcal S\cap (Y^\rho_1\times Y^\rho_1)\) has cost \(1\), a contradiction. So \(\mathcal R\) has cost strictly greater than \(1\) as wanted.
\end{proof}

\begin{remark}
	In the context of the above proposition, we do not know whether if \(\mathcal S\) is treeable, then \(\mathcal R\) has to be treeable as well. Indeed it is unknown whether there are non-treeable equivalence relations with treeable actions. Our result only exploits the fact that such actions do not exist when the cost of the acting equivalence relation is \(1\).
\end{remark}

\subsection{Solid ergodicity}

The notion of solid ergodicity derives from solidity, a von Neumann algebraic notion
with strong structural consequences that was discovered by Ozawa \cite{ozawaSolidNeumannAlgebras2004},
answering a question of Ge about free group factors
that can be found at the end of \cite{geApplicationsFreeEntropy1998}.
While solidity asks that all diffuse subalgebras have amenable relative commutant,
solid ergodicity only makes sense for von Neumann algebras of \pmp equivalence relations,
and asks that all diffuse subalgebras of \(\LL^\infty(X,\mu)\)
have amenable relative commutant.

In particular, all amenable equivalence relations are solidly ergodic (regardless of their ergodicity!),
but more interestingly Bernoulli shifts of non-amenable groups yield
solidly ergodic equivalence relations, as was discovered
by Chifan and Ioana \cite{chifanErgodicSubequivalenceRelations2010}.
Their result was later on subsumed by the work of Boutonnet on Gaussian actions
\cite{boutonnetSolidErgodicityGaussian2012}.
Chifan and Ioana obtained the following characterization of solid ergodicity
which we take as a definition
(see \cite[Prop.~6]{chifanErgodicSubequivalenceRelations2010}).
The terminology ``solidly ergodic'' was proposed by Gaboriau 
shortly after their paper, in Section~5 of the survey 
\cite{gaboriauOrbitEquivalenceMeasured2011}.

\begin{df}
	A \pmp equivalence relation \(\mathcal R\) is called \textbf{solidly ergodic} if for any subequivalence relation \(\mathcal R'\subseteq\mathcal R\), there is a partition \(X=\bigsqcup_{n\geq 0}X_n\) of \(X\) into \(\mathcal R'\)-invariant sets such that the restriction of \(\mathcal R'\) to \(X_0\) is amenable, and for every \(n\geq 1\) the restriction of \(\mathcal R'\) to \(X_n\) is ergodic and non-amenable.
\end{df}
Note that by \cite[Prop.~6]{chifanErgodicSubequivalenceRelations2010}, one could equivalently ask that for every \(n\geq 1\) the restriction of \(\mathcal R'\) to \(X_n\) is strongly ergodic.
\begin{prop}\label{prop: easier def solid}
	A \pmp equivalence relation \(\mathcal R\) is solidly ergodic iff for every non-amenable subequivalence relation \(\mathcal R'\), there is an \(\mathcal R'\)-invariant set \(A\subseteq X\) of positive measure such that the restriction of \(\mathcal R'\) to \(A\) is ergodic.
\end{prop}
\begin{proof}
	The direct implication is clear. Conversely,
	let \(\mathcal R'\) be any subequivalence relation of \(\mathcal R\), let \(X_0\subseteq X\) measurable and maximal among measurable sets \(Y\) such that the restriction of \(\mathcal R'\) to \(Y\) is amenable. Observe that \(X_0\) is \(\mathcal R'\)-invariant. By maximality the restriction of \(\mathcal R'\) to \(X\setminus X_0\) is everywhere non-amenable, which by hypothesis yields that its ergodic decomposition is discrete as wanted.
\end{proof}

Our main result implies the following rigidity result for homomorphisms taking values into full groups of solidly ergodic equivalence relations.

\begin{prop}\label{prop: S solid}
	Let \(\mathcal R\)  and \(\mathcal S\) be non-amenable \pmp equivalence relations on \((X,\mu)\)
	and \((Y,\nu)\) respectively, with \(\mathcal R\) ergodic.
	Let us assume that there is a non-trivial group homomorphism
	\(\rho\colon [\mathcal R]\rightarrow [\mathcal S]\).
	If \(\mathcal S\) is solidly ergodic, then \(\mathcal R\) is solidly ergodic as well
	and \(\rho\) is extension-like.
\end{prop}

\begin{lemma}\label{lem: factor of solid is solid}
	Let \(\mathcal R\) and \(\mathcal S\) be non-amenable \pmp ergodic equivalence relations on \((X,\mu)\) and \((Y,\nu)\) respectively. Let \(\pi\colon Y\rightarrow X\) be a class-bijective factor. If \(\mathcal S\) is solidly ergodic, then \(\mathcal R\) is.
\end{lemma}
\begin{proof}
	We use the characterization provided by Proposition \ref{prop: easier def solid}.
	Consider a non-amenable sub-equivalence relation \(\mathcal R'\subseteq \mathcal R\).
	It yields a sub-equivalence relation \(\mathcal S'\subseteq \mathcal S\) given by \((x,y)\in\mathcal S'\) if \((x,y)\in\mathcal S\) and \((\pi(x),\pi(y))\in\mathcal R'\).

	By construction \(\pi\) is a class-bijective factor map from \(\mathcal S'\) to \(\mathcal R'\). Lemma \ref{lem: amenable iff ext is} yields that \(\mathcal S'\) is not amenable. 
	Since \(\mathcal S\) is solidly ergodic, there is a measurable \(\mathcal S'\)-invariant subset \(A\subseteq Y\) such that \(\mathcal S'\) restricted to \(A\) is ergodic. 
	Let \(B\subseteq X\) be the essential image of \(A\) (obtained as the preimage 
	of the $\pi$-measurable envelope of $A$ in $Y$, see \cite[Thm.~29.12]{kechrisClassicaldescriptiveset1995}), \emph{i.e.}, 
	the unique up to measure zero measurable set such that
	\begin{enumerate}
		\item \(\nu(\pi^{-1}(B)\cap A)=\nu(A)\);
		\item for every positive measure subset \(B'\subseteq B\) one has that \(\nu(\pi^{-1}(B')\cap A)>0\).
	\end{enumerate}
	Then \(\mathcal R'\) restricted to \(B\) is ergodic: if \(B'\subseteq B\) is \(\mathcal R'\)-invariant, then \(\pi^{-1}(B')\cap A\) is \(\mathcal S'\)-invariant, so it must be either empty or contain \(A\), which implies that \(B'=\emptyset\) or \(B'=B\) as wanted. So we have shown that \(\mathcal R'\) can be restricted to an ergodic equivalence relation,
	which by Proposition \ref{prop: S solid} shows that \(\mathcal R\) is solidly ergodic as wanted.
\end{proof}

\begin{proof}[Proof of Proposition \ref{prop: S solid}]
	We begin by showing that for all \(n\geq 2\), the \pmp equivalence relation \(\mathcal R ^{\odot n}\) is not solid. 
	Indeed \(\mathcal R^n\) is not solidly ergodic because it contains the non-amenable equivalence relation \(\id_X\times \mathcal R^{n-1}\) which has diffuse ergodic decomposition. 
	Since \(\mathcal R^n\) is a class-bijective extension of \(\mathcal R^{\odot n}\), we conclude by the previous lemma that \(\mathcal R^{\odot n}\) is not solidly ergodic.

	Now, let \(\rho:[\mathcal R]\to[\mathcal S]\) be a non-trivial group homomorphism.
	Since \(\mathcal S\) is solidly ergodic, all its restrictions are solidly ergodic as well, and combining the non solid ergodicity of \(\mathcal R^{\odot n}\) for \(n\geq 2\) and Corollary \ref{cor: main thm pmp to pmp}, we see that \(\rho\) must be extension-like.
	Since \(\rho\) is not trivial and any restriction of \(\mathcal S\) has to be solidly ergodic as well, it then follows from Lemma \ref{lem: factor of solid is solid} that \(\mathcal R\) is solid as wanted.
\end{proof}

\subsection{Conclusions}

Let us recap the restrictions for non-trivial homomorphisms \(\rho:[\mathcal R]\to [\mathcal S]\) that we have obtained when \(\mathcal R\) is ergodic:
\begin{itemize}
	\item If \(\mathcal S\) is amenable, \(\mathcal R\) must be amenable as well;
	\item If both \(\mathcal R\) and \(\mathcal S\) are non-amenable and \(\mathcal S\) is treeable, then \(\mathcal R\) has cost \(>1\) and \(\rho\) is extension-like;
	\item If both \(\mathcal R\) and \(\mathcal S\) are non-amenable and \(\mathcal S\) is solidly ergodic, then \(\mathcal R\) is solidly ergodic as well and \(\rho\) is extension-like.
\end{itemize}

As mentioned in Remark \ref{rmk: haagerup vs T}, a further natural result would be to rule out non-trivial homomorphisms when \(\mathcal R\) has (T) and \(\mathcal S\) has the Haagerup property, but seems to require to develop these properties in the non-ergodic case, which we leave as a question.

\begin{question}
	Let \(\mathcal R\) be a \pmp equivalence relation with property (T), is it true that if \(\mathcal S\) is a class-bijective extension of \(\mathcal R\) with the Haagerup property, then both \(\mathcal R\) and \(\mathcal S\) are periodic (have only finite classes)?
\end{question}

It would be interesting to investigate non-singular boolean actions of full groups, both of \pmp equivalence relations, but also of non-singular equivalence relations.
When \(\mathcal R=\mathcal R_\Gamma\) comes from a free \pmp \(\Gamma\)-action, one can use a non-singular \(\Gamma\)-action to build a non-singular extensions of \(\mathcal R\) which is not p.m.p., so full groups of \pmp equivalence relations seem to often admit non \pmp boolean actions.
If one wants to follow a similar approach to ours, here is a first important question.

\begin{question}
	What are the quasi-invariant random subgroups of \(\mathfrak S_{2^\infty}\), i.e. the probability measures on \(\Sub(\mathfrak S_{2^\infty})\) which are quasi-preserved by the conjugacy action?
\end{question}




\section{Homomorphisms and actions of \pmp equivalence relations}\label{sec:boolean actions}

The remainder of the paper is devoted to the proof of Theorem~\ref{thmi:charaTfullgroups}, 
which states that for ergodic \pmp equivalence relations,
property (T) is equivalent to  
every boolean \pmp ergodic non-free action of their full group
being strongly ergodic. 

In order to define property (T), we first need to recall the theory 
of unitary representations of \pmp equivalence relations.
We will also need to connect those unitary representations to \pmp actions 
of \pmp equivalence relations via 
the Gaussian construction, and these \pmp actions are not 
directly of the (spatial) form that we have seen so far. 
In the first subsection, we thus develop the theory of booolean actions of \pmp equivalence relations
and show that it is essentially the same thing as a measure-preserving action
(i.e., a class-bijective extension), as defined 
in Section \ref{sec: boolean action from eq rel}. We also observe 
that boolean actions 
can be characterized abstractly in terms of the associated 
full group actions (Remark~\ref{rem: full actions}). 

Section~\ref{sec:urep} is then devoted to unitary representations of 
\pmp equivalence relations. We take the time to obtain 
a similar abstract characterization of these unitary representations 
in therms of the associated full group unitary representation 
(Theorem~\ref{thm: full unitary rep}). This characterization 
will be useful later on,
and we think it should be taken as the definition of a unitary representation
of a \pmp equivalence relation as it is much more algebraic and 
allows one to forget about measurability issues. 
But for now, let us spell out the classical framework for 
boolean actions and unitary representations of ergodic \pmp equivalence relations.

Any equivalence relation \(\mathcal R\) on a set \(X\) can be seen as a groupoid where the product is given by \((x_1,x_2)(x_2,x_3)=(x_1,x_3)\) whenever \((x_1,x_2)\) and \((x_2,x_3)\) belong to \(\mathcal R\).
This allow one to make sense of actions of \(\mathcal R\) on various structures. 


\begin{definition}
	Let \(G\) be a Polish group and let \(\mathcal R\) be a p.m.p.~equivalence relation. A \textbf{homomorphism}, or \textbf{cocycle} is a Borel map \(c: \mathcal R\to G\) such that for all \((x_1,x_2),(x_2,x_3)\in\mathcal R\),
	\[c(x_1,x_3)=c(x_1,x_2)c(x_2,x_3).\]
\end{definition}

We will be interested in the cases where \(G\) is either the automorphism group of a probability space \(\Aut(Z,\eta)\), the unitary group of a Hilbert space \(\mathcal U(\mathcal H)\) or the group of affine isometries of a real Hilbert space \(\Iso(\mathcal H)\).
As we will see, there is a deep connection between these three kinds
of actions stemming from natural embeddings between these three Polish groups. For now, let us stick to \(G=\Aut(Z,\eta)\).

\begin{remark}
The reader acquainted with groupoids will note that our viewpoint amounts to working only with actions on constant bundles, 
thus making our life easier. 
Note that these are not the right kinds of actions when working with non-ergodic groupoids, but we do not consider them here.
\end{remark}

\subsection{Boolean actions of \pmp equivalence relations}

\begin{df}
	A (p.m.p.) \textbf{boolean action} of a \pmp equivalence relation \(\mathcal R\) on \((X,\mu)\) is a homomorphism \(c:\mathcal R\to \Aut(Z,\eta)\), where \((Z,\eta)\) is a standard probability space, possibly with atoms.
\end{df}

Our aim is now to show that for ergodic \pmp equivalence relations, 
boolean actions are the same thing as \pmp actions as defined in Section \ref{sec: boolean action from eq rel}, 
which is the content of Proposition~\ref{prop:boolean yields pmp action} and Proposition~\ref{prop: pmp spatial R action from boolean}.

Let \(\mathcal R\) be a \pmp equivalence relation on the standard probability \((X,\mu)\). Consider another standard probability space \((Z,\eta)\) and let \(c\colon\mathcal R\rightarrow \Aut(Z,\eta)\) be a boolean action. Let \(\Gamma\leq [\mathcal R]\) such that \(\mathcal R\) is the equivalence relation generated by \(\Gamma\). Denote by \(\alpha\) the associated action of \(\Gamma\) on \((X,\mu)\). Then we would like to obtain an action \(\beta\) of \(\Gamma\) on \(X\times Z\) via \[\beta(\gamma)(x,z)=(\gamma x,c(\gamma x,x)z).\]
However, \(c(\gamma x,x)\) is only defined up to null-sets, therefore we have to make sure that such a \(\beta\) would be measurable and an action.

The upcoming Proposition \ref{proposition:l0measurable} will yield that, but first we need a lemma. Recall that given a standard probability space \((X,\mu)\) and a Polish space \(M\), \(\LL^0(X,\mu,M)\) denotes the space of measurable maps \(X\to M\) up to measure zero.
The latter is a Polish space for the topology of convergence in measure, see e.g.\ \cite[Sec.~19]{kechrisGlobalaspectsergodic2010}.
\begin{lemma}\label{lem: l0 malg}
	Let \((X,\mu)\) and \((Z,\eta)\) be standard probability spaces, possibly with atoms.
	Then we have a natural identification between \(\MAlg(X\times  Z,\mu\otimes \eta)\) and \(\LL^0(X,\mu,\MAlg(Z,\eta))\).
\end{lemma}
\begin{proof}
	In this proof, it is convenient for clarity to denote by \([A]\) the equivalence class of a Borel set in the corresponding measure algebra. For \(A\subseteq X\times Z\) and \(x\in X\), let \(A_x\coloneqq \{z\in Z\colon (x,z)\in A\}\).
	Given \(A\subseteq X\times Z\) Borel, the map \(x\mapsto \eta(A_x)\) is Borel, and so for \(B\subseteq Z\) Borel the map \(x\mapsto \mu(A_x\bigtriangleup B)\) is Borel.
	This means that every \(A\subseteq X\times Z\) Borel defines a Borel map \(F_A:x\mapsto [A_x]\) from \(X\) to \(\MAlg(Z,\eta)\).
	Denote by \([F_A]\) the equivalence class of \(F_A\) in \(\LL^0(X,\mu,\MAlg(Z,\eta))\).

	Observe that if \([A]=[A']\) then by Fubini's theorem \([F_A]=[F_{A'}]\)
	In particular, the map \(A\mapsto [F_A]\) quotients down to an injective map \(\Phi:\MAlg(X\times  Z,\mu\otimes \eta)\to \LL^0(X,\mu,\MAlg(Z,\eta))\), and we now only have to show it is surjective.

	Endow \(\LL^0(X,\mu,\MAlg(Z,\eta)\) with  the complete metric \[d(F,G)=\int_X \mu(F(x)\bigtriangleup G(x))d\mu(x).\]
	Applying Fubini's theorem once more, we see that \(\Phi\) is isometric.
	To see that \(\Phi\) is surjective, first observe that the space of maps taking countably many values is dense in \(\LL^0(X,\mu,\MAlg(Z,\eta))\).
	Given such a map \(F\) taking countably many values \([B_0],[B_1],\dots\), let \(A_n=F\inv[B_n]\) (which is well-defined up to a null set).
	Then the set \(A_f\coloneqq \bigcup_n A_n\times B_n\) satisfies \(\Phi(A_F)=F\).
	So \(\Phi\) has dense image, and since both \(d_{\mu\otimes \eta}\) and \(d\) are complete we conclude that \(\Phi\) is a surjective isometry as wanted.
\end{proof}

\begin{proposition}\label{proposition:l0measurable}
	Let \((X,\mu)\) and \((Y,\nu)\) be standard Borel spaces. Then for every Borel function \(F:X\rightarrow \Aut(Z,\eta)\), there is an element \(\tilde F\in\Aut(X\times Z,\mu\times \eta)\) such that for almost every \(x\), we have that \(\tilde F(x,z)=(x,F(x)z)\) for almost all \(z\in Z\).
\end{proposition}
\begin{proof}
	Denote by \(\cdot\) the natural action of \(\Aut(Z,\eta)\) on \(\MAlg(Z,\eta)\).
	Then for every \(f\in\LL^0(X,\mu,\MAlg(Z,\eta))\), we define \(F(f)\in\LL^0(X,\mu,\MAlg(Z,\eta))\) by \[F(f)(x)=F(x)\cdot f(x).\]
	Through the identification between \(\MAlg(X\times Z,\mu\otimes\eta\) and \(\LL^0(X,\mu,\MAlg(Z,\eta))\) provided by the previous lemma, \(F\) thus defines an automorphism \(T_F\) of \(\MAlg(X\times Z,\mu\otimes \eta)\).

	Let us denote by \(\tilde F\) an arbitrary lift of \(T_F\) to a measure-preserving bijection of \((X\times Z,\mu\otimes\eta)\), 
	which we write as \(T_F(x,z)=(\alpha(x,z),\beta(x,z))\).
	By construction \(T_F\) preserves elements of \(\MAlg(X\times Z,\mu\otimes \eta)\) of the form \(A\times Z\) where \(A\in\MAlg(X,\mu)\), so \(T_F\) almost preserves Borel sets of the form \(A\times Z\).
	If \((C_n)\) is a separating sequence of Borel subsets of \(X\), we conclude that after throwing away a null set, \(\tilde F\) preserves every set of the form \(C_n\times Z\), which implies that
	\(\alpha(x,z)=x\) for almost all \(x\in X\) and almost all \((x,z)\in X\times Z\).

	We finally conclude that \(\beta(x,\cdot)=F(x)\) for almost every \(x\)
	from the fact that up to measure zero, \(\beta(x,\cdot)\) has to be by construction a measure-preserving bijection of \(Z\) which is equal to \(F(x)\) in \(\Aut(Z,\eta)\).
\end{proof}

We can finally check that boolean actions yield \pmp actions.

\begin{proposition}\label{prop:boolean yields pmp action}
	Let \(\mathcal R\) be a p.m.p.\ equivalence relation on the standard probability \((X,\mu)\) and let \(c\colon \mathcal R\rightarrow \Aut(Z,\eta)\) be a Borel cocycle. Then there is an essentially unique \pmp equivalence relation \(\mathcal S\) on \((X\times Z,\mu\otimes\eta)\) such that for almost every \(((x,z),(x',z'))\in\mathcal S\) we have that \(z'=c(x,x')z\). Moreover the projection \(\pi:X\times Y\rightarrow X\) is a class-bijective factor map from \(\mathcal S\) to \(\mathcal R\).
\end{proposition}
\begin{proof}
	Let \(\alpha:\Gamma\curvearrowright (X,\mu)\) generate the equivalence relation \(\mathcal R\).
	We define a new \(\Gamma\)-action \(\beta\) on \((X\times Z,\mu\otimes\eta)\) by \(\beta(\gamma)(x,z)=(\alpha(\gamma) x,c(\alpha(\gamma) x,x)z)\). Proposition \ref{proposition:l0measurable} guarantees us that after possibly modifying it on a null set, \(\beta\) is a well-defined \pmp action such that the projection \(\pi:X\times Z\to X\) is  \(\Gamma\)-equivariant.

	Denote by \(\mathcal S\) the equivalence relation generated by \(\beta\).
	Since \(\pi\) is \(\Gamma\)-equivariant, the map \(\pi\) is factor map from \(\mathcal S\).
	To see that it is class-bijective, we need to show that it is injective when restricted to a \(\beta\)-orbit. By Remark \ref{rmk: support pres is class bij}, it suffices to show that if \(\alpha(\gamma)x=x\), then for all \(z\in Z\) we have \(\beta(\gamma)(x,z)=(x,z)\), which is clear since \(c(\alpha(\gamma)x,x)=c(x,x)=\id_Z\). Therefore the projection is a class-bijective factor map \(\mathcal S\to\mathcal R\).
\end{proof}

Using Rokhlin's skew product theorem, we now show that conversely, every \pmp action of an ergodic \pmp equivalence relation is a boolean action.

\begin{prop}\label{prop: pmp spatial R action from boolean}
	Let \(\mathcal R\) be a \pmp ergodic equivalence relation, let \(\mathcal S\) be a \pmp equivalence relation on \((Y,\nu)\) which is a class-bijective extension of \(\mathcal R\) via a map \(\pi: Y\to X\).
	Then one can decompose \((Y,\nu)\) as a product space \(Y=X\times Z\) endowed with the product measure \(\mu\otimes \eta\), where \(\eta\) is a probability measure on the standard Borel space \(Z\), possibly with atoms.

	Moreover in this decomposition \(\pi\) becomes the projection on the first coordinate and we get a cocycle \(c:\mathcal R\to \Aut(Z,\eta)\) such that for all \(z\in Z\)
	\(c(x,x')z\) is the unique \(z'\in Z\) such that \(((x,z),(x',z'))\in\mathcal S\) and \(\pi((x',z'))=x'\).
\end{prop}
\begin{proof}
	Let \(\Gamma\leq[\mathcal R]\) be a countable group such that \(\mathcal R=\mathcal R_\Gamma\).
	Denote by \(\beta\) the \(\Gamma\)-action on \((Y,\nu)\) provided by the injection \(\rho_\pi:[\mathcal R]\to[\mathcal S]\) as in Lemma \ref{lemma: def rho pi}, which is a class-bijective extension of the \(\Gamma\)-action on \((X,\mu)\) provided by the inclusion \(\Gamma\leq [\mathcal R]\).

	Since \(\mathcal R=\mathcal R_\Gamma\) and \(\mathcal R\) is ergodic, the \(\Gamma\)-action on \((X,\mu)\) provided by the inclusion \(\Gamma\leq [\mathcal R]\) is ergodic.
	By Rokhlin's skew product theorem (see \cite[Thm.~3.18]{glasnerErgodicTheoryJoinings2003}), we get the desired decomposition of \((Y,\nu)\) as a product \((X\times Z,\mu\otimes\eta)\) and a cocycle \(\tilde c: \Gamma\times X\to \Aut(Z,\nu)\) such that for all \((x,z)\in X\times Z\),
	\(\beta(\gamma)(x,z)=(\gamma x, \tilde c(\gamma,x)z)\).
	Since \(\pi\) is class-bijective, we get a well-defined cocycle \(c:\mathcal R\to\Aut(Z,\eta)\) by letting \(c(x,\gamma x)=\tilde c(x,\gamma)\).
\end{proof}

We finally make the following convenient definition, which will be followed in the next section by its analogue for unitary representations.

\begin{definition}\label{def: action full group}
	Given a boolean action \(\alpha:\mathcal R\to \Aut(Z,\eta)\) of
	a \pmp equivalence relation \(\mathcal R\) on \((X,\mu)\),
	we denote by
	\([\alpha]\) the boolean measure-preserving action of \([\mathcal R]\) on
	\((X\times Z, \mu\otimes\eta)\) corresponding to the class-bijective extension of \(\mathcal R\) given by the previous proposition. It is defined by: for all \(T\in [\mathcal R]\),
	\[
		[\alpha](T)(x,z)=(T(x),\alpha(x,T(x))z).
	\]
\end{definition}

\begin{remark}\label{rem: full actions}
	It follows from Lemma \ref{lem:pmp action from support pres} 
	and the correspondence between boolean and \pmp actions of equivalence relations that the boolean actions \([\alpha]\) of a full group \([\mathcal R]\) 
	are also characterized as those which factor onto the inclusion \(\iota:[\mathcal R]\into\Aut(X,\mu)\) in a support-preserving manner.
	The next section is mainly devoted to obtaining a similar characterization 
	for unitary representations.
\end{remark}

\subsection{Unitary representations of \pmp equivalence relations}\label{sec:urep}

\begin{definition}\label{df: various actions of eqrel}
	Let \(\mathcal R\) be an ergodic  p.m.p.~equivalence relation on a standard probability space \((X,\mu)\).
	A \textbf{unitary representation} of \(\mathcal R\)	is a homomorphism \(\pi:\mathcal R\to \mathcal U(\mathcal H)\), where \(\mathcal U(\mathcal H)\) is the unitary group of a separable Hilbert space \(\mathcal H\).
\end{definition}
Given unitary representation \(\pi: \mathcal R\to\mathcal U(\mathcal H)\) of a p.m.p.~equivalence relation \(\mathcal R\)
on \((X,\mu)\),
a \textbf{section} of the representation is a Borel map \(\xi: X\to \mathcal H\).
It is called \textbf{square-integrable} when \(\int_X \norm{\xi(x)}^2d\mu(x)\) is finite.
The space of square-integrable sections is naturally a Hilbert space for the scalar product
\(\la\xi,\eta\ra=\int_X \la \xi(x),\eta(x)\ra d\mu(x)\), and we denote this Hilbert space by \(\LL^2(X,\mu,\mathcal H)\).
Observe that \(\LL^2(X,\mu,\mathcal H)\) has a natural \(\LL^\infty(X,\mu)\)-module structure with \(\LL^\infty(X,\mu)\)
acting by multiplication on sections.
Moreover, one can associate to \(\pi\)
a unitary representation \([\pi]\) of the full group of \(\mathcal R\)  on \(\LL^2(X,\mu,\mathcal H)\): for all \(T\in[\mathcal R]\)
and all  \(\xi\in\LL^2(X,\mu,\mathcal H)\), we let
\[
	([\pi](T)\xi) (T(x))=\pi(x,T(x))\xi(x).
\]
We will now characterize algebraically unitary representations of the full group of the form \([\pi]\) in a manner analogous to Remark \ref{rem: full actions}. 
We need the following additional structure.

\begin{definition}
	An \textbf{\(\LL^\infty(X,\mu)\)-module} is a Hilbert space \(\mathcal K\) an a \emph{faithful} normal
	\(*\)-representation of \(\LL^\infty(X,\mu)\) in \(\mathcal B(\mathcal K)\).
\end{definition}

Given an \textbf{\(\LL^\infty(X,\mu)\)-module} \(\mathcal K\), we identify \(\LL^\infty(X,\mu)\) to its image in \(\mathcal B(\mathcal K)\).
For all \(A\subseteq X\), we denote by \(\mathcal K_A\) the subspace \(\chi_A\mathcal K\).

Our motivating example of \(\LL^\infty(X,\mu)\)-modules is provided by unitary representations of p.m.p.~equivalence relations: indeed the space of square-integrable sections is easily checked to be an
\(\LL^\infty(X,\mu)\)-module.

Here are the abstract properties of \([\pi]\) which will allow us to reconstruct \(\pi\).

\begin{definition}\label{df: full urep}
	Let \(\mathcal R\) be a p.m.p.~equivalence relation on \((X,\mu)\),
	let \(\rho: [\mathcal R]\to \mathcal U(\mathcal K)\) be a unitary representation on a Hilbert space \(\mathcal K\)
	equipped with an \(\LL^\infty(X,\mu)\)-module structure. Say that \(\rho\) is \textbf{full} when it satisfies the following two conditions:
	\begin{enumerate}[(i)]
		\item \label{cond: full urep equiv}For all \(A\subseteq X\) Borel, \(\rho(T)\chi_A=\chi_{T(A)}\rho(T)\);
		\item \label{cond: full urep is full}For all \(A\subseteq X\), and all \(T,U\in[\mathcal R]\) such that \(T_{\restriction A}=U_{\restriction A}\),
		      we have \(\rho(T)_{\restriction \mathcal K_A}=\rho(U)_{\restriction \mathcal K_A}\).
	\end{enumerate}
\end{definition}
Note that the first condition yields that \(\rho(T)\mathcal K_A=\mathcal K_{T(A)}\) for all \(T\in[\mathcal R]\),
while the second is equivalent to asking that \(\rho(T)_{\restriction\mathcal K_{\Fix(T)}}=\id_{\mathcal K_{\Fix(T)}}\)
for all \(T\in [\mathcal R]\).

\begin{theorem}
	\label{thm: full unitary rep}
	Let \(\mathcal R\) be an ergodic p.m.p.~equivalence relation on \((X,\mu)\).
	Every full unitary representation \(\rho:[\mathcal R]\to\mathcal U(\mathcal K)\) is unitarily equivalent to
	a unitary representation of the form \([\pi]\) for some unitary representation \(\pi\) of \(\mathcal R\).
\end{theorem}
\begin{proof}
	By \cite[Thm.~14.2.1]{kadisonFundamentalsTheoryOperator1997}, the \(\LL^\infty(X,\mu)\)-module
	\(\mathcal K\) can be decomposed
	as a direct integral of Hilbert spaces \((\mathcal H_x)_{x\in X}\). By condition
	\eqref{cond: full urep equiv} and ergodicity the dimension of \(\mathcal H_x\) is almost
	surely constant. We can thus assume \(\mathcal H_x\) is constant equal to  a fixed
	Hilbert space \(\mathcal H\) and so
	\(\mathcal K=\LL^2(X,\mu,\mathcal H)\).

	Denote by \(\lambda:[\mathcal R]\to \LL^2(X,\mu,\mathcal H)\) be the unitary representation defined by
	\((\lambda(T)\xi)(T(x))=\xi(x)\).
	Then for all \(A\subseteq X\) Borel, we have \(\lambda(T)\chi_A=\chi_{T(A)}\lambda(T)\).
	So by condition \eqref{cond: full urep equiv}, for every \(T\in [\mathcal R]\),
	the unitary \(\lambda(T)^*\rho(T)\) commutes with \(\LL^\infty(X,\mu)\) and is thus a decomposable operator
	\[
		\lambda(T)^*\rho(T)=\int_X c(T,x)d\mu(x),
	\]
	where \(c(T,x)\in \mathcal U(\mathcal H)\) for almost all \(x\in X\).
	Observe that by uniqueness of the decomposition of a decomposable operator, for any \(U\in [T]\) we have
	\(\lambda(U)\int_X c(T,x)d\mu(x)\lambda(U)^*=\int_Xc(T,U\inv(x))d\mu(x)\).
	When \(T_1,T_2\in [\mathcal R]\), we thus have
	\begin{align*}
		\lambda(T_1T_2)^*\rho(T_1T_2) & =\lambda(T_2)^*\lambda(T_1)^*\rho(T_1)\rho(T_2)      \\
		                              & =\lambda(T_2)^*\int_ X c(T_1,x)d\mu(x)\rho(T_2)      \\
		                              & =\int_X c(T_1,T_2(x))d\mu(x) \lambda(T_2)^*\rho(T_2) \\
		                              & =\int_Xc(T_1,T_2(x))c(T_2,x)d\mu(x),
	\end{align*}
	so for almost all \(x\in X\), \(c(T_1,T_2(x))c(T_2,x)=c(T_1T_2,x)\).
	We now pick a countable group \(\Gamma\) such that \(\mathcal R=\mathcal R_\Gamma\). By
	restricting to a full measure set, we may assume that the above cocycle relation holds for
	all \(\gamma\in\Gamma\) and all \(x\in X\).

	Furthermore by condition \eqref{cond: full urep is full}, we have that \(\rho(\gamma)\) is trivial when restricted to
	\(\mathcal K_{\Fix(\gamma)}\), so up to restricting again we have \(c(\gamma,x)=\id_{\mathcal H}\)
	for all \(\gamma\in\Gamma\) and all \(x\in \Fix(\gamma)\).
	It follows that we have a well defined unitary representation \(\pi:\mathcal R\to \mathcal U(\mathcal H)\)
	by letting \(\pi(x,\gamma\cdot x)=c(\gamma,x)\). By construction we have
	\(\rho(T)=\lambda(T)\int_X c(T,x)d\mu(x)\), so for all \(\xi\in \LL^2(X,\mu,\mathcal H)\) and all \(x\in X\)
	we have \[\rho(T)(\xi)(T(x))=\lambda(T)\int_Xc(T,y)d\mu(y)\xi(T(x))=c(T,x)\xi(x).\]
	Now if we pick \(\gamma\in\Gamma\) such that \(\gamma\cdot x=T(x)\) we have by condition \eqref{cond: full urep is full}
	\[c(T,x)\xi(x)=c(\gamma,x)\xi(x)=\pi(x,\gamma\cdot x)\xi(x)=\pi(x,T(x))\xi(x).\]
	We thus have \(\rho=[\pi]\) as wanted.
\end{proof}

\section{Property (T)}\label{sec: T}

The goal of this Section is threefold: we want to define 
property (T), connect it to strong ergodicity, and use this 
to show the direct implication in Theorem~\ref{thmi:charaTfullgroups}, 
namely that if an ergodic \pmp equivalence relation has property (T), 
then all the \pmp ergodic non-free boolean actions of its full group 
are strongly ergodic. Given Proposition~\ref{thm: main pmp equivalence relations}, 
we need to show that if 
$\mathcal R$ is ergodic and has (T), then for all $n\geq 1$,
every \pmp ergodic action of $\mathcal R^{\odot n}$ is strongly 
ergodic. 

This is done by noting the well-known fact that property (T) implies strong ergodicity
(Proposition~\ref{prop: T implies st ergo}),
and that property (T) passes to class-bijective extensions 
(Prop.~\ref{prop: (T) action induces (T) eq rel}).

Further, we need to know that if $\mathcal R$ has (T) then so does $\mathcal R^n$, 
and this is actually the result which takes up most of this section, 
as this requires us to have a more quantitative version of property (T)
(a Kazhdan pair with a quantitative estimate on the distance between
almost invariant vectors and invariant vectors). This is done in 
the two following sections, after recalling the definition of property (T)
for \pmp equivalence relations. The result that if $\mathcal R$
has (T) then so does $\mathcal R^{\odot n}$ is proven in Section~\ref{sec: perm properties of T}.
All this allows us to prove the direct implication from Theorem~\ref{thmi:charaTfullgroups}
in Section~\ref{sec:directimpl}.

\subsection{First definition and strong ergodicity}

Let \(\mathcal R\) be a measure-preserving equivalence relation on \((X,\mu)\).
Given a unitary representation \(\pi: \mathcal R\to \mathcal U(\mathcal H)\), a  \textbf{unit section}
is a section \(\xi:X\to\mathcal H\) satisfying
that for all \(x\in X\), \(\norm{\xi(x)}=1\).
%
%
We say that a sequence \((\xi_n)\) of unit sections is \textbf{almost invariant}
when the sequence of functions \((x,y)\in\mathcal R\mapsto \norm{\pi(x,y)\xi_n(x)-\xi_n(y)}\)
converges pointwise to zero, up to restricting to a full measure subset of \(X\).

\begin{remark}
	Equivalently, one could take a weaker definition for almost invariant unit sections
	by asking that the sequence of functions
	\((x,y)\mapsto \norm{\pi(x,y)\xi_n(x)-\xi_n(y))}\) converges to zero \emph{in measure}.
	Since every sequence converging in measure to zero has a subsequence which converges
	pointwise to zero, this does not affect the definition of property (T)
	(see \cite[Lem.~4.1]{anantharaman-delarocheCohomologyPropertyGroupoids2005}).
\end{remark}

\begin{definition}
	A p.m.p.~equivalence relation \(\mathcal R\) has \textbf{property (T)} if whenever a unitary representation \(\pi\) of \(\mathcal R\)
	has almost invariant unit sections, then it has an invariant unit section, i.e. we can find a section \(\xi\)
	such that for almost all \((x,y)\in\mathcal R\), we have \(\pi(x,y)\xi(x)=\xi(y)\).
\end{definition}

\begin{remark}
	An ergodic p.m.p.~equivalence relation \(\mathcal R\) has property (T) iff every unitary representation
	with almost invariant unit sections has a non-zero invariant section. Indeed such a section must have an invariant norm which is thus constant by ergodicity. This constant must be non zero, so by rescaling we
	obtain an invariant unit section.
\end{remark}

We now connect property (T) to the fundamental notion of strong ergodicity.
Given a p.m.p.~equivalence relation \(\mathcal R\) on a standard probability space \((X,\mu)\), a sequence
\((A_n)\) of measurable subsets of \(X\) is called \textbf{asymptotically invariant} if for all
\(T\in [\mathcal R]\), \(\mu(T(A_n)\bigtriangleup A_n)\to 0\), and \textbf{asymptotically trivial} if
\(\mu(A_n)(1-\mu(A_n))\to 0\).

\begin{definition}\label{df: strong ergo eq rel}
	A p.m.p.~equivalence relation \(\mathcal R\) on a standard probability space 
	\((X,\mu)\) is \textbf{strongly ergodic} if every asymptotically invariant sequence of subsets of \(X\) has to be asymptotically trivial.
\end{definition}

By considering, for an \(\mathcal R\)-invariant set \(A\), the constant sequence equal to \(A\), we see that
strong ergodicity implies ergodicity.

We now note the well-known fact (implicit in \cite{pichotTheorieSpectraleRelations2007}) that for an ergodic p.m.p.~equivalence relation,
property (T) implies strong ergodicity.

\begin{proposition}[Folklore]\label{prop: T implies st ergo}
	Let \(\mathcal R\) be an ergodic p.m.p.~equivalence relation with property (T). Then \(\mathcal R\) is strongly ergodic.
\end{proposition}
\begin{proof}
	Suppose by contradiction that \(\mathcal R\) has (T) but is not strongly ergodic.
	Since \(\mathcal R\) is ergodic but not strongly ergodic, by \cite[Thm.~2.2 and Rem.~2.5]{jonesAsymptoticallyInvariantSequences1987} there is a measure-preserving
	map \(\pi:(X,\mu)\to (\{0,1\}^\N,\mathcal B(1/2)^{\odot \N})\) such that \((\pi\times\pi)(\mathcal R)=\mathcal R_0\), where \(\mathcal R_0\) is the ergodic hyperfinite 
	\pmp equivalence relation as in Definition \ref{def:R0}.

	For every \(k\in\N\), let \(\mathcal S_k\) be the equivalence relation on \(\{0,1\}^\N\) defined
	by \((x_n)\mathcal S_k(y_n)\) if and only if \(x_n=y_n\) for all \(n\geq k\). We already observed that \(\mathcal R_0=\bigcup_{k\in\N}\mathcal S_k\), so if \(\mathcal R'_k\coloneqq (\pi\times \pi) \inv(\mathcal S_k)\) then
	\(\mathcal R=\bigcup_{k \in \N}\mathcal R'_k\).

	But since \(\mathcal R\) has property (T) it is not \emph{approximable} :
	there is \(k\in\N\) and a non-negligible Borel set
	\(A\subseteq X\) such that \(\mathcal R\cap (A\times A)=\mathcal R'_k\cap (A\times A)\)
	(see \cite[Prop.~16]{pichotTheorieSpectraleRelations2007}; see also \cite{gaboriauApproximationsStandardEquivalence2016} for more on approximability).
	Now observe that each \(\mathcal S_k\) has
	diffuse ergodic decomposition because it has only finite classes, hence so does each \(\mathcal R'_k\).
	However \(\mathcal R\) is ergodic, so we
	cannot have  \(\mathcal R\cap (A\times A)=\mathcal R'_k\cap (A\times A)\), which is the desired contradiction.
\end{proof}

Our aim is now to give a notion of Kazhdan pair (and triple) which provides us more quantitative versions
of property (T) as in \cite[Prop.~1.1.9]{bekkaKazhdanProperty2008}.
Such a characterization of property (T) will be given purely
in terms of full unitary representations of the full group, using implicitely Theorem \ref{thm: full unitary rep}.
This will rely crucially on the following spectral gap characterization of strong ergodicity due to Houdayer, Marrakchi and Verraedt \cite{houdayerStronglyErgodicEquivalence2019}.

\begin{thm}[Houdayer, Marrakchi, Verraedt]\label{thm: spec gap}
	Let \(\mathcal R\) be a strongly ergodic p.m.p.~equivalence relation. Then there is a finite
	set \(F\subseteq [\mathcal R]\) and \(\kappa>0\) such that for all \(\eta\in \LL^2(X,\mu)\),
	\[\norm{\eta-\int_X\eta(x)d\mu(x)}^2\leq \kappa\sum_{T\in F}\norm{\eta-\eta\circ T\inv}^2.\]
\end{thm}

When we have a finite set \(F\subseteq[\mathcal R]\) and \(\kappa>0\) as above, we say that \((F,\kappa)\) is
a \textbf{spectral gap pair} for \(\mathcal R\).
\subsection{Kazhdan pairs and proximity of invariant vectors} \label{sec:kazhdan pair}

As explained before,  one can associate to every unitary representation \(\pi\) of a p.m.p.~equivalence relation \(\mathcal R\)
a full unitary representation \([\pi]\) of its full group  on the space of square integrable sections: for all \(T\in[\mathcal R]\)
and all square-integrable section \(\xi\), we let
\[
	([\pi](T)\xi) (T(x))=\pi(x,T(x))\xi(x).
\]
Observe that a square-integrable section is \([\pi]\)-invariant if and only if it is a \(\pi\)-invariant section.
Given a unitary representation \(\pi:G\to\mathcal U(\mathcal H)\) of a group \(G\) and a finite subset \(F\) of \(G\), we say that
a vector \(\xi\in\mathcal H\) is \((F,\epsilon)\)\textbf{-invariant} if for all \(g\in F\),
\[\norm{\pi(g)\xi-\xi}<\epsilon \norm\xi.\]

\begin{definition}
	Let \(\mathcal R\) be a p.m.p.~equivalence relation, let \(F\) be a finite subset of \([\mathcal R]\) and let \(\epsilon>0\).
	Then \((F,\epsilon)\) is a \textbf{Kazhdan pair} for the equivalence relation \(\mathcal R\) if	whenever \(\pi\) 
	is a unitary representation of \(\mathcal R\) such that $[\pi]$ 
	admits an \((F,\epsilon)\)-invariant vector,
	then \([\pi]\) has a non-zero invariant vector.
\end{definition}

\begin{remark}
	Theorem~\ref{thm: full unitary rep} allow us 
	to restate Kazhdan pairs as follows: \((F,\epsilon)\) is a Kazhdan pair
	for $\mathcal R$ if and only if every \emph{full}
	unitary representation of the full group $[\mathcal R]$
	with $(F,\epsilon)$-invariant vectors must have 
	a non-zero invariant vector. 
\end{remark}

\begin{proposition}
	A p.m.p.~ergodic equivalence relation \(\mathcal R\) has property (T) if an only if it admits
	a Kazhdan pair.
\end{proposition}
\begin{proof}
	Suppose \((F,\epsilon)\) is a Kazhdan pair for \(\mathcal R\). By the dominated convergence theorem,
	if \((\xi_n)\) is a sequence of almost invariant unit sections for a unitary representation \(\pi\) of
	\(\mathcal R\), then for all \(T\in F\) we have \(\norm{\pi(T)\xi_n-\xi_n}\to 0\). In particular
	for \(n\) large enough, \(\xi_n\) is \((F,\epsilon)\)-invariant and thus we have a non-zero invariant
	section, which by ergodicity and rescaling yields an invariant unit section.

	Conversely, assume that \(\mathcal R\) has (T) but no Kazhdan pair.
	Then \(\mathcal R\) is strongly ergodic, so by Theorem \ref{thm: spec gap} it has a spectral gap pair
	\((F,\kappa)\).
	By the Feldman-Moore theorem,
	there is an increasing sequence \((F_n)\) of finite subsets of \([\mathcal R]\) so that
	\(\mathcal R\) is covered by the reunion over \(n\in\N\) of the graphs of the elements of \(F_n\).
	For every \(n\in\N\), \((F\cup F_n,\frac 1n)\) is not a Kazhdan pair so we find a unitary representation
	\(\pi_n\) of \(\mathcal R\) on \(\mathcal H_n\) so that \([\pi_n]\) has an \((F\cup F_n,\frac 1n)\)-invariant vector \(\xi_n\) but no non-zero invariant vector.

	For each \(n\in\N\), define \(\eta_n(x)=\norm{\xi_n(x)}\).
	By the reversed triangle inequality, we have for all \(T\in [\mathcal R]\):
	\begin{align*}
		\norm{[\pi](T)\xi_n-\xi_n}^2 & =\int_X\norm{\pi(T)\xi_n(x)-\xi_n(x)}^2d\mu(x)                         \\
		                             & \geq \int_X\left(\norm{\pi(T)\xi_n(x)}-\norm{\xi_n(x)}\right)^2d\mu(x) \\
		                             & = \int_X\abs{\eta_n(T\inv(x))-\eta_n(x)}^2d\mu(x).
	\end{align*}

	Since \((F,\kappa)\) is a spectral gap pair, we then
	have
	\[
		\norm{\eta_n-\int_X\eta_n(x)d\mu(x)}^2\leq \kappa \sum_{T\in F}\norm{\eta_n-\eta_n\circ T\inv}^2
		\leq \kappa \sum_{T\in F} \norm{\pi(T)\xi_n-\xi_n}^2\leq \frac{\kappa\abs F}{n^2}.
	\]
	Moreover since \(\norm{\eta_n}=1\) we deduce from the reversed triangle inequality 
	\[
		\abs{1-\int_X\eta_n(x)d\mu(x)}^2\leq
		\norm{\eta_n-\int_X\eta_n(x)d\mu(x)}^2\leq
		\frac{\kappa \abs F}{n^2},
	\]
	so by the triangle inequality
	\[\norm{\eta_n-1}\leq \frac{2\sqrt{\kappa\abs F}}n.\]
	Let \(v_n\in \mathcal H_n\) be a fixed unit vector.
	Define \(\xi'_n: X\to \mathcal H_n\) by
	\[
		\xi'_n(x)=
		\left\{
		\begin{array}{cl}
			\frac 1{\eta_n(x)}\xi_n(x) & \text{if }\eta_n(x)\neq 0 \\
			v_n                        & \text{ otherwise.}
		\end{array}
		\right.
	\]
	Each \(\xi'_n\) is by construction a unit section for \(\pi_n\).
	For all $x$ such that $\eta_n(x)\neq 0$, we have 
	\[
	\norm{\xi_n(x)-\xi'_n(x)}^2=\left( 1- \frac 1{\eta_n(x)}\right)^2\eta_n(x)^2=(1-\eta_n(x))^2,
	\]
	and the equality $\norm{\xi_n(x)-\xi'_n(x)}^2=(1-\eta_n(x))^2$ remains true if $\eta_n(x)=0$. We thus have
	\begin{align*}
		\norm{\xi_n-\xi'_n}^2 &=\int_{X}(1-\eta_n(x))^2d\mu(x)\\
		&=\norm{1-\eta_n}^2\\
		&\leq \frac{4\kappa\abs F}{n^2}.
	\end{align*}
	Since \(\xi_n\) is \((F_n,\frac 1n)\)-invariant, it follows from the triangle inequality that each \(\xi'_n\) is an \(\left(F_n, \frac 1n+\frac{4\sqrt{\kappa\abs F}}n\right)\)-invariant
	unit section.

	Consider the infinite direct sum \(\mathcal H\coloneqq\bigoplus_n \mathcal H_n\), we have a unitary representation \(\bigoplus_n \pi_n\) of \(\mathcal R\) on \(\mathcal H\).
	Each \(\xi'_n\) defines a unit section of \(\bigoplus_n \pi_n\)
	and since \(\mathcal R\) is the union of the graphs of the elements of \(F_n\), by \cite[Lem.~4.1]{anantharaman-delarocheCohomologyPropertyGroupoids2005} after taking a subsequence the sequence \((\xi'_n)\) is almost invariant. But since no \(\pi_n\) had an invariant section, neither does their direct sum, contradicting that \(\mathcal R\) had (T) as wanted.
\end{proof}

We can now give an even more quantitative version of property (T) by controlling how far our invariant
vector will be from the \((F,\epsilon)\)-invariant vector, obtaining the desired p.m.p.~equivalence relation
version of \cite[Prop.~1.1.9]{bekkaKazhdanProperty2008} (see also \cite[Thm.~20]{pichotTheorieSpectraleRelations2007}
for a related sequential and pointwise version of what follows).

\begin{proposition}\label{prop: quant (T)}
	Let \(\mathcal R\) be an ergodic p.m.p.~equivalence relation with a Kazhdan pair \((F_1,\epsilon)\) and a spectral gap pair \((F_2,\kappa)\). Suppose that \(\pi:\mathcal R\to \mathcal U(\mathcal H)\) is a unitary representation  and that we have a section \(\xi\in \LL^2(X,\mu,\mathcal H)\) which is \((F_1\cup F_2, \delta\epsilon)\)-invariant for some \(\delta>0\). Then \([\pi]\) admits an invariant vector which is at distance at most \(\delta(1+\epsilon\sqrt{\kappa\abs{F_2}})\norm\xi\) \mbox{from \(\xi\).}
\end{proposition}
\begin{proof}
	Let \(\mathcal K_1\) be the subspace consisting of \([\pi]\)-invariant vectors.
	Denote by \(p_1\) the orthogonal projection onto \(\mathcal K_1\).
	Let \(\mathcal K\) be be the \(\LL^\infty(X)\)-module spanned by \(\mathcal K_1\).
	Observe that the restriction of \([\pi]\) to \(\mathcal K\) is a multiple of \([\theta]\)
	where \(\theta\) is the trivial unitary representation  of \(\mathcal R\) (on \(\C\)).
	Write \(\xi=\xi'+\xi''\) with \(\xi'\in\mathcal K\) and \(\xi''\in \mathcal K^\perp\), then by construction
	\(p_1(\xi')\) is the invariant vector which is the closest to \(\xi\).

	Since \(\mathcal K^\perp\) is a \([\pi]\)-invariant \(\LL^\infty(X)\)-module without non-zero invariant vectors and
	\((F_1,\epsilon)\) is a Kazhdan pair, there is some \(T\in F_1\) such that
	\(\norm{\xi''-[\pi](T)\xi''}\geq \epsilon \norm{\xi''}\).  On the other hand,
	we have \(\norm{\xi''-[\pi](T)\xi''}\leq \norm{\xi-[\pi](T)\xi}<\epsilon\delta\norm{\xi}\), so
	\(\norm{\xi''}<\delta\norm\xi\).

	Finally, since the restriction of \([\pi]\) to \(\mathcal K\) is a multiple of \([\theta]\), we have
	\[\norm{\xi'-p_1(\xi')}^2\leq \kappa\sum_{T'\in F_2}\norm{\xi'-[\pi](T')\xi'}^2\leq \kappa\delta^2\epsilon^2\abs{F_2}\norm{\xi'}^2 \leq \kappa\delta^2\epsilon^2\abs{F_2}\norm{\xi}^2.\]
	Since \(p_1(\xi)=p_1(\xi')\), we finally have
	\begin{align*}\norm{\xi-p_1(\xi)}\leq \norm{\xi''}+\norm{\xi'-p_1(\xi')}\leq \delta(1+\epsilon \sqrt{\kappa\abs{F_2}})\norm\xi.
		 & \qedhere\end{align*}
\end{proof}

\subsection{Some permanence properties}\label{sec: perm properties of T}

We begin this section with a result which is a direct consequence of Anantharaman-Delaroche's Theorem 5.15 from \cite{anantharaman-delarocheCohomologyPropertyGroupoids2005}.
We provide a version of her proof
in terms of full
unitary representations and Kazhdan pairs for the convenience of the reader.

\begin{proposition}\label{prop: (T) action induces (T) eq rel}
	Let \(\mathcal R\) be an ergodic p.m.p.~equivalence relation on \((X,\mu)\) with property (T),
	suppose \(\mathcal S\) is a p.m.p. equivalence relation on \((Y,\nu)\) which is a class-bijective extension of \(\mathcal R\).
	Then \(\mathcal S\) has property (T).
\end{proposition}
\begin{proof}
	Let \(\pi: Y\to X\) denote the class-bijective factor map from \(\mathcal S\) to \(\mathcal R\),
	and let \( \rho_\pi: [\mathcal R]\to[\mathcal S]\) be the corresponding boolean action, which is a support-preserving extension of the inclusion \(\iota:[\mathcal R]\to \Aut(X,\mu)\) (see Lemma~\ref{lemma: def rho pi}).
	Let \((F,\epsilon)\) be a Kazhdan pair for \(\mathcal R\), we will actually show that \((\rho_\pi(F),\epsilon)\)
	is a Kazhdan pair for \(\mathcal S\).

	Let \([\pi]: [\mathcal S]\to \mathcal U(\mathcal H)\) be a full unitary representation of \(\mathcal S\) with
	a \((\tilde\alpha(K),\epsilon)\)-invariant vector, we view its space
	of square-integrable sections as an \(\LL^\infty(X\times Y)\)-module. In particular,
	it is an \(\LL^\infty(X)\)-module and it follows from the fact that \(\tilde\alpha\) is a support-preserving extension of \(\iota\) that
	\([\pi]\circ\tilde\alpha\) is a full unitary representation of \([\mathcal R]\).
	Since \([\pi]\) has a \((\tilde \alpha(F),\epsilon)\)-invariant vector and \((F,\epsilon)\) is a Kazhdan pair
	for \(\mathcal R\), there is a non-zero \(\pi(\tilde \alpha([\mathcal R]))\)-invariant vector \(\xi\in \LL^2(X\times Y,\mathcal H)\).
	By fullness and the fact that the full group generated by \(\tilde \alpha([\mathcal R])\) is equal to \([\mathcal S]\),
	this vector is actually \(\pi([\mathcal S])\)-invariant, which concludes the proof.
\end{proof}

We will now show that if \(\mathcal R\) is a p.m.p.~ergodic equivalence relation with property (T), then
\(\mathcal R^{\odot n}\) has property (T) as well. The following proposition already yields that \(\mathcal R^n\)
has property (T).

\begin{proposition}\label{prop: T carries over to products}
	Let \(\mathcal R\) and \(\mathcal S\) be two p.m.p.~ergodic equivalence relations on \((X,\mu)\) and \((Y,\nu)\) respectively.
	Suppose that both \(\mathcal R\) and \(\mathcal S\) have property (T).
	Then the equivalence relation \(\mathcal R\times\mathcal S\) on \(X\times Y\) is ergodic and has property (T).
\end{proposition}
\begin{proof}
	The ergodicity of \(\mathcal R\times\mathcal S\) is well-known, and follows from the fact that ergodic
	full groups act transitively on sets of the same measure.

	Let \((F_{\mathcal S},\epsilon_{\mathcal S})\) be a Kazhdan pair for \(\mathcal S\). By Proposition \ref{prop: quant (T)}, we find
	a Kazhdan pair \((F_{\mathcal R},\epsilon_{\mathcal R})\) such that given any full unitary representation of \([\mathcal R]\), if there is an \((F_{\mathcal R},\epsilon_{\mathcal R})\)-invariant unit vector, then there is an invariant vector at distance at most \(\min(\frac 12,\frac{\epsilon_{\mathcal S}}3)\) from it.
	Consider the two commuting inclusions \(\iota_{\mathcal R}: [\mathcal R]\to [\mathcal R\times\mathcal S]\) and \(\iota_{\mathcal S}:[\mathcal S]\to [\mathcal R\times\mathcal S]\), which are support-preserving for the projection on the first and second coordinate respectively.
	Let \(F\coloneqq \iota_{\mathcal R}(F_{\mathcal R})\cup \iota_{\mathcal S}(F_{\mathcal S})\) and \(\epsilon<\min(\frac{\epsilon_{\mathcal S}}3,\epsilon_{\mathcal R})\).
	We will show that \((F,\epsilon)\) is a Kazhdan pair for \(\mathcal R\times\mathcal S\).

	Let \([\pi]: [\mathcal R\times\mathcal S]\to \mathcal U(\mathcal H)\) be a full unitary representation with
	an \((F,\epsilon)\)-invariant unit vector \(\xi\), then  we can view \(\mathcal H\) both as an \(\LL^\infty(X)\)
	and as an \(\LL^\infty(Y)\)-module.
	Then \([\pi]\circ\iota_{\mathcal R}\) is a full unitary representation of \([\mathcal R]\) and \(\xi\) is \((F_{\mathcal R}, \epsilon_{\mathcal R})\)-invariant, so there is a \([\pi]\circ\iota_{\mathcal R}\)-invariant vector \(\eta\) at distance at most
	\(\min(\frac 12,\frac{\epsilon_{\mathcal S}}3)\) from \(\xi\), in particular \(\eta\) is not zero.
	By the triangle inequality, \(\eta\) is \((F,\epsilon+\frac{2\epsilon_{\mathcal S}}3)\) invariant.
	In particular, \(\eta\) is \((\iota_{\mathcal S}(F_{\mathcal S}),\epsilon_{\mathcal S})\)-invariant.

	Let \(\mathcal K\) denote the space of \([\pi]\circ\iota_{\mathcal R}\)-invariant vectors,
	note that \(\mathcal K\) is an \(\LL^\infty(Y)\)-module which is \([\pi]\circ\iota_{\mathcal S}\)-invariant.
	This module is non-zero by the previous paragraph and contains an \((F_{\mathcal S}),\epsilon_{\mathcal S})\)-invariant vector, so it contains a non-zero invariant vector for \([\pi]\circ\iota_{\mathcal S}\).
	We thus have found a non-zero vector which is both \([\pi]\circ\iota_{\mathcal R}\) and \([\pi]\circ\iota_{\mathcal S}\)
	invariant. By fullness and the fact that  \([\mathcal R\times\mathcal S]=[\iota_{\mathcal R}([\mathcal R])\cup \iota_{\mathcal S}([\mathcal S])]\), this vector is \([\pi]\)-invariant as wanted.
\end{proof}

Since \(\mathcal R^n\) is a class-bijective extension of \(\mathcal R^{\odot n}\), the following proposition is a consequence of the second item in \cite[Thm.~5.15]{anantharaman-delarocheCohomologyPropertyGroupoids2005} along with the previous lemma, but we provide a direct proof.

\begin{proposition}\label{prop: T for odot}
	Let \(\mathcal R\) be an ergodic p.m.p.~equivalence relation on \((X,\mu)\) and \(n \in \N\). Then \(\mathcal R^{\odot n}\) has property (T).
\end{proposition}
\begin{proof}
	Let \(\pi: \mathcal R^{\odot n}\to \mathcal U(\mathcal H)\) be a unitary representation with a sequence \((\xi_k)\)
	of almost invariant unit sections. Each \(\xi_k\) is thus an element of \(\LL^2(X^{\odot n},\mathcal H)\).
	We denote by \(\hat\pi: \mathcal R^n\to \mathcal U(\mathcal H)\) the unitary representation defined
	by \(\hat\pi((x_1,\dots,x_n),(y_1,\dots,y_n))=\pi([x_1,\dots,x_n],[y_1,\dots,y_n])\).
	Its space of square integrable sections is equal to \(\LL^2(X^n,\mathcal H)\), and the subspace of \(\mathfrak S_n\)-invariant sections naturally identifies to \(\LL^2(X^{\odot n},\mathcal H)\) : every \(\mathfrak S_n\)-invariant
	section \(X^n\to\mathcal H\) quotients down to a section \(X^{\odot n}\to \mathcal H\) and conversely
	every section \(\xi:X^{\odot n}\to \mathcal H\) yields a \(\mathfrak S_n\)-invariant section \(\hat \xi: X^n\to \mathcal H\) given by \(\hat\xi(x_1,\dots,x_n)=\xi([x_1,\dots,x_n])\).

	With these identifications in mind, the orthogonal projection onto
	\(\LL^2(X^{\odot n},\mathcal H)\) is  the map which takes \(\xi\in \LL^2(X^n,\mathcal H)\) to the average
	\(\frac 1{n!}\sum_{\sigma\in\mathfrak S_n} \sigma \xi,\)
	where \(\sigma\xi(x_1,\dots,x_n)=(x_{\sigma\inv(1)},\dots,x_{\sigma\inv(n)})\).
	Observe that \(\mathcal R^n\) has property (T) as a consequence of the previous lemma.
	The sequence \((\xi_k)\) is a sequence of almost invariant unit sections for \(\hat \pi\) and they belong to the subspace \(\LL^2(X^{\odot n},\mathcal H)\),
	hence by Proposition \ref{prop: quant (T)} there is \(\eta\in\LL^2(X^n,\mathcal H)\) which is \(\hat \pi\)-invariant, has norm \(1\)
	and is at distance \(<1\) from  \(\LL^2(X^{\odot n},\mathcal H)\).
	Its orthogonal projection \(\tilde \eta=\frac 1{n!}\sum_{\sigma\in\mathfrak S_n} \sigma \eta\) onto \(\LL^2(X^{\odot n},\mathcal H)\) is then non zero.
	Furthermore, observe that for every \(\sigma\in \mathfrak S_n,\), the section \(\sigma\eta\) is still
	\(\hat \pi\)-invariant.
	It follows that for every element \(([x_1,\dots,x_n],[y_1,\dots,y_n])\in\mathcal R^{\odot n}\), we have
	\begin{align*}
		\pi([x_1,\dots,x_n],[y_1,\dots,y_n])\tilde \eta(x_1,\dots,x_n) & = \frac 1{n!}\sum_{\sigma\in\mathfrak S_n} \hat\pi((x_1,\dots,x_n),(y_1,\dots,y_n))\sigma \eta(x_1,\dots, x_n) \\
		                                                               & =\frac 1{n!}\sum_{\sigma\in\mathfrak S_n} \sigma\eta(y_1,\dots,y_n)                                            \\
		                                                               & =\tilde \eta(y_1,\dots,y_n).
	\end{align*}
	So \(\tilde\eta\) is a non-zero invariant section for \(\pi\), which by ergodicity and rescaling yields
	an invariant unit section. We conclude that \(\mathcal R^{\odot n}\) has property (T) as wanted.
\end{proof}

\subsection{Proof of the direct implication in Theorem \ref{thmi:charaTfullgroups}}\label{sec:directimpl}

In this section, we put together the previous results so as to obtain the direct implication from Theorem \ref{thmi:charaTfullgroups}. 
We need a natural definition which mirrors Definition~\ref{df: strong ergo eq rel}.

\begin{definition}
	Let \(G\) be a group. Given a boolean action \(\rho:G\to\Aut(Y,\nu)\), a sequence \((A_n)\) of subsets of \(X\) is asymptotically \(\rho\)-invariant when for all \(g\in G\), \(\nu(\rho(g)A_n\bigtriangleup A_n)=0\), and  asymptotically trivial when \(\mu(A_n)(1-\mu(A_n))\to 0\). Finally \(\rho\) is \textbf{strongly ergodic} if every asymptotically \(\rho\)-invariant sequence is asymptotically trivial.
\end{definition}

\begin{theorem}
	Let \(\mathcal R\) be an ergodic p.m.p.~equivalence relation, suppose \(\mathcal R\) has property (T).
	Then every ergodic non-free boolean action of \([\mathcal R]\) is strongly ergodic.
\end{theorem}
\begin{proof}
	Let \(\rho: [\mathcal R]\to \Aut(Y,\nu)\) be p.m.p.~boolean ergodic non-free action of \([\mathcal R]\).
	By Theorem \ref{thm: main pmp equivalence relations} and ergodicity, there is some \(n\geq 1\) such that \(\rho\) comes from a p.m.p.~action of \(\mathcal R^{\odot n}\). Since \(\rho\) is ergodic, the action of \(\mathcal R^{\odot n}\) has to be ergodic
	as well, and since \(\mathcal R^{\odot n}\) has (T) the equivalence relation generated by \(\rho\)
	must have (T) as a consequence of Proposition \ref{prop: (T) action induces (T) eq rel}. It is thus strongly ergodic
	by Proposition \ref{prop: T implies st ergo}
	and we conclude that the action \(\alpha\) itself is strongly ergodic.
\end{proof}
\section{From strongly ergodic actions to property (T)}

In this section, we prove a version of the Connes-Weiss characterization of property (T) for groups,
but for p.m.p.~equivalence relations, allowing us to complete our proof of Theorem \ref{thmi:charaTfullgroups}.

\begin{theorem}
	Let \(\mathcal R\) be a p.m.p.~ergodic equivalence relation. Then \(\mathcal R\) has property (T)
	if and only if every ergodic p.m.p.~action of \(\mathcal R\) is strongly ergodic.
\end{theorem}

The direct implication was already used in the previous section,
and it is a direct consequence of Proposition \ref{prop: (T) action induces (T) eq rel} and Proposition \ref{prop: T implies st ergo}. 
Most of this section will thus be devoted to the converse, which is proved by
contraposition. Towards this, we 
follow the approach developped by Bekka-de la Harpe-Valette in 
\cite{bekkaKazhdanProperty2008}.
The thus need a
weaker version of property (T) in terms
of finite dimensional subrepresentations instead of invariant vectors \cite{bekkaKazhdanPropertyAmenable1993}. 
This is where affine actions come in, as in the group case.

\subsection{The Bekka-Valette characterization of property (T)}

\begin{definition}
	Let \(\mathcal R\) be a p.m.p.~equivalence relation and \(\pi\) be a unitary representation of \(\mathcal R\) on a Hilbert space \(\mathcal H\). \textbf{A finite dimensional subrepresentation of \(\pi\)} is a nontrivial \([\pi]\)-invariant subspace \(\mathcal K\) of \(\LL^2(X,\mu,\mathcal H)\) such that there exists a finite set \(F\subseteq \mathcal K\) with \(\LL^\infty(X,\mu) F\) dense in \(\mathcal K\).
\end{definition}

This section is devoted to the proof of the following theorem, which is the p.m.p.\ equivalence relation
version of the Bekka-Valette theorem \cite[Thm.~1]{bekkaKazhdanPropertyAmenable1993}.

\begin{theorem} \label{thm:stronger version of (T)}
	Let \(\mathcal R\) be an ergodic p.m.p.~equivalence relation. Then \(\mathcal R\) has property (T) if and only if every unitary representation of \(\mathcal R\) which has almost invariant unit sections contains a finite dimensional subrepresentation.
\end{theorem}

Note that the latter condition is a weaker version of property (T). Thus only right-to-left direction remains to be proven.

\begin{proposition}[{\cite[Chap.~4]{delaharpeProprieteKazhdanPour1989}}] \label{prop:construction of representations for affine actions}
	Let \(\mathcal H\) be a real affine Hilbert space.
	Then for any \(t>0\), there exists a unique complex Hilbert space \(\mathcal H_t\) and a continuous mapping \(\xi \mapsto \xi_t\) from \(\mathcal H\) to the unit sphere of \(\mathcal H_t\) such that for any \(\xi,\eta \in \mathcal H\), \(\langle\xi_t,\eta_t \rangle = \exp (-t \Vert \xi - \eta \Vert^2)\) and the image of \(\{\xi_t \colon \xi \in \mathcal H\}\) is total in \(\mathcal H_t\).
\end{proposition}

Let \(\mathcal R\) be a \pmp equivalence relation.
An \textbf{affine isometric action} of \(\mathcal R\) is a homomorphism \(\rho:\mathcal R\to \Iso(\mathcal H)\), where \(\Iso(\mathcal H)\) denotes the group of affine isometries of a separable real Hilbert space \(\mathcal H\).

For \(t>0\), we define a homomorphism \(\phi \mapsto \phi_t\) from the group of affine isometries \(\mathrm{Iso}(\mathcal H)\) of \(\mathcal H\) to the unitary group \(\mathcal U(\mathcal H_t)\) of \(\mathcal H_t\), 
where \(\phi_t\) is the only element of \(\mathcal U(\mathcal H_t)\) such that \(\forall \xi \in \mathcal H, \phi_t(\xi_t)=(\phi(\xi))_t\). 
This homomorphism is continuous if we endow both $\Iso(\mathcal H_t)$ and $\mathcal U(\mathcal H_t)$ with the topology 
of pointwise convergence for the Hilbert space norm. 
To see this, first recall that this topology coincides with the topology of pointwise weak convergence on $\mathcal U(\mathcal H_t)$.	
Then let \((\phi_n)\) be a sequence in \(\mathrm{Iso}(\mathcal H)\) converging to \(\phi\), for all \(\xi, \eta \in \mathcal H\), we have \(\langle (\phi_n)_t(\xi_t), \eta_t \rangle = \langle (\phi_n(\xi))_t, \eta_t \rangle = \exp(-t \Vert \phi_n(\xi) - \eta \Vert^2) \underset{n \rightarrow \infty}{\longrightarrow} \exp(-t \Vert \phi(\xi) - \eta \Vert^2)\).

Therefore, any affine isometric action \(\rho\) of an equivalence relation \(\mathcal R\) on \(\mathcal H\) gives rise to a unitary representation \(\rho_t\) on \(\mathcal H_t\).
We now give a cocycle-free proof of some results from \cite{anantharaman-delarocheCohomologyPropertyGroupoids2005}.

\begin{lemma}[{\cite[Lem.~3.20]{anantharaman-delarocheCohomologyPropertyGroupoids2005}}] \label{lem: bounded everywhere}
	Let \(\mathcal R\) be an ergodic p.m.p.~equivalence relation. Let \(\rho\) be an affine isometric action of \(\mathcal R\) on a real Hilbert space \(\mathcal H\), suppose that there is \(A\subseteq X\) of positive measure and a Borel section \(\xi:X\to\mathcal H\) such that for all \(x\in A\), we have
	\[\sup_{y\in\mathcal R_x \cap A}\norm{\rho(x,y)\xi(x)-\xi(y)}<+\infty.\]
	Then there exists a section \(\xi' \colon X \to \mathcal H\) such that for almost all \(x\in X\), we have
	\[\sup_{y\in\mathcal R_x}\norm{\rho(x,y)\xi'(x)-\xi'(y)}<+\infty.\]
\end{lemma}
\begin{proof}
	Fix a sequence \((T_n)\) in \([\mathcal R]\) such that the graphs of the \(T_n\) cover \(\mathcal R\). We define a Borel map \(\theta \colon X \to X\) by letting \(\theta \colon x \mapsto T_nx\) where \(n\) is the least integer such that \(T_nx \in A\). Such an integer exists by ergodicity of \(\mathcal R\).

	Now let \(\xi' \colon x \mapsto \rho(\theta x,x)\xi(\theta x)\). For \((x,y) \in \mathcal R\), we compute
	\begin{align*}
		\norm{\rho(x,y)\xi'(x) - \xi'(y)} & = \norm{\rho(x,y)\rho(\theta x,x)\xi(\theta x) - \rho(\theta y,y)\xi(\theta y)} \\
		                                  & = \norm{\rho(\theta x,y)\xi(\theta x) - \rho(\theta y,y)\xi(\theta y)}          \\
		                                  & = \norm{\rho(y,\theta y)\rho(\theta x,y)\xi(\theta x) - \xi(\theta y)}          \\
		                                  & = \norm{\rho(\theta x,\theta y)\xi(\theta x) - \xi(\theta y)}
	\end{align*}
	which is bounded when \(y\) goes through \(\mathcal R_x\).
\end{proof}

\begin{lemma}[{\cite[Thm.~3.19]{anantharaman-delarocheCohomologyPropertyGroupoids2005}}] \label{lem: bounded section implies fixed section}
	Let \(\mathcal R\) be an ergodic p.m.p.~equivalence relation. 
	Let \(\rho\) be an affine isometric action of \(\mathcal R\) on a real Hilbert space \(\mathcal H\) 
	and let \(\xi: X\to\mathcal H\) be a Borel section, suppose that
	for all \(x\in X\), we have \(\sup_{y\in\mathcal R_x}\norm{\rho(x,y)\xi(x)-\xi(y)}<+\infty\). Then \(\rho\) admits an invariant section.
\end{lemma}
\begin{proof}
	For all \((x,y)\in\mathcal R\), since \(\rho(y,x)\) is an isometry which is the inverse of \(\rho(x,y)\), we have
	\(\sup_{y\in\mathcal R_x}\norm{\rho(y,x)\xi(y)-\xi(x)}<+\infty\).
	Using a countable group \(\Gamma\) such that \(\mathcal R=\mathcal R_\Gamma\), we see that the section \(\eta\)
	which takes \(x\in X\) to the circumcenter of the closed convex hull of
	\(\{\rho(y,x)\xi(y): y\in \mathcal R_x\}\) is Borel.
	Such a section is easily checked to be invariant by \(\rho\) since \(\rho(x',x)\) takes \(\rho(y,x')\xi(y)\) to
	\(\rho(y,x)\xi(y)\).
\end{proof}

\begin{proposition}\label{prop:invariant sections affine vs unitary}
	Let \(\mathcal R\) be an ergodic p.m.p.~equivalence relation. 
	Let \(\rho\) be an affine isometric action of \(\mathcal R\) on a real Hilbert space \(\mathcal H\) and \(t>0\). 
	Then \(\rho\) admits an invariant section if and only if the unitary representation \(\rho_t\) has an invariant unit section.
\end{proposition}

\begin{proof}
	If \(\xi\) is an invariant section for \(\rho\) then it is straightforward that \(\xi_t \colon x \mapsto \xi(x)_t\) is an invariant unit section for \(\rho_t\).

	Conversely, if there exists an invariant unit section \(\zeta\) for \(\rho_t\), suppose by contradiction that \(\rho\) does not admit any invariant section.
	As in the proof of Theorem 4.12 in \cite{anantharaman-delarocheCohomologyPropertyGroupoids2005}, we fix a dense sequence \((h_n)_{n \geq 1} \in \mathcal H\) and for \(m,n \in \N^*\), let \(E_{m,n} = \{x \in X \colon |\langle (h_m)_t, \zeta(x) \rangle| \geq \frac{1}{n} \}\). Since the span of the \((h_m)_t\) for \(m \geq 1\) is dense in \(\mathcal H_t\) and \(\zeta\) is a unit section, we have \(\bigcup_{m,n \geq 1} E_{m,n} = X\) so there exist \(m,n \in \N^*\) such that \(\mu(E_{m,n}) >0\).

	Applying Lemma \ref{lem: bounded everywhere} and Lemma \ref{lem: bounded section implies fixed section} to the constant section \(x \mapsto h_m\), we get that for almost any \(x \in E_{m,n}\),
	\[\sup_{y \in \mathcal R_x \cap E_{m,n}} \norm{\rho(x,y) h_m - h_m} = +\infty.\]
	Therefore, for almost any \(x \in E_{m,n}\), we can fix a sequence \((y_k)\) in \(\mathcal R_x \cap E_{m,n}\) such that \((\norm{\rho(x,y_k) h_m - h_m})_k\) tends to \(+ \infty\). In particular the sequence of vectors \((\rho(x,y_k) h_m)_k\) tends to \(\infty\) in \(\mathcal H\) as \(k\to+\infty\).

	Now for every \(\xi \in \mathcal H\) we have
	\[|\langle (\rho(x,y_k) h_m)_t, \xi_t \rangle| = \exp(-t \norm{\rho(x,y_k) h_m -\xi}^2) \underset{k \rightarrow \infty}{\longrightarrow} 0.\]
	Since the image of \(\mathcal H\) is total in \(\mathcal H_t\), we conclude that the sequence \(\left( (\rho(x,y_k) h_m)_t \right)_{k \in \N}\) weakly converges to \(0\) in \(\mathcal H_t\).
	But on the other hand, for every \(k \in \N\) we have
	\[|\langle (\rho(x,y_k)h_m)_t , \zeta(x) \rangle| = |\langle (h_m)_t, \rho_t(y_k,x)\zeta(x) \rangle| = |\langle (h_m)_t, \zeta(y_k) \rangle| \geq \frac{1}{n}.\]
	This contradicts the fact that \(\left( (\rho(x,y_k) h_m)_t \right)_{k \in \N}\) weakly converges to \(0\), which finishes the proof.
\end{proof}

\begin{lemma} \label{lem:sum of arbitrarily small representations has almost invariant vectors}
	Let \(\rho\) be an affine isometric action of a p.m.p.~equivalence relation \(\mathcal R\) on a real Hilbert space \(\mathcal H\). Consider the family \((\mathcal H_t, \rho_t)_t\) associated to \(\rho\) as in Proposition \ref{prop:construction of representations for affine actions}. Let \((t_n)_{n \in \N}\) be a sequence of positive reals converging to \(0\). Then \(\bigoplus_{n \in \N} \rho_{t_n}\) has almost invariant unit sections.
\end{lemma}

\begin{proof}
	For \(n \in \N\), let \(\xi_n\) be the constant section with value \(0_{t_n}\).
	Then \(\xi_n\) is a unit section and moreover for every \((x,y) \in \mathcal R\),
	\begin{align*}
		\Vert \rho_{t_n}(x,y) \xi_n (x) - \xi_n (y) \Vert^2 &
		= \Vert (\rho(x,y)(0))_{t_n} - 0_{t_n} \Vert^2      &                                                                                                                  \\
		                                                    & = \Vert (\rho(x,y)(0))_{t_n} \Vert^2 + \Vert 0_{t_n} \Vert^2 - 2 \langle (\rho(x,y)(0))_{t_n}, 0_{t_n} \rangle & \\
		                                                    & = 2 \left(1 - \exp(-t_n \Vert \rho(x,y)(0) \Vert^2) \right).                                                    &
	\end{align*}
	It easily follows that \((\xi_n)\), seen as a sequence of sections in \(\bigoplus_{n \in \N} \mathcal H_{t_n}\) is a sequence of almost invariant unit sections for \(\bigoplus_{n \in \N} \rho_{t_n}\).
\end{proof}

\begin{lemma} \label{lem:representations of the square of affine actions}
	Let \(\rho\) be an affine isometric action of a p.m.p.~equivalence relation \(\mathcal R\) on a real Hilbert space \(\mathcal H\). Let us denote \(\rho\oplus \rho\) the diagonal action on \(\mathcal H \oplus \mathcal H\). Then for \(t>0\), \((\rho\oplus \rho)_t = \rho_t \otimes \rho_t\).
\end{lemma}

\begin{proof}
	Fix \(t>0\). We define a map \(\Psi_t \colon \mathcal H \oplus \mathcal H \to \mathcal H_t \otimes \mathcal H_t\) by the formula \(\Psi_t(\xi,\eta) \coloneqq \xi_t \otimes \eta_t\).

	First, for \(\xi,\xi',\eta,\eta' \in \mathcal H\) we have
	\begin{align*}
		\langle \Psi_t(\xi\oplus \eta), \Psi_t(\xi'\oplus\eta') \rangle & = \exp(-t \Vert \xi - \xi' \Vert^2) \exp(-t \Vert \eta - \eta' \Vert^2) & \\
		                                                                & = \exp(-t \Vert (\xi\oplus\eta) - (\xi'\oplus\eta') \Vert^2)            &
	\end{align*}
	and it is clear that the image of \(\Psi_t\) is total in \(\mathcal H_t \otimes \mathcal H_t\). Moreover, for
	\begin{align*}
		\Psi_t(\rho\oplus \rho(x,y)(\xi\oplus\eta)) & = \Psi_t(\rho(x,y)\xi\oplus\rho(x,y)\eta)              & \\
		                                            & = (\rho(x,y)\xi)_t \otimes (\rho(x,y)\eta)_t           & \\
		                                            & = \rho_t(x,y) \xi_t \otimes \rho_t(x,y) \eta_t         & \\
		                                            & = (\rho_t \otimes \rho_t)(x,y) (\xi_t \otimes \eta_t)  & \\
		                                            & = (\rho_t \otimes \rho_t)(x,y) \Psi_t(\xi\oplus \eta). &
	\end{align*}
	By uniqueness of the construction in Proposition \ref{prop:construction of representations for affine actions}, 
	it follows that \(\Psi_t(\xi\oplus\eta) = (\xi\oplus\eta)_t\) and \((\rho\oplus \rho)_t = \rho_t \otimes \rho_t\).
\end{proof}

Recall that given a complex Hilbert space \(\mathcal H\),
we get a conjugate Hilbert space \(\overline{\mathcal  H}\) which is as a set equal to \(\mathcal H\), but endowed with a new structure where scalar multiplication by a complex number \(c\in\C\) becomes \(\xi\mapsto \bar c\xi\) and the scalar product  becomes
\(\la \xi,\eta\ra_{\overline{\mathcal H}}=\la \eta,\xi\ra=\overline{\la \xi,\eta\ra}\).
Given a unitary representation \(\pi\) on \(\mathcal H\), we thus get a unitary representation \(\bar \pi\) on \(\overline{\mathcal H}\). Note that its invariant vectors are the same as those of \(\pi\).

\begin{lemma}[\cite{gardellaComplexityConjugacyOrbit2017}, Lemma 3.16.(2)] \label{lem:invariant vectors of tensor}
	Let \(\mathcal R\) be a p.m.p equivalence relation and let \(\pi\) be a unitary representation of \(\mathcal R\). The following assertions are equivalent:
	\begin{enumerate}
		\item \(\pi \otimes \overline{\pi}\) has a non-zero invariant  section,
		\item there is a unitary representation \(\sigma\) of \(\mathcal R\) such that \(\pi \otimes \sigma\) has a non-zero invariant  section,
		\item \(\pi\) contains a finite dimensional subrepresentation.
	\end{enumerate}
\end{lemma}

Finally, let us introduce the affine version of property (T), following Anantharaman-Delaroche:

\begin{definition}[{Property (FH), \cite{anantharaman-delarocheCohomologyPropertyGroupoids2005}, Section 4.3}]
	A p.m.p.~equivalence relation \(\mathcal R\) is said to have property (FH) if every affine isometric action of \(\mathcal R\) on a separable real Hilbert space admits an invariant section.
\end{definition}

\begin{thm}[{\cite[Theorems 4.8 and 4.12]{anantharaman-delarocheCohomologyPropertyGroupoids2005}}]\label{thm: T eq FH}
	Let \(\mathcal R\) be an ergodic \pmp equivalence relation. Then \(\mathcal R\) has property (T) if and only if \(\mathcal R\) has property (FH).
\end{thm}

We are now ready to prove Theorem \ref{thm:stronger version of (T)}.

\begin{proof}
	Recall that left-to-right direction is trivial. For the other direction, suppose that every unitary representation of \(\mathcal R\) which has almost invariant unit sections contains a finite dimensional subrepresentation and let us show that \(\mathcal R\) has property (FH). We will use repeatedly Lemma \ref{lem:invariant vectors of tensor}.

	Let \(\rho\) be an affine isometric action of \(\mathcal R\) on a separable real Hilbert space \(\mathcal H\). Let us consider the representation \(\pi \coloneqq \bigoplus_{n \in \N^*} \rho_{1/n}\). Then \(\pi\) admits almost invariant unit sections by Lemma \ref{lem:sum of arbitrarily small representations has almost invariant vectors}, and therefore \(\pi\) contains a finite dimensional subrepresentation, or in other words, \(\pi \otimes \overline{\pi}\) has an non-zero invariant  section.

	But \(\pi \otimes \overline{\pi} = \bigoplus_{n,m \in \N^*} \rho_{1/n} \otimes \overline{\rho_{1/m}}\)  so there must be \(n_0,m_0 \in \N^*\) such that we have a non-zero invariant section for \(\rho_{1/n_0} \otimes \overline{\rho_{1/m_0}}\). Then \(\rho_{1/n_0}\) contains a finite dimensional subrepresentation and finally \(\rho_{1/n_0} \otimes \overline{\rho_{1/n_0}}\) has a non-zero invariant  section, which by the definition of \(\overline{\rho_{1/n_0}}\) is equivalent to \(\rho_{1/n_0} \otimes \rho_{1/n_0}\) admitting a non-zero invariant  section.

	Applying Lemma \ref{lem:representations of the square of affine actions} and Proposition \ref{prop:invariant sections affine vs unitary} we get that \(\rho \oplus \rho\) has an invariant section. Therefore \(\rho\) has an invariant section and \(\mathcal R\) has property (FH).
\end{proof}

\begin{remark}
	We have developed almost all the tools needed to prove that 
	(T) is equivalent to (FH) (Theorem~\ref{thm: T eq FH}).
	Indeed, that (T) implies (FH) for \pmp equivalence relations follows 
	from a simpler version of the above proof of Theorem~\ref{thm:stronger version of (T)}: 
	given an arbitrary affine isometric action $\rho$ on $\mathcal H$, the 
	direct sum $\bigoplus_{n\geq 1} \rho_{1/n}$ has almost invariant unit sections
	by Lemma~\ref{lem:sum of arbitrarily small representations has almost invariant vectors}, thus it has a non zero invariant section. 
	One of its orthogonal projection on $\mathcal H_{1/n}$ is a non 
	zero invariant section, and thus $\rho$ admits an invariant section
	by Proposition~\ref{prop:invariant sections affine vs unitary}.

	Conversely, one needs to use the fact that an isometric action 
	is the same data as a unitary representation $\pi$ plus a $\pi$-cocycle, 
	and that this cocycle is a coboundary if and only if the isometric action 
	has a fixed point. 
	Now suppose by contradiction $\mathcal R$ has (FH) but not (T), 
	let $\pi$ be a unitary representation with almost invariant sections 
	but no non-zero invariant section. 
	The absence of non-zero invariant sections for $\pi$ 
	implies that the map 
	which take a section $\xi$ to the corresponding coboundary 
	$\beta(\xi):(x,y)\in\mathcal R \mapsto \pi(x,y)\xi(y)-\xi(x)$
	is injective. 
	It is moreover a continuous linear map
	if we endow both the space of unit sections and the space of cocycles with 
	the Polish (in particular, Fréchet) topology of convergence in measure, 
	thus it is a homeomorphism onto its image by the closed 
	graph theorem. 
	However, if we now consider a 
	sequence of almost invariant unit sections 
	$(\xi_n)$, then $\beta(\xi_n)\to 0$ in measure but $\xi_n\not\to 0$
	in measure as each $\xi_n$ is a unit section, a contradiction.
	For details, see \cite[Thm.~4.8]{anantharaman-delarocheCohomologyPropertyGroupoids2005}.
\end{remark}
\subsection{A reminder on Gaussian actions} \label{sec:gaussian actions}

We now briefly present the classical construction of Gaussian actions associated to unitary representations,
which given our definitions of unitary representations and of boolean actions
of p.m.p.~equivalence relations (Definition \ref{df: various actions of eqrel}) work exactly as in the group case.
This can be phrased as the construction of a homomorphism from the orthogonal group of a separable real Hilbert space to \(\Aut(Z,\eta)\), but the details of this construction are important for the proofs.
Here we simply state the facts which are relevant to us; the unacquainted reader should first read \cite[Appendix E]{kerrErgodicTheoryIndependence2016}.

First we define the symmetric Fock space of a real Hilbert space \(\mathcal H\). Let \(n \in \N\) and \(\sigma \in \mathfrak S_n\). Then \(\sigma\) induces an orthogonal operator \(O_\sigma\) on \(\mathcal H^{\otimes n}\) defined on simple tensors by \(O_\sigma(h_1 \otimes \dots \otimes h_n) = h_{\sigma^{-1}(1)} \otimes \dots \otimes h_{\sigma^{-1}(n)}\).
We let the \textbf{n-th symmetric space of \(\mathcal H\)} be the subspace of \(\mathcal H^{\otimes n}\) of elements that are fixed by all \(O_\sigma\) for \(\sigma \in \mathfrak S_n\),
and we denote it by \(\mathcal H^{\odot n}\).
For \(n=0\), we take the convention that \(\mathcal H^{\odot 0} = \R\).
The \textbf{symmetric Fock space of \(\mathcal H\)} is the direct sum \(S(\mathcal H) \coloneqq \bigoplus_{n \in \N} \mathcal H^{\odot n}\).

Any orthogonal representation \(\pi\) on \(\mathcal H\) then induces orthogonal representations \(\pi^{\odot n}\) on \(\mathcal H^{\odot n}\) and thus an orthogonal representation \(S(\pi)\coloneqq\bigoplus_{n\in\N}\pi^{\odot n}\) on \(S(\mathcal H)\). These constructions also make sense in the complex case for unitary representations, and we will use the same notation.

Now let \(\mathcal R\) be an ergodic p.m.p.~equivalence relation on \((X,\mu)\) and let \(\pi\) be a unitary representation of \(\mathcal R\) on a Hilbert space \(\mathcal H\).
Then one can associate to \(\pi\) a \pmp boolean  action \(\alpha_\pi\) on a standard atomless probability space \((Z,\eta)\) called the \textbf{Gaussian action associated to \(\pi\)}  so that its Koopman representation on the orthogonal of constant functions \(\kappa^0_{\alpha_\pi}\) contains \(\pi\) as a direct summand as follows.
One starts by considering \(\pi\) as an orthogonal representation \(\pi_\R\) by viewing \(\mathcal H\) as a real Hilbert space.

One constructs a \emph{Gaussian Hilbert space} \(\widetilde{\mathcal H} \subset \LL^2(Z,\eta,\R)\) of dimension \(\mathrm{dim}\ \mathcal H\) that does not contain the constant functions, such that the symmetric Fock space \(S(\widetilde{\mathcal H})\) of \(\widetilde{\mathcal H}\) is isometrically isomorphic to \(\LL^2(Z,\eta,\R)\) via an isomorphism sending \(\widetilde{\mathcal H}^{\odot 0}\) on \(\R \1_Z\) (see \cite[Sec.~E.3]{kerrErgodicTheoryIndependence2016}).
By the means of fixed isometric isomorphisms, let us identify \(\widetilde{\mathcal H}\) and \(\mathcal H\) as well as \(S(\widetilde{\mathcal H})\) and \(\LL^2(Z,\eta,\R)\).
We can then consider \(S(\pi_\R)\) as an orthogonal representation on \(\LL^2(Z,\eta,\R)\).

The construction of \(S(\pi_\R)\) then ensures that \(S(\pi_\R)\) preserves the subset \(\mathrm{MAlg}(Z,\eta) \subseteq \LL^2(Z,\eta,\R)\) and acts by automorphism on it. Call \(\alpha_\pi\) the co-restriction of \(S(\pi_\R)\) to \(\mathrm{Aut}(Z,\eta)\). Then \(\alpha_{\pi}\) is a Borel homomorphism and so it is a boolean measure-preserving \(\mathcal R\)-action on \((Z,\eta)\).

Furthermore, by definition, the orthogonal Koopman representation of \(\alpha_\pi\) on \(\LL^2(Z,\eta,\R)\) is \( S(\pi_\R)\) . Finaly the unitary Koopman representation of \(\alpha_\pi\) on \(\LL^2(Z,\eta,\C)\) is the complexification of the orthogonal Koopman representation. Since the complexification of \(\pi_\R\) is \(\pi\oplus\bar \pi\) (see \cite[Prop.~E.1]{kerrErgodicTheoryIndependence2016}), we get that
\[\kappa_{\alpha_\pi}=S(\pi\oplus\bar \pi).\]
Moreover, the Koopman representation on the orthogonal of constant functions is  \(\kappa^0_{\alpha_\pi} = \bigoplus_{n \in \N^*} (\pi\oplus\bar \pi)^{\odot n}\), in particular \(\pi\) is a direct summand of \(\kappa^0_{\alpha_\pi}\).


\subsection{The Connes-Weiss characterization of property (T)}

Now we can prove the following theorem, which extends a theorem of Connes and Weiss from discrete countable groups to \pmp equivalence relation \cite{connesPropertyAsymptoticallyInvariant1980}. A version for \pmp groupoids can be found in the thesis of the second-named author \cite{giraudClassificationActionsErgodiques2018}.

\begin{theorem}[Connes-Weiss for \pmp equivalence relations] \label{thm:connes-weiss}
	Let \(\mathcal R\) be a p.m.p.~equivalence relation. Then the following assertions are equivalent:
	\begin{enumerate}[(i)]
		\item \label{item: T}\(\mathcal R\) has property (T).
		\item \label{item: erg is strerg}Every ergodic measure-preserving \(\mathcal R\)-action is strongly ergodic.

	\end{enumerate}
\end{theorem}

\begin{proof}
	The implication  \eqref{item: T}\(\Rightarrow\)\eqref{item: erg is strerg} is a direct consequence of  Proposition \ref{prop: T implies st ergo} and Proposition \ref{prop: (T) action induces (T) eq rel}. We thus have to show \eqref{item: erg is strerg}\(\impl\)\eqref{item: T}, which we do via the contrapositive.

	Let us thus suppose that \(\mathcal R\) does not have property (T) and construct an ergodic action which is not strongly ergodic.

	By Theorem \ref{thm:stronger version of (T)}, there is a unitary representation \(\pi\) of \(\mathcal R\) which has almost invariant unit sections but does not admit any finite dimensional subrepresentation. We construct the Gaussian action \(\alpha_\pi\)
	on \( (Z,\eta) \)
	associated to \(\pi\) (see Section \ref{sec:gaussian actions}). Since \(\pi\) does not admit any finite dimensional subrepresentation,  neither does \(\pi\oplus\bar \pi\), and hence neither do \((\pi\oplus\bar \pi)^{\otimes n}\) for \(n \in \N^*\) by Lemma \ref{lem:invariant vectors of tensor}. In particular, for \(n \in \N^*\), \((\pi\oplus\bar \pi)^{\odot n}\) does not admit any finite dimensional subrepresentation. It follows again from Lemma \ref{lem:invariant vectors of tensor} that for any \(n,m \in \N^*\), \((\pi\oplus\bar \pi)^{\odot n} \otimes \overline{(\pi\oplus\bar \pi)^{\odot m}}\) does not have a non-zero invariant section.

	So it follows from the above description that the Koopman representation on the orthogonal of constant functions \(\kappa^0_{\alpha_\pi}\) of \(\mathcal R\) satisfies
	\(\kappa^0_{\alpha_\pi} \otimes \overline{\kappa^0_{\alpha_\pi}} = \bigoplus_{n,m \in \N^*} (\pi\oplus\bar \pi)^{\odot n} \otimes \overline{(\pi\oplus\bar \pi)^{\odot m}}\).
	Applying Lemma \ref{lem:invariant vectors of tensor} once more, we see that \(\kappa^0_{\alpha_\pi}\) has no finite dimensional subrepresentation, in particular it has no non-zero invariant section.

	Now the Koopman representation of the boolean action \([\alpha_\pi]\) of \([\mathcal R]\) on \(\LL^2(X\times Z,\mu\otimes\eta)\) is equal to the direct sum of the standard representation 
	\(\kappa_\iota\) on \(\LL^2(X,\mu)\) (corresponding to fiberwise constant functions) and \([\kappa^0_{\alpha_\pi}]\) (corresponding to functions whose integral over each fiber is zero). 
	By the previous paragraph, the latter has no non-zero invariant vectors. 
	It then follows from the ergodicity of \(\mathcal R\) that the only invariant vectors of our direct sum are given by constant functions in \(\LL^2(X\times Z,\mu\otimes\eta)\), 
	so the \(\mathcal R\)-action \(\alpha_\pi\) is ergodic\footnote{This argument can actually be upgraded to show that the action is weakly mixing, see \cite[Sec.~3.10]{gardellaComplexityConjugacyOrbit2017} for the definition.} by Proposition \ref{prop:ergo via L2}. \\

	We then show that the action is not strongly ergodic by constructing a sequence of asymptotically \([\alpha_\pi]\)-invariant Borel sets \(A_n \subset X \times Z\) which is not asymptotically trivial.

	Let \((\xi_n)\) be a sequence of almost invariant unit sections for \(\pi\) viewed as an orthogonal representation. Recall that through the construction of the Gaussian action \(\alpha_{\pi}\), we identified \(\mathcal H\) to a Gaussian Hilbert space, so for all \(x \in X\), \(\xi_n(x,\cdot)\) is a centered gaussian random variable of variance \(1\). Write \(\xi_n^x\) for \(\xi_n(x,\cdot)\) and for \(x \in X\), let \(A_n^x \subset Z\) be the set \(\{z\in Z \colon \xi_n^x(z) \geq 0\}\). Since \(\xi_n^x\) is centered, we have \(\nu(A_n^x)=\frac{1}{2}\).

	Fix \((x,x') \in \mathcal R\). We write \(\pi(x,x')\xi_n^x = \cos \theta_{n,x,x'} \xi_n^{x'} + \sin \theta_{n,x,x'} \eta_{n,x,x'}\) in \(\mathcal H\), for some \(\theta_{n,x,x'} \in [0,\pi]\) and \(\eta_{n,x,x'} \in \left( \xi_n^{x'} \right) ^\perp\).
	Then \(\xi_n^{x'}\) and \(\eta_{n,x,x'}\) are orthogonal gaussian random variables of variance \(1\) so they are independent and their joint distribution is a probability measure \(m\) on \(\R^2\) which is the product of two gaussian measures (in particular, it is rotation invariant). Last, \((\xi_n)\) is a sequence of almost invariant sections so for \(n\) large enough, \(\theta_{n,x,x'} \in [0, \frac{\pi}{2})\).
	It follows that
	\begin{align*}
		\eta(\alpha_\pi(x,x') A_n^x \bigtriangleup A_n^{x'})
		=\	 & \eta \big( \{z\in Z\colon \pi(x,x')\xi_n^x(z) \geq 0 \text{ and } \xi_n^{x'}(z) <0\}                                   \\
		    & \cup \{z\in Z\colon \pi(x,x')\xi_n^x(z) <0 \text{ and } \xi_n^{x'}(z) \geq 0\} \big)                                   \\
		=\	 & m \left( \{(a,b) \in \R^2 \colon \cos \theta_{n,x,x'} a + \sin \theta_{n,x,x'} b \geq 0 \text{ and } a<0\} \right)     \\
		    & + m \left( \{(a,b) \in \R^2 \colon \cos \theta_{n,x,x'} a + \sin \theta_{n,x,x'} b < 0 \text{ and } a \geq 0\} \right)
	\end{align*}
	Using Figure \ref{fig: angle} and the fact that these inequalities are about the signs of the scalar product of \((a,b)\) with \((\cos \theta_{n,x,x'},\sin\theta_{n,x,x'})\) and \((1,0)\) respectively, we see that the first term is equal to the measure of the upper left gray part, while the second term is equal to the measure of the lower right gray part.

	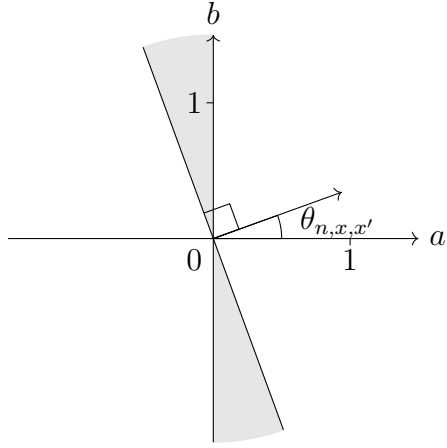
\begin{figure}[htpb]
		\begin{center}
			\begin{tikzpicture}[scale=1.8]
				\filldraw[fill=gray!20,draw=gray!20] (0,0) -- (0,1.5) arc (90:110:1.5cm) -- cycle;
				\filldraw[fill=gray!20,draw=gray!20] (0,0) -- (0,-1.5) arc (-90:-70:1.5cm) -- cycle;
				\draw[->] (-1.5,0) -- (1.5,0) node[anchor=west] {\(a\)};
				\draw[->] (0,-1.5) -- (0,1.5) node[anchor=south] {\(b\)};
				\draw (0,0) node[anchor=north east] {\(0\)};
				\draw (1,-0.05) -- (1,-0) node[anchor= north] {\(1\)};
				\draw (-0.05,1) -- (0,1) node[anchor= east] {\(1\)};
				\draw[->] (0,0) -- (20:1cm);
				\draw (0,0) --(110:1.5cm);
				\draw (0,0) --(-70:1.5cm) ;
				\draw arc (0:20:1cm);
				\draw (0.91,0.12) node {\(\theta_{n,x,x'}\)};
				\draw (0,0) -- (0.5,0) arc (0:20:5mm) -- cycle;
				\draw[rotate=20] (0,0) -- (0.2,0) -- (0.2,0.2) -- (0,0.2);
			\end{tikzpicture}
		\end{center}
		\caption{The \(m\)-measure of the gray part is equal to \(\eta(\alpha_\pi(x,x')A_n^x \bigtriangleup A_n^{x'})\)}\label{fig: angle}
	\end{figure}
	It then follows from rotation invariance that
	\[ \nu(\alpha_\pi(x,x') A_n^x \bigtriangleup A_n^{x'})=\frac{\theta_{n,x,x'}}\pi.\]

	Let now \(A_n = \{(x,z) \in X \times Z \colon z \in A_n^x\}\), note that \(A_n\) is measurable as a consequence of Lemma \ref{proposition:l0measurable}. Then \(\mu \otimes \eta(A_n) = \frac{1}{2}\) so the sequence \((A_n)\) is not asymptotically trivial.
	Let \(T \in [\mathcal R]\), then
	\begin{align*}
		\mu \otimes \eta([\alpha_\pi](T) A_n \bigtriangleup A_n) & = \int_X \eta(\alpha_\pi(T^{-1}x,x) A_n^{T^{-1}x} \bigtriangleup A_n^x)\ d\mu(x) \\
		                                                         & = \frac{1}{\pi} \int_X \theta_{n,T^{-1}x,x}\ d\mu(x)
	\end{align*}
	which converges to \(0\) by dominated convergence theorem. Therefore \((A_n)\) is asymptotically invariant.
	We conclude that \(\alpha_\pi\) is not strongly ergodic, although it is ergodic as shown before. Note that \(\alpha_\pi\) was constructed as a boolean \(\mathcal R\)-action, but by Proposition \ref{prop:boolean yields pmp action} we can equally view it as a \pmp \(\mathcal R\)-action as defined in Section \ref{sec: boolean action from eq rel}. This finishes the proof.
\end{proof}

\subsection{Proof of the indirect implication in Theorem \ref{thmi:charaTfullgroups}}

\begin{theorem}
	Let \(\mathcal R\) be a p.m.p equivalence relation. Suppose that all non-free ergodic boolean actions of \([\mathcal R]\) are strongly ergodic. Then \(\mathcal R\) has property (T).
\end{theorem}

\begin{proof}
	We use our characterization of property (T) obtained in Theorem \ref{thm:connes-weiss}. Let \(\pi:(Y,\nu)\to(X,\mu)\) be an ergodic  measure-preserving \(\mathcal R\)-action on \((Y,\nu)\), i.e.~\(\pi\) is a class-bijective extension of \(\mathcal R\) by a \pmp ergodic equivalence relation \(\mathcal S\) on \((Y,\nu)\) . Then the corresponding boolean action \(\rho_\pi\) of \([\mathcal R]\) (see Definition \ref{lemma: def rho pi}) is ergodic  and far from being free. Indeed, by Corollary \ref{cor:every set is some support} we can take \(T \in [\mathcal R]\) such that \(\mu(\supp T) \notin \{0,1\}\), and then \(\supp [\alpha](T) = \pi\inv(\supp T)\) has measure \(\mu(\supp T)\).

	By our assumption the non-free ergodic boolean action \(\rho_\pi\) is strongly ergodic, which by definition means that the equivalence relation \(\mathcal S\) is strongly ergodic. So every ergodic \(\mathcal R\)-action is strongly ergodic, and we conclude that \(\mathcal R\) has property (T).
\end{proof}

\bibliographystyle{alphaurl}
\bibliography{bib}

\end{document}